\documentclass[a4paper,11pt,twoside,reqno]{amsart}
\usepackage{fullpage}
\usepackage[T1]{fontenc}
\usepackage{multicol}
\usepackage[dvipsnames]{xcolor}
\usepackage[colorlinks=true,linkcolor=MidnightBlue,citecolor=MidnightBlue,urlcolor=MidnightBlue,pagebackref=true]{hyperref}
\usepackage{enumerate,bm}
\usepackage{amsmath,amssymb,amsfonts,amsthm,mathtools}
\usepackage{booktabs,array}
\usepackage{etoolbox}
\hypersetup{pdftitle={Landau's problems through ratios of L-functions},pdfsubject={Common arithmetic kernels and four conditional prime-counting asymptotics}}
\numberwithin{equation}{section}
\theoremstyle{plain}
\newtheorem{theorem}{Theorem}[section]
\newtheorem{lemma}[theorem]{Lemma}
\newtheorem{proposition}[theorem]{Proposition}
\newtheorem{corollary}[theorem]{Corollary}
\newtheorem{hypothesis}[theorem]{Hypothesis}
\theoremstyle{definition}

\theoremstyle{remark}

\newtheorem{introhypothesis}{Hypothesis}[section]

\makeatletter
\patchcmd{\@starttoc}
  {\@startsection{}\@M\z@{\linespacing\@plus\linespacing}%
    {.5\linespacing}{\centering\contentsnamefont}{#2}}
  {}{} {\PackageError{toc}{Could not remove the contents heading}{}}
\patchcmd{\@starttoc}
  {\addvspace{32\p@\@plus14\p@}}
  {}{} {\PackageError{toc}{Could not remove trailing contents space}{}}
\makeatother

\renewcommand{\le}{\leqslant}
\renewcommand{\leq}{\leqslant}
\renewcommand{\ge}{\geqslant}
\renewcommand{\geq}{\geqslant}

\newcommand{\Z}{\mathbb Z}
\newcommand{\R}{\mathbb R}

\newcommand{\F}{\mathbb F}
\newcommand{\Q}{\mathbb Q}

\newcommand{\ind}{\mathbf 1}

\newcommand{\Disc}{\operatorname{Disc}}

\DeclareMathOperator{\supp}{supp}
\DeclareMathOperator{\Res}{Res}

\author{Haozhe Gou}
\address{School of Mathematics, Shandong University, Jinan 250100, China}
\email{haozhegou@gmail.com}

\author{Victor Y. Wang}
\address{Institute of Mathematics, Academia Sinica, Taipei 106319, Taiwan}
\email{vywang@as.edu.tw}

\author{Max Wenqiang Xu}
\address{Yau Mathematical Sciences Center, Tsinghua University, Beijing 100084, China}
\address{Beijing Institute of Mathematical Sciences and Applications (BIMSA), Beijing, China}
\email{maxxu1729@gmail.com}

\title{Landau's problems through ratios of $L$-functions}
\date{}
\begin{document}
\begin{abstract}
Very recently, Lee announced a proof of the Twin Prime Conjecture
assuming the Ratios Conjecture for the Riemann zeta function.
We obtain conditional asymptotic formulas with power-saving errors for each of Landau's four prime problems.
In the function field setting, using Sawin's estimates for factorization functions in squarefree progressions and a discriminant cohomology argument, we obtain unconditional Goldbach and prime-pair asymptotics over sufficiently large fixed odd fields, including prime fields.
Similar results may be possible for other problems in analytic number theory.
\end{abstract}
\maketitle
\thispagestyle{empty}
\setcounter{tocdepth}{1}
\begin{center}
  \contentsnamefont Contents
\end{center}
{\footnotesize
\begin{multicols}{2}
\tableofcontents
\end{multicols}
}

\section{Introduction}
In his invited address at the 1912 International Congress of Mathematicians in Cambridge, Edmund Landau surveyed the progress and outstanding challenges in the theory of prime numbers and the Riemann zeta function. He singled out four fundamental problems that he regarded as beyond the reach of the methods available at the time. These problems, now known as Landau's four problems, are the Goldbach Conjecture, the Twin Prime Conjecture, the existence of a prime between consecutive squares, and the infinitude of primes of the form \(n^2+1\). For further historical background on these problems and an account of subsequent progress, see for instance \cite{Pintz2009,HelfgottTernary,MaynardTakagi,TurnageButterbaugh,GrimmeltMerikoski} and the references within.

It has long been known that Montgomery's Pair Correlation Conjecture
can be derived from a uniform version of the Hardy--Littlewood Prime Pair Conjecture
\cite{Montgomery,GM,MontgomerySoundararajan}.
Very recently, Lee~\cite{Lee} announced a result in the other direction,
showing that a non-uniform version of the Hardy--Littlewood Prime Pair Conjecture
can be derived from the Ratios Conjecture
of Conrey--Farmer--Zirnbauer~\cite{CFZ}.
Other recent consequences of the Ratios Conjecture can be found in
\cite{JT,BGWSums,WX},
with the work of Wang--Xu \cite{WX} being in retrospect closely related to Lee's work.
In the present paper, we establish asymptotics
for all four Landau problems
under the Ratios Conjecture.\footnote{Strictly speaking we also assume GRH, but as in Lee's work, it is possible to formulate the Ratios Conjecture in a way that does not require GRH.}
In the function field setting,
we then establish new unconditional Goldbach and prime-pair results,
inspired by recent progress on the Ratios Conjecture for harmonic families
by Sawin \cite{Sawin}.
Recent progress on the Ratios Conjecture for geometric families
can be found in \cite{BDPW,MPPRW,WangRatios,LLOP},
and progress on negative moments in various families can be found in
\cite{BFK,FloreaSmallShifts,BFZeta,BFM,BuiFloreaHeycock}.

%Landau's four prime problems concern Goldbach representations, fixed prime differences, primes between consecutive squares, and prime values of $n^2+1$. 
Landau's four problems involve rather different arithmetic constraints. 
The purpose of this paper is to relate them through ratio statistics of three classical families of $L$-functions:
% three natural families of $L$-functions give explicit arithmetic kernels, and suitable tests of those kernels recover the four prime counts.
Dirichlet $L$-functions,
Riemann zeta shifts,
and Hecke $L$-functions on $\Q(i)$.
The ratios recipe of Conrey--Farmer--Zirnbauer~\cite{CFZ}, building on the moments recipe of \cite{CFKRS}, predicts averages of quotients of products of $L$-functions.
% For the continuous zeta family, Lee~\cite{Lee} already obtains the Hardy--Littlewood asymptotic for every fixed even prime difference from a restricted uniform ratios formula, including a version without RH. The connection between zero correlations and primes in short intervals is classical~\cite{GM,HB,LPZ}. We place these implications alongside Goldbach and quadratic prime values, identify the arithmetic relation between them, and specify the uniformity needed in each ratios input.

%We state the three ratios inputs together, with their uniformity and error terms, and then give a separate theorem for each prime problem. The quantitative statements use power-saving errors; the weaker qualitative variants are recorded in Appendix~\ref{wk:section}.
We state the three ratios hypotheses and their uniformity ranges before giving the four theorems. Asymptotic formulas under weaker ratios hypotheses are given in Appendix~\ref{wk:section}.

For $h\in\mathbb Z_{\ge1}$, define
\begin{equation}\label{eq:singular}
 \mathfrak S(h)=\prod_p\sigma_p(h),\qquad
 \sigma_p(h)=\begin{cases}p/(p-1),&p\mid h,\\p(p-2)/(p-1)^2,&p\nmid h.\end{cases}
\end{equation}
Thus $\mathfrak S(h)=0$ for odd $h$, and for even $h$ it is the positive Hardy--Littlewood singular series. For the quadratic problem, write $\chi=\chi_{-4}$ and put
\begin{equation}\label{eq:quadratic-series}
 \mathfrak S_E=\prod_{p>2}\left(1-\frac{\chi(p)}{p-1}\right)
 =\frac4\pi\prod_{p\equiv1\ (4)}\frac{p(p-2)}{(p-1)^2}
              \prod_{p\equiv3\ (4)}\frac{p^2}{p^2-1}.
\end{equation}
The first product is ordered by increasing primes; the second expression is absolutely convergent.

\subsection{Ratios Conjectures}
We state three Ratios Conjectures below, for primitive
Dirichlet $L$-functions, the Riemann zeta function, and angular Hecke
$L$-functions over $\mathbb Q(i)$. The untwisted Dirichlet sum runs over
all primitive characters modulo $q$. The Hecke average is unweighted.
Section~\ref{sec:prelim} gives the Euler products in the conjectured main
terms. Section~\ref{ci:section} defines $\mathcal M_q$, the sum of the two
main terms for the Dirichlet conjecture. At a removable singularity, we
evaluate the sum of the main terms by analytic continuation. On the
indicated disks, the errors are holomorphic in the numerator shifts.
Fix $\eta_*=1/16$ and $\rho=1/100$.

\begin{introhypothesis}[{Ratios of primitive Dirichlet $L$-functions ($\mathrm{RC}_D(\omega_D)$)}]
\label{hyp:int-D}
Fix $0<\omega_D\le 1/2$. There exist constants $\eta_0>0$ and $B\ge0$
such that for every fixed $0<\eta<\eta_0$ and every $\varepsilon>0$,
\[
 \sum_{\chi\bmod q}^{*}
 \frac{L(a,\chi)L(b,\overline\chi)}{L(s,\chi)L(t,\overline\chi)}
 =\mathcal M_q(a,b;s,t)
 +O_{\eta,\varepsilon}\bigl(q^{1-\omega_D+\varepsilon}(1+T)^B\bigr),
\]
uniformly for all conductors $q>1$, including composite conductors,
and all shifts satisfying
\[
 \Re a=\Re b=\tfrac12,\qquad
 \Re s=\Re t=\tfrac12+\eta,\qquad
 T=\max\{|\Im a|,|\Im b|,|\Im s|,|\Im t|\}\le 100\log(2q).
\]
The starred sum is over primitive characters modulo $q$.
The saving $\omega_D$ is independent of $\eta$.
\end{introhypothesis}
\begin{introhypothesis}[{Ratios of Zeta functions ($\mathrm{RC}_Z(\omega_Z)$)}]
\label{hyp:int-Z}
Let $0<\omega_Z\le1/2$, and put $s=1/2+a$ and $t_0=1/2+b$.
We use the $2/2$ main terms of Conrey, Farmer and
Zirnbauer~\cite[Conjecture~5.1 and Section~6.4]{CFZ}, denoted by
$\mathcal M_Z^{\mathrm{ex}}$ in~\eqref{co:Zmodel-exact}.
For every $\varepsilon>0$, assume
\[
 \int_0^Q\frac{\zeta(s+u+iy)\zeta(t_0+v-iy)}
                 {\zeta(s+iy)\zeta(t_0-iy)}\,dy
 =\mathcal M_Z^{\mathrm{ex}}(Q;a,b;u,v)
  +O_\varepsilon(Q^{1-\omega_Z+\varepsilon}),
\]
uniformly as $Q\to\infty$ for
\[
 \frac1{\log Q}\le\Re a,\Re b\le\eta_*,
 \qquad \max\{|\Im a|,|\Im b|\}\le Q^{1/2},
 \qquad |u|,|v|\le\min\{\rho,\tfrac12\Re a,\tfrac12\Re b\}.
\]
\end{introhypothesis}

\begin{introhypothesis}[Ratios of angular Hecke $L$-functions ($\mathrm{RC}_H(\omega_H)$)]\label{hyp:int-H}
Let $0<\omega_H\le1/2$. For $L_k(s)=L(s,\xi_k)$, with
$\xi_k((\alpha))=(\alpha/\overline\alpha)^{2k}$ on $\mathbb Z[i]$, assume
\[
 \frac1Q\sum_{1\le k\le Q}\frac{L_k(1/2+a)}{L_k(1/2+c)}
 =\mathcal M_H(Q;a,c)+O_\varepsilon(Q^{-\omega_H+\varepsilon}),
\]
uniformly for
\[
 \frac2{\log Q}\le\Re c\le\eta_*,\qquad
 |\Im c|\le Q^{3/4},\qquad |a-c|\le\rho,
\]
where $\mathcal M_H$ is the sum of the main terms given in Section~\ref{co:angular-input}.
\end{introhypothesis}

The errors in Hypotheses~\ref{hyp:int-Z} and~\ref{hyp:int-H} are
uniform in the imaginary parts of the shifts on the stated ranges.
This uniformity is used to remove smoothing without a loss in the exponent.
Under RH, Lemma~\ref{co:Z-reduction} gives the power form of the zeta
main terms on the range of Hypothesis~\ref{hyp:int-Z}.
For the angular family we use the main terms of
Waxman~\cite[Conjecture~3.3 and Section~3.6]{Waxman}.
The choices $\omega_Z=\omega_H=1/2$ give the square-root error bounds
for these two families.
Wang--Xu~\cite[Conjecture~2.1]{WX} conjecture a power saving for
prime conductors. Hypothesis~\ref{hyp:int-D} assumes the formula
for all primitive conductors and specifies the exponent $\omega_D$.
The case $\omega_D=1/2$ is the square-root error bound for this family.
Any fixed $\omega_D>0$ suffices for Theorem~\ref{thm:goldbach}.

%The last two error bounds are uniform in the imaginary parts of the shifts throughout the stated ranges. This uniformity is used to remove smoothing without a loss in the exponent. The upper bound $Q^{3/4}$ may be replaced by $Q^\theta$ for any fixed $1/2<\theta<1$.
%The uniform conjectures of Conrey--Snaith~\cite[Conjecture 2.1 and (2.11)]{CS} and Waxman~\cite[Conjecture 3.3 and Section 3.6]{Waxman} predict square-root error bounds, corresponding to $\omega_Z=\omega_H=1/2$. On the other hand, Wang--Xu~\cite[Conjecture 2.1]{WX} simply records a power saving for prime conductors. Hypothesis~\ref{hyp:int-D} assumes the formula for all primitive conductors and specifies the exponent $\omega_D$.
%The case $\omega_D=1/2$ is the square-root error bound for this family, but we emphasize that any fixed $\omega_D>0$ suffices for our Goldbach result below (Theorem~\ref{thm:goldbach}).

\subsection{Four conditional results over $\mathbb Z$}
Each theorem below assumes only its own ratios input and the indicated Riemann Hypothesis. The constants in the error terms may depend on the fixed ratios parameters.\footnote{For brevity, throughout the paper we refer to the
asymptotic formulas in Theorems~\ref{thm:goldbach}, \ref{thm:prime-pairs},
~\ref{thm:squares}, and~\ref{thm:quadratic-main} as ``Goldbach'', ``prime pairs'', ``Legendre'', and ``Prime values of $n^2+1$'', respectively. }

\begin{theorem}[Goldbach]\label{thm:goldbach}
Assume GRH for Dirichlet $L$-functions, and $\mathrm{RC}_D(\omega_D)$. For every $\varepsilon>0$ and every even $N$, we have
\begin{equation}\label{eq:main-G}
 \sum_{m+n=N}\Lambda(m)\Lambda(n)
 =N\mathfrak S(N)+O_\varepsilon(N^{1-\omega_D+\varepsilon}).
\end{equation}
%The same formula holds with the sum restricted to primes and weights $\log m\log n$. The ordered unweighted count satisfies
%\begin{equation}\label{eq:main-G-unweighted}
% \#\{(p_1,p_2):p_1+p_2=N\}
% =\mathfrak S(N)\int_2^{N-2}\frac{dt}{\log t\log(N-t)}
 %+O_\varepsilon(N^{1-\omega_D+\varepsilon}).
%\end{equation}
\end{theorem}

\begin{theorem}[Prime pairs]\label{thm:prime-pairs}
Assume RH for $\zeta$ and $\mathrm{RC}_Z(\omega_Z)$. For every fixed integer $h\ge1$ and every $\varepsilon>0$, we have
\begin{equation}\label{eq:main-T}
 \sum_{n\le X}\Lambda(n)\Lambda(n+h)
 =\mathfrak S(h)X+O_{h,\varepsilon}(X^{1-\omega_Z+\varepsilon}).
\end{equation}
%For even $h$, partial summation gives
%\[
% \#\{p\le X:p+h\text{ prime}\}
% =\mathfrak S(h)\int_2^X\frac{dt}{\log t\log(t+h)}
% +O_{h,\varepsilon}(X^{1-\omega_Z+\varepsilon}).
%\]
%Both errors are $O_{h,\varepsilon}(X^{1/2+\varepsilon})$ for the square-root input. The case $h=2$ gives twin primes.
\end{theorem}

\begin{theorem}[Legendre]\label{thm:squares}
Assume RH for $\zeta$ and $\mathrm{RC}_Z(\omega_Z)$. For every $\varepsilon>0$, we have 
\begin{align}
 \psi((n+1)^2)-\psi(n^2)
 &=2n+1+O_\varepsilon(n^{1-\omega_Z/3+\varepsilon}).
 \label{eq:main-Lpsi}
%  \pi((n+1)^2)-\pi(n^2)
% &=\int_{n^2}^{(n+1)^2}\frac{dt}{\log t}
%   +O_\varepsilon(n^{1-\omega_Z/3+\varepsilon}).\label{eq:main-L}
\end{align}
Moreover $p_{j+1}-p_j\ll_\varepsilon p_j^{(1-\omega_Z)/(2-\omega_Z)+\varepsilon}$.
For $\omega_Z=1/2$, the error is $O_\varepsilon(n^{5/6+\varepsilon})$ and the prime-gap exponent is $1/3+\varepsilon$.
\end{theorem}

\begin{theorem}[Prime values of $n^2+1$]\label{thm:quadratic-main}
Assume GRH for the angular Hecke family and $\mathrm{RC}_H(\omega_H)$. For every $\varepsilon>0$, we have
\begin{equation}\label{eq:main-E}
 \sum_{n\le X}\Lambda(n^2+1)
 =\mathfrak S_E X+O_\varepsilon(X^{1-\omega_H+\varepsilon}).
\end{equation}
%and
%\[
% \#\{n\le X:n^2+1\text{ prime}\}
% =\mathfrak S_E\int_1^X\frac{dt}{\log(t^2+1)}
% +O_\varepsilon(X^{1-\omega_H+\varepsilon}).
%\]
%Both errors are $O_\varepsilon(X^{1/2+\varepsilon})$ for the square-root input.
\end{theorem}

%The integral main terms retain the lower-order logarithmic terms, which are larger than a power-saving error. Their leading equivalents are $\mathfrak S(N)N/(\log N)^2$, $\mathfrak S(h)X/(\log X)^2$, $n/\log n$, and $\mathfrak S_E X/(2\log X)$, respectively.
The proofs are in Sections~\ref{ug:section}--\ref{sec:quadratic}. Appendix~\ref{wk:section} gives asymptotic formulas under Hypotheses~\ref{hyp:WZ} and~\ref{hyp:WH} together with the relevant Riemann Hypothesis. These ratios hypotheses allow errors of size $o(1)$ after normalization and shifts of logarithmic size.

\subsection{Unconditional function-field results}

Write $\mathcal M_j$ for the set of monic polynomials in $\F_q[T]$ of degree $j$, with $\mathcal M_0=\{1\}$.
All primes and divisors in $\F_q[T]$ are assumed to be monic. For $D\ne0$ define the singular series
\begin{equation}\label{ff:singular}
 \mathfrak S_q(D)=
 \prod_{P\mid D}(1-|P|^{-1})^{-1}
 \prod_{P\nmid D}\left(1-\frac1{(|P|-1)^2}\right),
 \qquad |P|=q^{\deg P},
\end{equation}
as in \cite[(6.10)]{SSAnnals}.
The product over $P\nmid D$ converges absolutely.

\begin{theorem}[Polynomial Goldbach]\label{ff:reflection}
Fix an odd prime power $q>10^{200}$. Uniformly over monic $N$ of degree $d\ge3$, we have
\begin{equation*}
 \sum_{F\in\mathcal M_{d-1}}\Lambda(F)\Lambda(N-F)
   =q^{d-1}\mathfrak S_q(N)+O_q(q^{23d/24}).
\end{equation*}
\end{theorem}

\begin{theorem}[Polynomial prime pairs]\label{ff:twins}
Fix an odd prime power $q>10^{200}$. For $n\ge1$, uniformly over $0\ne h\in\F_q[T]$ of degree less than $n$, we have
\begin{equation*}
 \sum_{F\in\mathcal M_n}\Lambda(F)\Lambda(F+h)
   =q^n\mathfrak S_q(h)+O_q(q^{23n/24}).
\end{equation*}
\end{theorem}

A proof sketch is given in Section~\ref{SKETCH-Fpt},
and full proofs are given in Section~\ref{ff:section}.
In view of the results of Sawin and Shusterman~\cite{SSAnnals},
our leading-order asymptotics are only new for $q=p,p^2$.
Our error terms may be new for other values of $q$ as well.

\subsection{Overview of the paper}
The proof separates the arithmetic information in a ratios prediction from the operation used to recover a prime count. Differentiating the ratios with respect to the numerator shifts gives logarithmic derivatives, whose Dirichlet series have von Mangoldt coefficients. The resulting main terms contain explicit Euler products. Character orthogonality, Fourier inversion, and estimates for sums over zeros then relate these main terms to the prime counts.

Table~\ref{tab:ratios} lists the ratios and the methods used in the proofs. Here $2/2$ denotes two numerator and two denominator factors, and $1/1$ denotes one of each.
\begin{table}[htbp]
\centering
\small
\begin{tabular}{@{}>{\raggedright\arraybackslash}p{.15\linewidth}>{\raggedright\arraybackslash}p{.29\linewidth}>{\raggedright\arraybackslash}p{.47\linewidth}@{}}
\toprule
Problem & {Ratios} & {Method}\\
\midrule
Goldbach & {Primitive Dirichlet $L$-functions, $2/2$} & {The parity twist and summation over primitive conductors}\\[3pt]
Prime pairs & {Riemann zeta function, $2/2$} & {Fourier inversion for a fixed difference}\\[3pt]
Legendre & {Riemann zeta function, $2/2$} & {Pair correlation and second moments of primes in short intervals}\\[3pt]
Quadratic values & {Angular Hecke $L$-functions, $1/1$} & {Poisson summation over the Hecke characters}\\
\bottomrule
\end{tabular}
\caption{{Ratios of $L$-functions and the methods used to estimate the prime counts.}}
\label{tab:ratios}
\end{table}

The arguments for Goldbach, prime pairs, and primes in short intervals use the functions $D_1$ and $C_1$ from Section~\ref{sec:kernels}. We define them here
\begin{align}\label{in:common-series}
 D_1(z)&=\sum_{h\ge1}\frac{\mathfrak S(h)}{h^z}\qquad(\Re z>1),\\
 C_1(z)&=\Xi(z)D_1(z)
 =\zeta(1+z)\zeta(1-z)
       \prod_p\left(1-\frac{(p^{-z}-1)^2}{(p-1)^2}\right),\notag
\end{align}
where $\Xi(z)=2(2\pi)^{-z}\Gamma(z)\cos(\pi z/2)$ is the archimedean factor from the zeta functional equation. A Fourier test selects the coefficient $\mathfrak S(h)$ for a fixed prime difference. The interval argument uses a positive test of the continued kernel $C_1$.

For Goldbach, summing over primitive conductors gives the exact projection (and this will be proved in Proposition~\ref{co:conductor})
\begin{equation}\label{in:projector}
 \frac1{\varphi(N)}\sum_{q\mid N}\varphi^*(q)q^{-z}D_q(z)
     =\sum_{N\mid h}\frac{\mathfrak S(h)}{h^z},\qquad \Re z>1.
\end{equation}
where $D_q(z)$ is defined in \eqref{eq:Dq}.
Here $\varphi^*(q)$ counts primitive characters, with $\varphi^*(1)=1$. Additive support then isolates $h=N$. Thus the prime-sum and prime-difference arguments recover coefficients of the same singular-series Dirichlet series by different tests.

The angular counterpart is $\mathcal D(z)=\sum_{b\ge1}S_b b^{-z}$, where $S_b$ is the Bateman--Horn singular series for $n^2+b^2$; see~\eqref{eqc:D-series}.
The following approximations illustrate the parallel
between the multiplicative scale of a fixed difference $h$
and a fixed vertical coordinate $i$:
\begin{equation*}
 X\log(1+h/n)=\frac h{n/X}+O_h(X^{-1}),\qquad
 X\arg(n+i)=\frac1{n/X}+O(X^{-2})\qquad(n\asymp X).
\end{equation*}
Both tests distinguish an integer coefficient at reciprocal spatial scale. The interval argument instead applies a positive test to the continued zeta kernel; it does not isolate a single singular-series coefficient.

For the Goldbach problem, character orthogonality introduces
$\chi(-1)$.
{Applying the functional equations to both numerator
factors gives a second untwisted ratios average, with numerator
arguments $(1-b,1-a)$. The two averages together determine the
difference between the even and odd contributions. The coefficients
in this linear combination can grow exponentially when the imaginary
parts of the numerator arguments are large and of opposite signs.
Rational product weights give sufficient Mellin decay to control
this growth and obtain asymptotics for the weighted prime sums.
To obtain a power-saving error for the sharp Goldbach count, we use
the fixed entire test}
\[
 w(v)=e^{20(1-v)}\left(\frac{\sin\pi(v-1)}{\pi(v-1)}\right)^2,
 \qquad w(1)=1,\quad w(k)=0\quad(k=2,3,\ldots).
\]
The weight $w((m+n)/N)$ selects the target sum exactly after character orthogonality. Its double Mellin transform contains a beta factor which cancels the exponential loss of the parity inverse on the required logarithmic window. Outside that window the actual signed moment is used directly. Section~\ref{ug:exact-section} proves the resulting power-saving formula without approximating a sharp cutoff by rational product weights. Integrating the same argument against bounded real power weights removes the prime logarithms with the same saving.

{Theorem~\ref{ci:extension} supplies the numerator disks from critical-line character data}. Its interpolation loss tends to zero with the fixed denominator displacement, so the final error is $N^{1-\omega_D+\varepsilon}$.
For the prime pairs problem and quadratic values problem, the uniform height-independent ratios error and bounded-variation cutoffs give Mellin $L^1$ norms of only logarithmic size. Proposition~\ref{qt:sharp} therefore removes smoothing with no fixed loss in the exponent. For intervals, the localized second moment in Section~\ref{si:sec} gives a quantitative monotonicity argument; optimizing the averaging width yields the saving $\omega_Z/3$ in the square parameter $n$.
Appendix~\ref{wk:section} retains the moving-line argument with $\eta_X=4/\log X$ for the weaker qualitative variants.

\subsection{Proof sketch of function-field results}
\label{SKETCH-Fpt}

In this section we sketch the proof of Theorems~\ref{ff:reflection} and~\ref{ff:twins}.
The full proofs are given in Section~\ref{ff:section}.

Fix an odd prime power $q>10^{200}$. We count pairs of monic irreducible polynomials in $\F_q[T]$ with prescribed sum or difference, keeping $q$ fixed as the degrees of the polynomials grow. For this purpose it suffices to estimate a sum with two von Mangoldt weights. The Type I/II divisor decomposition of Sawin and Shusterman~\cite[\S6]{SSAnnals} reduces the problem to a combination of (I) a level of distribution beyond $1/2$ for primes in squarefree progressions, and (II) progress on two-point M\"obius correlations.
For (I) we use work of Sawin \cite{Sawin},
and for (II) we use a sieve-theoretic argument inspired by recent work on the Twisted Chowla Conjecture \cite{ACD26}.

The main subtlety in (II) is that our $q$ is fixed,
whereas the paper \cite{ACD26} takes $q\to \infty$.
To handle our fixed $q$,
we need to prove good Betti number bounds
and to conduct a cohomological sieve in a fixed positive characteristic.
Our Betti number bounds build on a
complete-intersection estimate of Wan--Zhang~\cite[Theorem 1.2.2, arXiv version]{WanZhang}
(although earlier work of Katz \cite{Katz} would also suffice for large enough fixed $q$),
together with a low-degree parameterization of
spaces of square-free polynomials
inspired by Sawin and Shusterman~\cite[Appendix B, especially Lemma B.10]{SSBetti}.
On the other hand, in a linear range of homology degrees,
the general algebraic sieve framework of Das and Howe~\cite{DasHowe}
lets us degenerate our space of square-free polynomials (a discriminant complement)
to a space whose signed homology can be computed.

Essentially, \cite{DasHowe} lets us introduce extra averaging into our problem, while controlling the error term of the averaging process.
The question is then how to average.
The approach of \cite{ACD26} would be to degenerate to characteristic zero,
where topological methods of Aumonier \cite{Aumonier} would apply,
at least for constant coefficients.
We instead work purely in positive characteristic, and degenerate
our discriminant complement to the complement of a suitable resultant,
where certain derivatives are decoupled from certain polynomials,
allowing for signed homology to be computed,
and in fact shown to vanish,
in a linear range of degrees.

It may be possible to generalize our work
to obtain new qualitative and quantitative results on
M\"obius cancellation and square-free values
over fairly general univariate or multivariate polynomial sequences
in the function field setting, even when $q$ is prime.
Existing results on these problems can be found in work of
Sawin--Shusterman \cite{SSQuadratic} and Poonen \cite{Poonen}.
To keep our paper clean and readable, we do not pursue such generalizations here.

\subsection{Further remarks}
The common kernels give a set of sufficient hypotheses for the four integer asymptotics. They do not assert reverse implications between the prime conjectures, nor do they replace the three family-dependent ratios inputs by a single zeta formula. Appendix~\ref{rc:section} gives a further application to two-point M\"obius and Liouville correlations. That application uses a separate reciprocal-ratios hypothesis, with additional uniformity in small square twists.

\subsection*{Notation}
The von Mangoldt function is $\Lambda(n)$
and the second Chebyshev function is $\psi(x) = \sum_{n\le x} \Lambda(n)$.
Over $\F_q[T]$ we use $\Lambda(P^j)=\deg P$ for monic primes $P$ and $j\ge1$.
For $f\in C_c^\infty((0,\infty))$, the Mellin transform is
\begin{equation*}
 \widetilde f(s)=\int_0^\infty f(u)u^{s-1}\,du.
\end{equation*}
Contours $(c)$ run upward on $\Re s=c$. Fourier conventions are stated with each detector. Implied constants may depend on fixed test functions and hypothesis parameters. Uniformity in conductors, summation lengths, and degrees is stated explicitly.

\subsection*{Acknowledgments}

VYW was supported by NSTC Project Grant 114-2115-M-001-010-MY2.

\subsection*{AI statement}

After we saw Lee's preprint on the Twin Prime Conjecture \cite{Lee},
which uses a case of the Ratios Conjecture
analogous to that in the work of Harper and Wang--Xu on the M\"{o}bius Legendre Conjecture \cite{HarperBTSC,WX,HarperICM},
we wondered whether there might exist a unified conditional treatment of Landau's problems,
as well as the Power-saving Two-point Chowla Conjecture,
via the Ratios Conjecture.
Subsequently, in the function field setting,
we wondered whether recent progress on the Ratios Conjecture could lead to new unconditional progress on such problems.
The dynamical systems
GPT-6 Pro
and SAIR Open Math OS
were used to convert our rough ideas into rigorous proofs and understanding,
and to assist in the writing process,
building on the existing literature on
the Ratios Conjecture, Chowla's Conjecture, and \'etale cohomology.
In particular, in the function field setting,
we suggested (1) the use of Sawin's paper \cite{Sawin},
(2) the use of resolutions, compactifications, or boundary sieves,
(3) the role of Chowla's Conjecture as an intermediate target,
and (4) the use of double covers to estimate Betti numbers with signed coefficients.
It seems to us that either (2) or (3)
prompted the AI
to build on useful sieve-theoretic ideas from
recent work on the Twisted Chowla Conjecture \cite{ACD26}.
The paper \cite{ACD26} involves passage from characteristic zero
to positive characteristic.
However, to enhance the diversity of the literature,
we chose to present an argument working purely in positive characteristic.

\section{{Ratios Conjectures and arithmetic factors}}\label{sec:prelim}\label{sec:kernels}
Fix $0<\eta<1/4$ and write $\sigma=1/2+\eta$. For the proofs with fixed Mellin contours we take $0<\omega\le1/2$ and $\omega>2\eta$. The hypotheses in the introduction allow $\eta$ to be arbitrarily small. We also give weaker hypotheses with contours approaching the critical line in Appendix~\ref{wk:section}.

\subsection{{Euler products for the Dirichlet and zeta ratios}}
For $q\geq1$, let $\zeta_q(s)=\zeta(s)\prod_{p\mid q}(1-p^{-s})$, and define
\begin{align}
 A_q(a,b;c,d)&=\prod_{p\nmid q}
 \frac{(1-p^{-c-d})(1-p^{-a-d}-p^{-b-c}+p^{-c-d})}
 {(1-p^{-a-d})(1-p^{-b-c})},\label{eq:Akernel}\\
 Z_q(a,b;c,d)&=\frac{\zeta_q(a+b)\zeta_q(c+d)}{\zeta_q(a+d)\zeta_q(b+c)}
                  A_q(a,b;c,d).\label{eq:Zkernel}
\end{align}
{The variables $a,b,c,d$ denote arguments of the $L$-functions. The Euler product for $A_q$ converges absolutely and locally uniformly in the regions used below. For each $p\nmid q$, the local factor of $Z_q$ is the constant term in the Laurent expansion of}
\begin{equation*}
 \frac{(1-p^{-c}z)(1-p^{-d}z^{-1})}
 {(1-p^{-a}z)(1-p^{-b}z^{-1})}.
\end{equation*}
{In the region of absolute convergence this gives}
\begin{equation*}
 Z_q(a,b;c,d)=\prod_{p\nmid q}
 \frac{1-p^{-a-d}-p^{-b-c}+p^{-c-d}}{1-p^{-a-b}},
\end{equation*}
and~\eqref{eq:Zkernel} supplies the meromorphic continuation needed here.

For $\Re z>0$, put
\begin{align}
 F_p(z)&=\frac{1-2/p+p^{-1-z}}{(1-1/p)^2},\label{eq:Fp}\\
 D_q(z)&=\left(\frac q{\varphi(q)}\right)^2\zeta(z)
      \prod_{p\nmid q}F_p(z)\prod_{p\mid q}(1-p^{z-1}),\label{eq:Dq}\\
 \Xi(z)&=2(2\pi)^{-z}\Gamma(z)\cos(\pi z/2),\label{eq:Xi}\\
 C_q(z)&=\Xi(z)D_q(z)
 =\left(\frac q{\varphi(q)}\right)^2\zeta(1-z)
      \prod_{p\nmid q}F_p(z)\prod_{p\mid q}(1-p^{z-1}),\label{eq:Cq}\\
 P_q(z)&=\sum_{(n,q)=1}\frac{\Lambda(n)^2}{n^{1+z}}.\label{eq:Pq}
\end{align}
{The second expression in~\eqref{eq:Cq} follows from the functional equation of the Riemann zeta function.} Since $F_p(z)=1+O(p^{-1-\Re z}+p^{-2})$, the infinite product is absolutely convergent for $\Re z>0$. For $q=1$, $D_1$ has a simple pole of residue one at $z=1$ and no other pole in this half-plane. If $q>1$, the finite product cancels this pole.

To pass from ratios to von Mangoldt coefficients, {differentiate once with respect to each numerator shift.} For denominator arguments $s,t$ with $\Re s,\Re t>1/2$, {the derivative of the first main term is}
\begin{equation*}
 \left.\partial_u\partial_v Z_q(s+u,t+v;s,t)\right|_{u=v=0}
 =P_q(s+t-1).
\end{equation*}
 {The term from the two functional equations contains $Z_q(1-t-v,1-s-u;s,t)$, which has the factor}
\begin{equation*}
 \zeta_q(1-u)^{-1}\zeta_q(1-v)^{-1}
 =\left(\frac q{\varphi(q)}\right)^2uv+O(|uv|(|u|+|v|)).
\end{equation*}
 {The mixed derivative of $Z_q(1-t-v,1-s-u;s,t)$ at $u=v=0$ is therefore} $C_q(z)$, with $z=s+t-1$. Any derivative falling on an additional holomorphic factor leaves one of the two zeros uncanceled and contributes nothing. This is why the functional-equation factors remain undifferentiated in  {the formulas for the Dirichlet-character and zeta averages}.

\subsection{ {Singular series and Ramanujan sums}}
For $D_1$ as defined in~\eqref{eq:Dq}, the case $q=1$ of
Proposition~\ref{co:conductor} gives
\begin{equation}\label{co:series}
 D_1(z)=\sum_{h\geq1}\frac{\mathfrak S(h)}{h^z}
 \qquad(\Re z>1).
\end{equation} 
 {For the continuation near $z=0$, put $Z=\zeta'/\zeta$ and we claim that}
\begin{align}
 P_1(z)&=Z'(1+z)-\sum_p\left(\frac{\log p}{p^{1+z}-1}\right)^2,\label{co:Pcontinuation}\\
 C_1(z)&=\zeta(1+z)\zeta(1-z)\mathcal A(z),\qquad
 \mathcal A(z)=\prod_p\left(1-\frac{(p^{-z}-1)^2}{(p-1)^2}\right).
 \label{co:Ccontinuation}
\end{align}
The prime sum in~\eqref{co:Pcontinuation} and
the Euler product  $\mathcal A(z)$
in~\eqref{co:Ccontinuation} converge absolutely and uniformly
on compact subsets of $\Re z>-1/2$. The local identity
\begin{equation*}
 (1-p^{-1-z})F_p(z)=1-\frac{(p^{-z}-1)^2}{(p-1)^2}
\end{equation*}
proves~\eqref{co:Ccontinuation}. Expanding the two local geometric series proves~\eqref{co:Pcontinuation}; summing the local densities in~\eqref{eq:singular} proves~\eqref{co:series},  including the factor at two. Thus the prime-difference and zero-correlation calculations use the same pair $P_1,C_1$.

 {The coefficients of $D_q$ involve the Ramanujan sum $c_q(h)$ and the singular series with the factors at primes dividing $q$ omitted. For $q\geq1$, write}
\begin{equation*}
 \mathfrak S_q(h)=\prod_{p\nmid q}\sigma_p(h),\qquad
 c_q(h)=\sum_{d\mid(q,h)}d\mu(q/d),\qquad
 K_q(h)=\frac q{\varphi(q)^2}c_q(h)\mathfrak S_q(h).
\end{equation*}
We set $\varphi^*(1)=1$ and for $q>1$, $\varphi^*(q)$ counts primitive characters modulo $q$.

\begin{proposition}[{Dirichlet coefficients and the sum over conductors}]\label{co:conductor}
For $\Re z>1$,
\begin{equation}\label{co:Dseries}
 D_q(z)=\sum_{h\geq1}K_q(h)(h/q)^{-z},
\end{equation}
and for all positive integers $N,h$,
\begin{equation}\label{co:selection}
 \sum_{q\mid N}\varphi^*(q)K_q(h)
   =\varphi(N)\mathfrak S(h)\mathbf1_{N\mid h}.
\end{equation}
Consequently~\eqref{in:projector} holds.
\end{proposition}
\begin{proof}
First note the local identity
\begin{equation*}
 (1-p^{-z})\sum_{v\geq0}\sigma_p(p^v)p^{-vz}=F_p(z)
\end{equation*}
which also holds at $p=2$ with $F_2(z)=2^{1-z}$. Hence $\sum_{k\geq1}\mathfrak S_q(k)k^{-z}=\zeta(z)\prod_{p\nmid q}F_p(z)$. Expand the finite product in~\eqref{eq:Dq}, obtaining
\begin{equation*}
 D_q(z)=\left(\frac q{\varphi(q)}\right)^2
 \sum_{d\mid\operatorname{rad}(q)}\frac{\mu(d)}d
 \sum_{k\geq1}\mathfrak S_q(k)(k/d)^{-z}.
\end{equation*}
Set $h=qk/d$. The coefficient of $(h/q)^{-z}$ is
\begin{equation*}
 \left(\frac q{\varphi(q)}\right)^2\mathfrak S_q(h)
 \sum_{\substack{d\mid q\\q/d\mid h}}\frac{\mu(d)}d=K_q(h),
\end{equation*}
which proves~\eqref{co:Dseries}.

For even $h$, divide~\eqref{co:selection} by $\mathfrak S(h)$ and use multiplicativity. For $p^a\Vert N$, with $v=v_p(h)$, the local sum is
\begin{equation*}
 1+\sum_{j=1}^a\varphi^*(p^j)
       \frac{p^jc_{p^j}(h)}{\varphi(p^j)^2\sigma_p(h)}.
\end{equation*}
Here $c_{p^j}(h)$ is $\varphi(p^j)$ for $j\leq v$, $-p^{j-1}$ for $j=v+1$, and zero thereafter; also $\varphi^*(p)=p-2$ and $\varphi^*(p^j)=p^{j-2}(p-1)^2$ for $j\geq2$. If $v=0$, the first nonconstant summand is $-1$. If $v\geq1$, the sum through $\min(v,a)$ is $\varphi(p^{\min(v,a)})$, and the next summand cancels it when $v<a$. This gives $\varphi(p^a)\mathbf1_{p^a\mid h}$, also at two since $\varphi^*(2)=0$.

For odd $h$, terms with odd $q$ vanish at the local factor two. If $2\Vert q$, $\varphi^*(q)=0$; if $4\mid q$, then $c_q(h)=0$. Thus both sides vanish. Finally insert~\eqref{co:Dseries} and use~\eqref{co:selection} to obtain~\eqref{in:projector}.
\end{proof}

\subsection{ {Coefficient extraction by Mellin inversion}}
 {To recover the Dirichlet coefficients, we move the contour into the half-plane of absolute convergence. A test function of integral zero cancels a possible pole at $z=1$.}

\begin{lemma}\label{co:extraction}
Let $0<c<1<d$. Suppose $D$ is meromorphic in $c\leq\Re z\leq d$, has polynomial vertical growth, and has no pole there except possibly a simple pole at one. Suppose $D(z)=\sum_{j\geq1}a_jj^{-z}$ converges absolutely on $(d)$. For $f,h\in C_c^\infty((0,\infty))$ with $\int h=0$,
\begin{equation}\label{co:Mellin-extraction}
 \frac1{2\pi i}\int_{(c)}D(z)\widetilde f(z)\widetilde h(z)\,dz
 =\sum_{j\geq1}a_j\int_0^\infty f(u)h(j/u)\,\frac{du}{u}.
\end{equation}
If the product of the supports of $f$ and $h$
contains just one positive integer $j_0$, only $a_{j_0}$ contributes.
\end{lemma}
\begin{proof}
The equality $\widetilde h(1)=0$ cancels the possible pole.  {The rapid decay of the Mellin transforms permits a shift to $(d)$, where the Dirichlet series converges absolutely. Mellin inversion for the multiplicative convolution gives}~\eqref{co:Mellin-extraction}. The support assertion follows directly.
\end{proof}
For prime differences, $a_j=\mathfrak S(j)$.  {For the Hecke ratios, the coefficients are the singular series $S_j$ for $n^2+j^2$ in~\eqref{eqc:D-series}. In the Goldbach argument the identity~\eqref{in:projector} restricts the first sequence to multiples of $N$.}

\subsection{ {Ratios of the Riemann zeta function}}\label{co:weak-inputs}
{
Write $\zeta(s)=\chi_\zeta(s)\zeta(1-s)$, where
\[
 \chi_\zeta(s)=\pi^{s-1/2}
             \frac{\Gamma((1-s)/2)}{\Gamma(s/2)}.
\]
For $s=1/2+a$ and $t_0=1/2+b$, the two main terms in
Hypothesis~\ref{hyp:int-Z} are
\begin{align}\label{co:Zmodel-exact}
 \mathcal M_Z^{\mathrm{ex}}(Q;a,b;u,v)
 ={}&QZ_1(s+u,t_0+v;s,t_0)\\
 &+ Z_1(1-t_0-v,1-s-u;s,t_0)
    \int_0^Q\chi_\zeta(s+u+iy)\chi_\zeta(t_0+v-iy)\,dy.\notag
\end{align}\par}

{ Write $\mathcal M_Z$ for the main terms with the power approximation
to the product of the functional-equation factors,}
\begin{align}\label{co:Zmodel}
 \mathcal M_Z(Q;a,b;u,v)
 ={}&QZ_1(s+u,t_0+v;s,t_0)\\
 &+(2\pi)^{a+b+u+v}\frac{Q^{1-a-b-u-v}}{1-a-b-u-v}
       Z_1(1-t_0-v,1-s-u;s,t_0).\notag
\end{align}
{
On $a+b+u+v=0$, both prefactors of the reflected kernel equal $Q$.
For the exact integral this follows from $\chi_\zeta(s)\chi_\zeta(1-s)=1$.
The arithmetic kernels have opposite residues there, so each sum
of main terms is regular.

\begin{lemma}\label{co:Z-reduction}
Assume RH. Uniformly on the range in Hypothesis~\ref{hyp:int-Z},
for every $\varepsilon>0$,
\[
 \mathcal M_Z^{\mathrm{ex}}(Q;a,b;u,v)-\mathcal M_Z(Q;a,b;u,v)
 \ll_\varepsilon (2+|\Im a|+|\Im b|)Q^\varepsilon.
\]
Thus under RH one can use $\mathcal M_Z$ in place of
$\mathcal M_Z^{\mathrm{ex}}$ in Hypothesis~\ref{hyp:int-Z}
with the same error bound.
\end{lemma}

\begin{proof}
Put $w=a+b+u+v$ and $L=2+|\Im a|+|\Im b|$. The numerator disks
give $\Re(a+u),\Re(b+v)>0$ and $1/\log Q\le\Re w<1/4$.
On $0<y\le4L$, the exact product of functional-equation factors is
$O(1)$ by Stirling's formula and boundedness on compact sets.
The integral of $|(y/2\pi)^{-w}|$ on this range is also $O(L)$.

For $y\ge4L$ write
\[
 \frac{\chi_\zeta(s+u+iy)}{\chi_\zeta(1-t_0-v+iy)}
 =(y/2\pi)^{-w}R(y).
\]
The quotient on the left is the exact product in~\eqref{co:Zmodel-exact}.
Its modulus is $O(y^{-\Re w})$, so $R(y)\ll1$.
On the segment from $1-t_0-v+iy$ to $s+u+iy$, Stirling's formula
gives, uniformly for points $\xi$ on the segment,
\[
 \frac{\chi_\zeta'}{\chi_\zeta}(\xi)
 =-\log(y/2\pi)+O(L/y),\qquad
 i\left(\frac{\chi_\zeta'}{\chi_\zeta}\right)'(\xi)
 =-1/y+O(L/y^2).
\]
Integrating along the segment, whose displacement is $w$, gives
$R(y)-1\ll |w|L/y$ and $R'(y)\ll |w|L/y^2$.
For the first estimate, if $R(y)=e^E$ for the logarithm obtained
by this integration, use $|e^E-1|\le |E|\max(1,|R(y)|)$.
Integration by parts on $[4L,Q]$ with the primitive
$(2\pi)^w y^{1-w}/(1-w)$ bounds the integrated difference by
$O(L\log(2Q))$, since $|w|/|1-w|\ll1$.

In the reflected kernel the zeta quotient is
\[
 \frac{\zeta(1-w)\zeta(1+a+b)}{\zeta(1-u)\zeta(1-v)}.
\]
The Euler product in~\eqref{eq:Akernel} is uniformly bounded here,
as its local factors are $1+O(p^{-2+2\rho})$.
The two reciprocal zeta factors are bounded on the fixed disks
about one, and $\zeta(1+a+b)\ll\log Q$ by absolute convergence.
For $\zeta(1-w)$, use the Lindel\"of bound under RH when
$|\Im w|\ge1$ and the pole bound $O(1+|w|^{-1})$ otherwise.
Thus the reflected kernel is $O_\varepsilon(Q^\varepsilon)$.
Multiplying by the integrated difference proves the lemma.
Since $L\ll Q^{1/2}$ and $\omega_Z\le1/2$, the difference is absorbed
by the error in Hypothesis~\ref{hyp:int-Z}.
\end{proof}

The qualitative hypothesis below uses the power form on a logarithmic
height range and on fixed numerator disks.\par}

\begin{hypothesis}[ {Zeta ratios with an $o(Q)$ error}]\label{hyp:WZ}
For every fixed $A>0$, uniformly for
\begin{equation*}
 \frac1{\log Q}\leq\Re a,\Re b\leq\frac8{\log Q},\qquad
 |\Im a|,|\Im b|\leq A\log Q,\qquad |u|,|v|\leq\rho,
\end{equation*}
one has
\begin{equation}\label{co:weak-Z}
 \int_0^Q\frac{\zeta(s+u+iy)\zeta(t_0+v-iy)}
                   {\zeta(s+iy)\zeta(t_0-iy)}\,dy
       =\mathcal M_Z(Q;a,b;u,v)+o(Q).
\end{equation}
The two main terms are interpreted together at removable singularities.
\end{hypothesis}

{
The power form is the one in Conrey--Snaith~\cite[Conjecture~2.1]{CS}.
Under RH the sum of the main terms and the error in
Hypothesis~\ref{hyp:WZ} are holomorphic on the numerator disks.
The fixed radius allows Cauchy differentiation of the $o(Q)$ error
with a constant loss.\par}

\subsection{ {Ratios of angular Hecke $L$-functions}}\label{co:angular-input}
 {The Hecke characters below depend on the argument of a Gaussian integer. Summing over the character index with a suitable weight allows us to isolate a fixed imaginary part of a Gaussian integer.}

Let $F=\mathbb Q(i)$ and let $\chi=\chi_{-4}$. For a nonzero integral ideal $\mathfrak a=(\alpha)$ of $\mathbb Z[i]$, put
\begin{equation*}
 \xi_k(\mathfrak a)=\left(\frac{\alpha}{\overline\alpha}\right)^{2k},
 \qquad L_k(s)=L(s,\xi_k),\qquad k\in\mathbb Z.
\end{equation*}
The character is independent of the generator. Conjugation gives $L_{-k}=L_k$, and $L_k$ is entire for $k\ne0$. These are the angular Hecke $L$-functions considered by Waxman~\cite{Waxman}.

For complex $a,c$, define
\begin{align}
 A_3(a,c)&=\prod_{p\equiv3\ (4)}
 \frac{1-p^{-2-4c}}{1-p^{-2-2a-2c}},\label{eqc:A3}\\
 A_1(a,c)&=\prod_{p\equiv1\ (4)}
 \frac{(1-2p^{-1-a-c}+p^{-1-2c})(1-p^{-1-2c})}
 {(1-p^{-1-a-c})^2},\label{eqc:A1}\\
 G(a,c)&=\frac{\zeta(1+2a)L(1+2c,\chi)}
 {\zeta(1+a+c)L(1+a+c,\chi)}A_1(a,c)A_3(a,c).
 \label{eqc:G}
\end{align}
The products converge absolutely and locally uniformly  {in neighborhoods of $a=c$ and $a=-c$}. The factor at $2$ is included in the zeta quotient.

 {The two main terms in the conjecture for the average of ratios are}
\begin{equation*}
 \mathcal M_H(Q;a,c)=G(a,c)+\frac1{1-2a}
                           \left(\frac\pi{2Q}\right)^{2a}G(-a,c).
\end{equation*}
\begin{hypothesis}[{Hecke ratios with an $o(1)$ error}]\label{hyp:WH}
For some fixed $A>0$, uniformly for
\begin{equation*}
 \frac2{\log Q}\leq\Re c\leq\frac8{\log Q},\qquad
 |\Im c|\leq A\log Q,\qquad |a-c|\leq\rho,
\end{equation*}
one has
\begin{equation}\label{co:weak-H}
 \frac1Q\sum_{1\leq k\leq Q}\frac{L_k(1/2+a)}{L_k(1/2+c)}
       =\mathcal M_H(Q;a,c)+o(1).
\end{equation}
 {The sum of the main terms and the error are holomorphic in the numerator shift on the indicated disk.}
\end{hypothesis}

 {Waxman's conjecture~\cite[Conjecture 3.3]{Waxman} gives these main terms with an error of $O_\varepsilon(Q^{-1/2+\varepsilon})$ for the normalized average. We use the uniformity in the shifts discussed there. At $a=0$ the two main terms have opposite residues, since the factor multiplying $G(-a,c)$ equals one there. Under the relevant Riemann Hypotheses, the remaining reciprocal factors are holomorphic on the numerator disks. Cauchy's formula applies to the sum of the main terms and to the error.}

\begingroup
\begin{corollary}[{Consequences of square-root error bounds for ratios}]\label{co:standard}
Under the relevant Riemann Hypotheses and Hypotheses~\ref{hyp:int-D}--\ref{hyp:int-H} with $\omega_D=\omega_Z=\omega_H=1/2$, the Goldbach, fixed-difference, and quadratic counts, weighted or unweighted with the stated main terms, have errors $O_\varepsilon(N^{1/2+\varepsilon})$, $O_{h,\varepsilon}(X^{1/2+\varepsilon})$, and $O_\varepsilon(X^{1/2+\varepsilon})$, respectively. Moreover, \eqref{eq:main-Lpsi} holds with error $O_\varepsilon(n^{5/6+\varepsilon})$.
\end{corollary}
 For the Dirichlet-character ratios this corollary assumes the
square-root error at every primitive conductor. More generally,
Theorem~\ref{thm:goldbach} gives an error of
$O_\varepsilon(N^{1-\omega_D+\varepsilon})$ for any fixed positive
$\omega_D$ in Hypothesis~\ref{hyp:int-D}.
\endgroup

\subsection{ {Bounds in a zero-free half-plane}}
 {For Mellin contours with fixed real part} it suffices to assume a common $\beta<\sigma$ such that the relevant $L$-functions have no zeros in $\Re s>\beta$. RH or GRH gives $\beta=1/2$. A standard logarithmic-derivative estimate in a closed zero-free half-plane gives, for primitive Dirichlet characters,
\begin{equation}\label{eq:logbound}
 \frac{L'}L(\sigma+it,\chi)\ll_{\sigma-\beta}\log(2q(2+|t|)).
\end{equation}
The analogous estimates hold for $\zeta$ away from its pole and for the angular Hecke family with conductor parameter $k$. One obtains these bounds from the local zero expansion of $L'/L$: zeros at bounded distance contribute $O(\log(q(2+|t|)))$, and their denominators are bounded below by  {the fixed distance $\sigma-\beta$ from the zero-free boundary}; the remaining terms obey the same bound. We use also a fixed polynomial bound for $1/L$ on these lines, obtained by integrating the zero-free logarithm from a line of absolute convergence. These estimates justify  {shifting the Mellin contours to a line of fixed real part}.  {Appendix~\ref{wk:section} uses RH or GRH to move the contours towards the critical line. See~\cite[Chapter 5]{IK} for the standard bounds.}
% \textcolor{blue}{do we want to mention the recent QRH work?--like what $\beta$ we can choose here}

 {By using Cauchy's formula, we can bound derivatives of the holomorphic error terms with respect to the numerator shifts. A radius of $r/\log Q$ introduces a power of $\log Q$, whereas a fixed radius introduces only a constant. Moreover, any fixed power of $\log Q$ can be absorbed into $Q^\varepsilon$.}

\subsection{ {Smooth weighted averages}}\label{sec:weight}
 {We use integration by parts to pass from integrals over $[0,U]$ to integrals with a smooth weight $g(t/T)$. If $g$ vanishes to sufficiently high order at zero, the contribution from the range of small $U$ is negligible.} Let $z=a+b$, $\Re a=\Re b=\eta$, and $L=2+|\Im a|+|\Im b|$. Suppose
\begin{equation}\label{eq:primitive-moment}
 \int_0^U (t)\,dt
 =U A(z)+(2\pi)^z\frac{U^{1-z}}{1-z}B(z)
  +O_\varepsilon(U^{1-\omega+\varepsilon}L^B)
\end{equation}
whenever $L\leq cU^\theta$, for some fixed $\theta>0$. Suppose {the integrand and the displayed main terms have polynomial growth} outside this range. If $g$ is even and Schwartz, vanishes to sufficiently high order $2M$ at zero, and $\widetilde g(w)=\int_0^\infty g(v)v^{w-1}\,dv$, then
\begin{align}\label{eq:weighted-moment}
 \int_0^\infty g(t/T)V_{a,b}(t)\,dt
 ={}&T\widetilde g(1)A(z)+T(2\pi/T)^z\widetilde g(1-z)B(z)\\
 &+O_\varepsilon(T^{1-\omega+\varepsilon}L^{B'})
    +O(T^{-2M}L^{C_M}).\notag
\end{align}
Here $B'$ is fixed and $C_M$ depends on $M$. To prove this,  {subtract the displayed main terms from the integral over $[0,U]$ and integrate by parts against} $T^{-1}g'(U/T)$. Above $U_0=C\max(2,L^{1/\theta})$ use~\eqref{eq:primitive-moment}; below $U_0$ use polynomial growth and $g'(v)\ll v^{2M-1}$. The contribution below $U_0$ is bounded by $T^{-2M}$ times a fixed power of $L$. Boundary terms vanish. On integrating against rapidly decreasing Mellin transforms, the last term in~\eqref{eq:weighted-moment} is $O(T^{-2M})$  {with an implied constant depending on finitely many seminorms of the test functions.}

\begingroup
\begin{lemma}[Removal of smoothing]\label{lem:unsmooth}
 {Suppose a smoothed count has an error $O(X^{1-\delta})$, with implied constant $O(v^{-A})$ for smooth majorants and minorants with transition width $v$. Here $\delta>0$ and $A\ge0$ are fixed. The error for the corresponding unsmoothed count is}
\[
 O_\varepsilon(X^{1-\delta/(A+1)+\varepsilon}).
\]
This applies equally to $\Lambda(n)\Lambda(n+h)$, to Goldbach  {when such estimates with compactly supported weights are available}, and to $\Lambda(n^2+1)$.
\end{lemma}
\begin{proof}
On a fixed dyadic interval, the transition contains $O(vX+1)$ integers. Each summand is $O(\log^2 X)$ in the first two cases and $O(\log X)$ in the quadratic case. Smooth majorants and minorants thus give errors
$O(vX\log^2X+X^{1-\delta}v^{-A}+\log^2X)$; the Goldbach main-term transition also uses $\mathfrak S(N)\ll_\varepsilon N^\varepsilon$. Choose $v=X^{-\delta/(A+1)}$ and sum dyadic blocks when appropriate. In the applications we may decrease $\delta$ if necessary so that $\delta/(A+1)<1$.
\end{proof}
 {Proposition~\ref{qt:sharp} retains the exponent of the smoothed estimate by using an error bound uniform in the imaginary parts of the shifts. In the Goldbach proof we use a weight that isolates $m+n=N$ exactly.}
\endgroup

\section{Goldbach from character ratios}
\label{ug:section}
In this section, we give the proof of Theorem \ref{thm:goldbach}.
 The Goldbach sum can be studied by imposing the congruence $m\equiv-n\bmod N$ and then isolating the terms with $m+n=N$. For $(mn,N)=1$, by the  orthogonality of characters we have
\begin{equation}\label{ug:orthogonality}
 \mathbf1_{N\mid m+n}
 =\frac1{\varphi(N)}\sum_{\chi\bmod N}
     \chi(-1)\chi(m)\overline{\chi(n)}.
\end{equation}
Note that the factor $\chi(-1)$ in~\eqref{ug:orthogonality} is not present in our standard Ratios Conjecture. We therefore  use the functional equations first to recover the parity-weighted average from two untwisted ratios averages. Then the Goldbach singular series comes from by applying the Mellin inversion and summation over primitive conductors. The main difficulty here is that the coefficients in this linear combination may grow exponentially when the imaginary parts of the numerator arguments have opposite signs. We first use rational Mellin weights to obtain a qualitative asymptotic formula, and then introduce a weight that isolates $m+n=N$ exactly and preserves the power-saving error. Finally, we derive the required estimates for the numerator shifts from Hypothesis~\ref{hyp:int-D} by interpolation.

Fix $0<\eta<1/4$, $\sigma=1/2+\eta$, and $0<r<\min(\eta/4,(1/4-\eta)/4)$. Write
\begin{equation*}
 R_q(a,b;c,d)=\sum_{\chi\bmod q}^{*}
 \frac{L(a,\chi)L(b,\overline\chi)}
 {L(c,\chi)L(d,\overline\chi)}.
\end{equation*}
All ratios formulas assumed in this section are untwisted and include both parities. Let $M_e(q)$ count primitive characters of parity $\chi(-1)=(-1)^e$, put $\varphi^*(q)=M_0(q)+M_1(q)$, and write
\begin{equation*}
 X_{q,e}(s)=\left(\frac q\pi\right)^{1/2-s}
 \frac{\Gamma((1-s+e)/2)}{\Gamma((s+e)/2)}.
\end{equation*}
Set
\begin{equation*}
 H_q(a,b)=\sum_{e=0,1}M_e(q)X_{q,e}(a)X_{q,e}(b).
\end{equation*}
 {Thus $H_q$ is the sum of the products of the functional equation factors over all primitive characters of conductor $q$.}

\begin{lemma}[Untwisted Dirichlet ratios]\label{ug:ratios}\label{ci:extension}
Assume GRH for Dirichlet $L$-functions and Hypothesis~\ref{hyp:int-D}, with $\eta>0$ sufficiently small.
For some fixed $2\eta<\omega\le1/2$ and $B\ge0$, uniformly for all $q>1$, $\Re s=\Re t=\sigma$, $|\Im s|,|\Im t|\le10\log(2q)$, and $|u|,|v|\le r$, one has
\begin{align}\label{ug:ordinary}
 R_q(a,b;s,t)
 ={}&\varphi^*(q)Z_q(a,b;s,t)
       +H_q(a,b)Z_q(1-b,1-a;s,t)\\
 &+O_\varepsilon\bigl(q^{1-\omega+\varepsilon}
       (1+|\Im s|+|\Im t|)^B\bigr),\notag
\end{align}
for both $(a,b)=(s+u,t+v)$ and $(a,b)=(1-t-v,1-s-u)$. The error is holomorphic in $u,v$ and uniform on  {the indicated closed disks}.
\end{lemma}

These are two choices of the numerator arguments in the untwisted Ratios Conjecture. The sum includes primitive characters of both parities. The numerator arguments have real parts in $(1/4,3/4)$ and the denominator arguments have real parts greater than $1/2$. We assume the formula for composite as well as prime conductors. Wang--Xu~\cite{WX} formulate their conjecture for prime conductors. In the ratios recipe~\cite{CFZ}, the first main term comes from the diagonal and the second from applying the functional equation to both numerator $L$-functions. The sum of the products of the two functional equation factors is $H_q$.

Since the saving exponent tends to $\omega_D$ as $\eta\to0$, we may choose $\eta$ sufficiently small so that $\omega>2\eta$.

\subsection{ {The parity twist from untwisted ratios}}
 {Applying the two functional equations gives a second linear combination of the sums over even and odd characters. The root numbers cancel because the characters are conjugate. We then bound the coefficients obtained by solving the two equations.}

For $a=s+u$, $b=t+v$, write
\begin{equation*}
 R=R_q(a,b;s,t),\qquad
 R^\vee=R_q(1-b,1-a;s,t),\qquad
 h_e=X_{q,e}(1-b)X_{q,e}(1-a).
\end{equation*}
For a primitive character $\chi$ of parity $e$, the
functional equations for $\chi$ and $\overline\chi$ give
\[
 L(1-b,\chi)L(1-a,\overline\chi)
 =h_eL(a,\chi)L(b,\overline\chi),
\]
since the root numbers cancel. Let $R_e$ denote
the contribution to $R$ from characters of parity
$e\in\{0,1\}$. We then have
\[
 R=R_0+R_1,\qquad
 R^\vee=h_0R_0+h_1R_1,
\]
and hence
\begin{equation}\label{ug:parity-recovery}
 \sum_{\chi\bmod q}^{*}\chi(-1)
 \frac{L(a,\chi)L(b,\overline\chi)}
 {L(s,\chi)L(t,\overline\chi)}
 =\frac{2R^\vee-(h_0+h_1)R}{h_0-h_1}.
\end{equation}

Put $z=a+b-1$. The duplication and reflection
formulas give
\begin{equation}\label{ug:hdifference}
 h_0-h_1=-4q^z(2\pi)^{-z-1}
       \Gamma(a)\Gamma(b)\sin(\pi z/2),
\end{equation}
where we have used
$\cos(\pi(a+b)/2)=-\sin(\pi z/2)$.
In particular, $h_0-h_1$ does not vanish on
the indicated disks.

Using $X_{q,e}(w)X_{q,e}(1-w)=1$, the two
instances of~\eqref{ug:ordinary} give the
following main terms for $R$ and $R^\vee$:
\[
\begin{aligned}
R:\quad&
 \varphi^*(q)Z_q(a,b;s,t)
 +\left(\frac{M_0(q)}{h_0}
       +\frac{M_1(q)}{h_1}\right)
       Z_q(1-b,1-a;s,t),\\
R^\vee:\quad&
 (M_0(q)h_0+M_1(q)h_1)Z_q(a,b;s,t)
 +\varphi^*(q)Z_q(1-b,1-a;s,t).
\end{aligned}
\]

The coefficients in~\eqref{ug:parity-recovery}
may grow exponentially in the imaginary parts
of the numerator arguments. Indeed,
\[
 \frac{h_0+h_1}{h_0-h_1}
 =-\frac{\cos(\pi(a-b)/2)}
         {\sin(\pi(a+b-1)/2)}.
\]
For $a=\sigma+iT$ and $b=\sigma-iT$, this becomes
\(
 -{\cosh(\pi T)}/{\sin(\pi\eta)},
\)
which grows exponentially as $T\to\infty$.
We therefore estimate the errors after
multiplication by Mellin transforms of the weights.

\subsection{Rational weights}
 {The Mellin transforms of smooth compactly supported functions decay faster than any power, which is insufficient to compensate for the exponential growth of the coefficients in~\eqref{ug:parity-recovery}. We use rational functions whose Mellin transforms contain two gamma factors. One compensates for that growth and the other gives exponential decay.}

Fix an integer $A>3$ and consider finite linear combinations of
\begin{equation}\label{ug:rational-basis}
 f_j(u)=\frac{u^j}{(1+u)^{A+j}},\qquad j=0,1,2,\ldots.
\end{equation}
Their Mellin transforms are
\begin{equation}\label{ug:rational-Mellin}
 \widetilde f_j(s)
 =\frac{\Gamma(s+j)\Gamma(A-s)}{\Gamma(A+j)}
 =\frac{\Gamma(s)(s)_j\Gamma(A-s)}{\Gamma(A+j)},
 \qquad -j<\Re s<A.
\end{equation}
Here $(s)_j$ is the rising factorial. The factor $\Gamma(A-s)$ supplies the additional decay.

Let $E$ and $E^\vee$ denote the error terms
obtained by applying~\eqref{ug:ordinary} to
$R$ and $R^\vee$, respectively.
Substituting these formulas into
\eqref{ug:parity-recovery}, we obtain the
main terms in~\eqref{ug:signed-model},
together with an error term
\[
 \frac{2E^\vee-(h_0+h_1)E}{h_0-h_1}.
\]

For rational weights $f,g$ from~\eqref{ug:rational-basis},
Stirling's formula and~\eqref{ug:hdifference} give,
uniformly on the numerator disks,
\begin{equation}\label{ug:stable-testing}
 \left|\frac{\widetilde f(s)\widetilde g(t)}
 {h_0-h_1}\right|
 +\left|\frac{(h_0+h_1)
 \widetilde f(s)\widetilde g(t)}{h_0-h_1}\right|
 \ll_{f,g,\eta,r}
 (1+|\Im s|+|\Im t|)^C
 e^{-c(|\Im s|+|\Im t|)}
\end{equation}
for some fixed $c>0$, uniformly in $q$.

To see this, write $x=\Im s$, $y=\Im t$ and
$z=a+b-1$, where $a=s+u$ and $b=t+v$.
The gamma ratios $\Gamma(s+j)/\Gamma(a)$
and $\Gamma(t+k)/\Gamma(b)$ have polynomial
growth, since $u,v$ are bounded. Moreover,
\[
 |\Gamma(A-s)\Gamma(A-t)|
 \ll (1+|x|+|y|)^C
 e^{-\pi(|x|+|y|)/2},
\]
and
\[
 \left|\sin\frac{\pi z}{2}\right|^{-1}
 \ll_{\eta,r}e^{-\pi|x+y|/2}.
\]
Since $q^{-\Re z}\le1$, these estimates and
\eqref{ug:hdifference} give
\[
 \left|\frac{\widetilde f(s)\widetilde g(t)}
 {h_0-h_1}\right|
 \ll (1+|x|+|y|)^C
 e^{-\pi(|x|+|y|)/2}.
\]
For the second term, use
\[
 \frac{h_0+h_1}{h_0-h_1}
 =-\frac{\cos(\pi(a-b)/2)}
         {\sin(\pi z/2)}.
\]
It follows that
\[
 \left|\frac{h_0+h_1}{h_0-h_1}\right|
 \ll_{\eta,r}
 e^{\frac{\pi}{2}(|x-y|-|x+y|)}.
\]
Together with the decay of the two Mellin
transforms and $|x-y|\le |x|+|y|$,
this proves~\eqref{ug:stable-testing}.

We apply Cauchy's formula in $u,v$ before
the two Mellin integrations. Since the disk
radii are fixed, \eqref{ug:stable-testing}
makes the resulting error integrable, giving
\begin{equation}\label{ug:tested-error}
 O_{f,g,\varepsilon}(q^{1-\omega+\varepsilon})
\end{equation}
for each conductor $q$.

Here the ratios formula is applied on
$|\Im s|,|\Im t|\le5\log(2q)$.
Outside this range we estimate the actual
signed logarithmic-derivative sum and its
main terms separately. Under GRH these
are bounded by
\[
 O_\eta\bigl(q^{1+\varepsilon}
 (1+|\Im s|+|\Im t|)^C\bigr).
\]
The exponential decay of the Mellin
transforms makes the remaining integral
$O_{f,g,\varepsilon}(q^{1-5\pi+\varepsilon})$.
After summing over conductors, dividing
by $\varphi(N)$, and including the factor
$N^{s+t}$ from Mellin inversion, this
contributes
$O_{f,g,\varepsilon}(N^{2\eta+\varepsilon})$.

\subsection{ {The weighted sum over multiples of the target integer}}

\begingroup
We first compute the main term of the parity-weighted ratios and then
apply it to the sum over all multiples of $N$. The calculation will
also apply to the weight in Section~\ref{ug:exact-section}.

Write $\Delta(q)=M_0(q)-M_1(q)$. Substituting the two main terms for $R$
and $R^\vee$ into~\eqref{ug:parity-recovery}, and using
$h_e^{-1}=X_{q,e}(a)X_{q,e}(b)$, gives the parity-weighted main term
\begin{align}\label{ug:signed-model}
 &\Delta(q)Z_q(a,b;s,t)
 +
 \bigl(M_0(q)X_{q,0}(a)X_{q,0}(b)
       -M_1(q)X_{q,1}(a)X_{q,1}(b)\bigr)Z_q(1-b,1-a;s,t).
\end{align}
Differentiate once in each of $u$ and $v$ at $u=v=0$. By the kernel
identities, %\textcolor{red}{cite the exact equations} in Section~\ref{sec:kernels}, 
the resulting main term is
\endgroup
\begin{equation}\label{ug:LDmodel}
 \Delta(q)P_q(z)+\big(M_0(q)X_{q,0}(s)X_{q,0}(t)
            -M_1(q)X_{q,1}(s)X_{q,1}(t)\big)C_q(z),
 \quad z=s+t-1.
\end{equation}
\begingroup
The two reciprocal zeta factors in the {reflected kernel} vanish in the
respective shifts. Hence the mixed derivative of that term acts only on
the kernel.
M\"obius inversion of character orthogonality gives
\[
 \Delta(q)=\sum_{d\mid(q,2)}\mu(q/d)\varphi(d),\qquad |\Delta(q)|\le2.
\]
In the coefficient of $C_q(z)$ in~\eqref{ug:LDmodel}, the sum of the
two gamma products has coefficient $\Delta(q)/2$, and their difference
has coefficient $\varphi^*(q)/2$. On $\Re s=\Re t=\sigma$,
$P_q(s+t-1)\ll_\eta 1$, while $C_q(s+t-1)$ is
$O_{\eta,\varepsilon}\bigl(q^\varepsilon(1+|\Im s|+|\Im t|)^C\bigr)$
by~\eqref{eq:Cq}. Stirling's formula bounds each gamma product by
$q^{-2\eta}(1+|\Im s|+|\Im t|)^C$ up to a constant. Thus both terms
proportional to $\Delta(q)$ have size
$O_{\eta,\varepsilon}\bigl(q^\varepsilon(1+|\Im s|+|\Im t|)^C\bigr)$.
\endgroup

Duplication and reflection give the identity
\begin{equation}\label{ug:beta}
 \frac{X_{q,0}(s)X_{q,0}(t)-X_{q,1}(s)X_{q,1}(t)}2\,\Xi(z)
   =q^{-z}\frac{\Gamma(1-s)\Gamma(1-t)}{\Gamma(1-z)}.
\end{equation}
Indeed
$X_{q,e}(s)=2q^{1/2-s}(2\pi)^{s-1}
 \Gamma(1-s)\sin(\pi(s+e)/2)$.
The difference of the sine products is $\sin(\pi z/2)$, and
$\Gamma(z)\sin(\pi z)/\pi=1/\Gamma(1-z)$.

\begingroup
Let $\mathcal W(s,t)$ be the double Mellin transform of a continuous
weight $W(u,v)$ on $(0,\infty)^2$, and put
\[
 F(x)=\int_0^x W(u,x-u)\,du .
\]
We use weights satisfying, for some $0\le\lambda,\mu\le1/2$,
\[
 |W(u,v)|\ll u^{-\lambda}v^{-\mu}(1+u)^{-A}(1+v)^{-A}.
\]
Suppose $\mathcal W$ has rapid vertical decay on closed substrips of
$\lambda<\Re s\le5/4$ and $\mu<\Re t\le5/4$, and $\widetilde F$
has rapid decay on the closed strip from $\min(2\eta,1/4)$ to $3/2$.
Suppose also that the integral of the absolute error in the
parity-weighted logarithmic-derivative formula with main
term~\eqref{ug:LDmodel}, against
$|\mathcal W(s,t)|\,|ds\,dt|$ on $(\sigma)\times(\sigma)$ is
$O_\varepsilon(q^{1-\omega+\varepsilon})$, as in~\eqref{ug:tested-error}.
Then for even $N$,
\begin{equation}\label{ug:weighted-multiples}
 \begin{aligned}
 N^{-\lambda-\mu}
 \sum_{\substack{m,n\ge1\\N\mid m+n}}
   \Lambda(m)\Lambda(n)W(m/N,n/N)
 ={}&N^{1-\lambda-\mu}\sum_{k\ge1}\mathfrak S(kN)F(k)
 +O_\varepsilon(N^{1-\omega+2\eta+\varepsilon}).
 \end{aligned}
\end{equation}
The implied constant depends on the bounds for the weight.
\endgroup

\begin{proof}[Proof of~\eqref{ug:weighted-multiples}]
\begingroup
Apply~\eqref{ug:orthogonality} and group the characters modulo $N$ by
their primitive conductors. A pair with $(mn,N)>1$ can contribute
only when $m$ and $n$ are powers of the same prime $p\mid N$.
In the normalized sum these pairs contribute $O(\log^3(2N))$, since
$m^{-\lambda}n^{-\mu}\le1$ and
\[
 \sum_{a\ge1}(1+p^a/N)^{-A}\ll_A 1+\frac{\log N}{\log p}.
\]

For a primitive character $\chi$ of conductor $q$, the logarithmic
derivative of the character it induces modulo $N$ differs from
$L'/L(s,\chi)$ by
$\sum_{p\mid N,\,p\nmid q}\chi(p)\log p/(p^s-\chi(p))$.
In double Mellin inversion, begin on $\Re s=\Re t=5/4$
and expand the product of the two logarithmic derivatives.
In terms containing an omitted Euler factor,
move its line to $\Re s=\lambda+\varepsilon/4$ or
$\Re t=\mu+\varepsilon/4$, taking $\varepsilon>0$ small. The Euler
factor sum on either line is $O_\varepsilon(\log(2N))$. Move any
remaining primitive logarithmic derivative to $\Re s=\sigma$ or
$\Re t=\sigma$ and apply the GRH bound~\eqref{eq:logbound}. After
multiplication by $N^{-\lambda-\mu}$ and averaging over characters,
these integrals contribute $O_\varepsilon(N^{\sigma+\varepsilon})$.
If a pole at $1$ is crossed, it comes from the principal character
alone. Its residue in an omitted-factor term contributes
$O_\varepsilon(N^\varepsilon)$ after division by $\varphi(N)$,
using $N/\varphi(N)\ll_\varepsilon N^\varepsilon$.

For the principal character, shifting the two primitive
logarithmic derivatives to $(\sigma)$ gives a double residue
$N^{2-\lambda-\mu}\mathcal W(1,1)/\varphi(N)$. The integrals with one
residue or with both variables on $(\sigma)$ contribute
$O_\varepsilon(N^{\sigma+\varepsilon})$ after division by
$\varphi(N)$ and normalization.
\endgroup

\begingroup
For $q>1$, use~\eqref{ug:LDmodel} and~\eqref{ug:beta} and set
$z=s+t-1$. The inner integral is
\begin{equation}\label{ug:Barnes}
 \frac{1}{2\pi i}\int_{(\sigma)}
 \mathcal W(s,1+z-s)
 \frac{\Gamma(1-s)\Gamma(s-z)}{\Gamma(1-z)}\,ds
 =\widetilde F(z),\qquad \Re z=2\eta .
\end{equation}
Indeed, the Mellin--Barnes formula writes $(u+v)^{z-1}$ as
the integral of the gamma quotient in~\eqref{ug:Barnes} times
$u^{s-1}v^{z-s}$. Integrating against $W(u,v)\,du\,dv$ proves
the identity. Since $\Re z<\sigma<1$, the gamma factors and the
bound on $W$ make the interchange absolute.

The two terms proportional to $\Delta(q)$ contribute
$O_\varepsilon(N^{2\eta+\varepsilon})$ in total. The integrated ratios
errors contribute $O_\varepsilon(N^{1-\omega+2\eta+\varepsilon})$,
since $\sum_{q\mid N}q^{1-\omega+\varepsilon}
\ll_\varepsilon N^{1-\omega+2\varepsilon}$ and the Mellin factor has
modulus $N^{2\sigma}$ on the two lines. We have used
$\sigma\le1-\omega+2\eta$ to absorb the character-replacement and
principal-character errors. Thus the normalized sum
in~\eqref{ug:weighted-multiples} equals
\[
 \frac{N^{2-\lambda-\mu}}{\varphi(N)}\mathcal W(1,1)
 +\frac{N^{1-\lambda-\mu}}{\varphi(N)}
  \sum_{\substack{q\mid N\\q>1}}\varphi^*(q)
  \frac{1}{2\pi i}\int_{(2\eta)}(N/q)^zD_q(z)\widetilde F(z)\,dz
 +O_\varepsilon(N^{1-\omega+2\eta+\varepsilon}).
\]
For $q>1$, the finite Euler factors cancel the pole of $D_q$ at $1$,
so we may move the contours to $\Re z=3/2$. On this line, write
the sum as the full sum over $q\mid N$ minus the $q=1$ term.
The full sum
is $N^{1-\lambda-\mu}\sum_{k\ge1}\mathfrak S(kN)F(k)$ by
\eqref{in:projector} and Mellin inversion. In the subtracted $q=1$
integral, move the contour to $\Re z=1/4$. The pole of $D_1$ at
$1$ has residue one, and $\widetilde F(1)=\mathcal W(1,1)$.
Its residue therefore cancels the principal-character term, and the
new integral contributes $O_\varepsilon(N^{1/4+\varepsilon})$.
\endgroup
\end{proof}

\begingroup
For the rational weights in~\eqref{ug:rational-basis}, apply
\eqref{ug:weighted-multiples} with $W(u,v)=f(u)g(v)$ and
$\lambda=\mu=0$. Its double Mellin transform is
$\widetilde f(s)\widetilde g(t)$, and the required integrated error
is~\eqref{ug:tested-error} together with the tail estimate that follows
it. The additive convolution
\[
 (f*g)(x)=\int_0^x f(u)g(x-u)\,du
\]
is $O(x)$ at zero and $O((1+x)^{-A})$ at infinity.
The same bounds hold after repeated application of $x\,d/dx$, by writing
the convolution as $x\int_0^1 f(xu)g(x(1-u))\,du$ and splitting at
$u=1/2$. Integration by parts gives rapid vertical decay of
$\widetilde{f*g}$ on closed substrips of $-1<\Re z<A$. Consequently,
with
\[
 \mathcal B_N(f,g)=\sum_{\substack{m,n\ge1\\N\mid m+n}}
       \Lambda(m)\Lambda(n)f(m/N)g(n/N),
\]
we obtain
\begin{equation}\label{ug:all-multiples}
 \mathcal B_N(f,g)
 =N\sum_{k\ge1}\mathfrak S(kN)(f*g)(k)
       +O_{f,g,\varepsilon}(N^{1-\omega+2\eta+\varepsilon}).
\end{equation}
\endgroup

\subsection{ {Approximation by rational weights}}
To pass to compactly supported weights, take
$w_A(u)=(1+u)^{-A}$. For even $N$ we have
\begin{equation}\label{ug:singular-domination}
 \mathfrak S(kN)\le\mathfrak S(N)\tau(k),\qquad k\ge1.
\end{equation}
Indeed each new odd prime divisor of $k$ contributes at most two. Splitting the convolution integral at $k/2$ gives
\begin{equation*}
 (w_A*w_A)(k)\ll_A(1+k)^{-A}.
\end{equation*}
Since $A>3$, the series in~\eqref{ug:all-multiples} for $w_A,w_A$
is $O_A(\mathfrak S(N))$. The error is $o(N)$ for sufficiently
small fixed $\varepsilon$, as $\omega>2\eta$. Using
$\mathfrak S(N)\gg1$ for even $N$, we obtain
\begin{equation}\label{ug:dominating-mass}
 \mathcal B_N(w_A,w_A)\ll_A N\mathfrak S(N).
\end{equation}

For every continuous compactly supported function $f_0$ on $[0,\infty)$ and every $\epsilon>0$, there is  {a rational function $f$} of the above form such that
\begin{equation}\label{ug:weighted-approx}
 |f(u)-f_0(u)|\le\epsilon w_A(u),\qquad u\ge0.
\end{equation}
Apply the Weierstrass approximation theorem to
$(1-x)^{-A}f_0(x/(1-x))$ on $[0,1]$, extending it by zero at $x=1$.
A polynomial $P$ approximating this function within $\epsilon$
gives $f(u)=(1+u)^{-A}P(u/(1+u))$.

Apply~\eqref{ug:weighted-approx} to both $f_0$ and $g_0$.
For $0<\epsilon\le1$, the difference between the two product
weights is bounded pointwise by
$O_{f_0,g_0,A}(\epsilon w_A(u)w_A(v))$.
By~\eqref{ug:dominating-mass} this changes the count by
$O_{f_0,g_0,A}(\epsilon N\mathfrak S(N))$.
The same bound for the change in the main term follows from
\eqref{ug:singular-domination} and the convolution bound for $w_A$.
Apply~\eqref{ug:all-multiples} to the fixed rational approximants,
then let $N\to\infty$ and $\epsilon\to0$ to obtain
\begin{equation}\label{ug:compact-tests}
 \mathcal B_N(f_0,g_0)
 =N\sum_{k\ge1}\mathfrak S(kN)(f_0*g_0)(k)
       +o_{f_0,g_0}(N\mathfrak S(N)).
\end{equation}
If both functions are supported in $(0,1)$, the count is supported
on $m+n=N$ and the main term is
$N\mathfrak S(N)\int_0^1 f_0(u)g_0(1-u)\,du$.

\begingroup
To recover the full von Mangoldt sum, fix $0<\delta<1/4$ and take
continuous cutoffs with values in $[0,1]$. Let $f_-$ be supported
in $[\delta,1-\delta]$ and equal to one on
$[2\delta,1-2\delta]$. Let $f_+$ be equal to one on $[0,1]$ and
supported in $[0,1+\delta]$. Then
\[
 \mathcal B_N(f_-,f_-)
 \le \sum_{m+n=N}\Lambda(m)\Lambda(n)
 \le \mathcal B_N(f_+,f_+).
\]
At $k=1$ the two convolutions lie in $[1-4\delta,1]$
and are equal to one, respectively. The upper convolution
at $k=2$ is $O(\delta)$ and vanishes at all $k\ge3$.
Since $\mathfrak S(2N)=\mathfrak S(N)$,
\eqref{ug:compact-tests} bounds the normalized count between
$1-O(\delta)+o_\delta(1)$ and $1+O(\delta)+o_\delta(1)$.
Letting $\delta\to0$ gives the qualitative asymptotic.
\endgroup

\subsection{ {A power-saving error for the Goldbach count}}\label{ug:exact-section}

\begingroup
A weight depending on the sum of the two variables can select
$m+n=N$ exactly. Put
\begin{equation}\label{ug:exact-weight}
 w(v)=e^{20(1-v)}\left(\frac{\sin\pi(v-1)}{\pi(v-1)}\right)^2,
 \qquad W(u,v)=u^{-\lambda}v^{-\mu}w(u+v),
 \quad 0\le\lambda,\mu\le\omega.
\end{equation}
With the removable value $w(1)=1$, the function vanishes at
all integers $k\ge2$. Therefore
\begin{equation}\label{ug:exact-selection}
 N^{-\lambda-\mu}
 \sum_{\substack{m,n\ge1\\N\mid m+n}}
 \Lambda(m)\Lambda(n)W(m/N,n/N)
 =\sum_{m+n=N}\frac{\Lambda(m)\Lambda(n)}{m^\lambda n^\mu}.
\end{equation}

The function $w$ is entire, is $O(|v|^2)$ at zero and decays
exponentially at infinity in the closed sector $|\arg v|\le1$,
since $20\cos1>2\pi\sin1$. Rotating the Mellin contour to
$\arg v=\operatorname{sgn}(y)$ for $y\ne0$ gives
\begin{equation}\label{ug:w-decay}
 |\widetilde w(b+iy)|\ll e^{-|y|}
\end{equation}
uniformly for $b$ in compact subsets of $(-2,\infty)$.
The double Mellin transform of $W$ is
\begin{equation}\label{ug:double-Mellin}
 \mathcal W(s,t)=
 \frac{\Gamma(s-\lambda)\Gamma(t-\mu)}
      {\Gamma(s+t-\lambda-\mu)}
 \widetilde w(s+t-\lambda-\mu),
 \qquad \Re s>\lambda,\quad \Re t>\mu .
\end{equation}
Indeed, integrating first over the fraction of $u+v$ contributed
by $u$ gives the beta factor. By Stirling's formula and
\eqref{ug:w-decay},
\[
 |\mathcal W(s,t)|
 \ll (1+|\Im s|+|\Im t|)^C e^{-|\Im s|-|\Im t|}
\]
on the closed substrips used in~\eqref{ug:weighted-multiples}.
Here the exponential from the gamma quotient is
$\exp\{-\tfrac{\pi}{2}(|\Im s|+|\Im t|-|\Im s+\Im t|)\}$.
The constants are uniform for $0\le\lambda,\mu\le\omega$ on lines
a fixed distance from $\Re s=\lambda$ and $\Re t=\mu$.
The exponential decay of $w$ on the positive real axis gives the
pointwise bound on $W$ required in~\eqref{ug:weighted-multiples}.

For $a=s+u$, $b=t+v$ on the numerator disks of
Lemma~\ref{ug:ratios}, \eqref{ug:hdifference} and the formula
for $(h_0+h_1)/(h_0-h_1)$ give
\begin{equation}\label{ug:exact-inverse}
 \left|\frac{\mathcal W(s,t)}{h_0-h_1}\right|
 +\left|\frac{(h_0+h_1)\mathcal W(s,t)}{h_0-h_1}\right|
 \ll (1+|\Im s|+|\Im t|)^C e^{-|\Im s+\Im t|}.
\end{equation}
In the first quotient, the ratios of the two numerator gamma
factors to $\Gamma(a)\Gamma(b)$ have polynomial growth. The
exponential from $\Gamma(s+t-\lambda-\mu)^{-1}$ cancels that
from $\sin(\pi(a+b-1)/2)^{-1}$, and
$q^{-\Re(a+b-1)}\le1$. In the second quotient the extra
exponential is
\[
 \exp\!\left\{\frac{\pi}{2}
 \bigl(|\Im s-\Im t|-|\Im s|-|\Im t|\bigr)\right\}\le1.
\]
The estimates are uniform on the disks because
$\Re(a+b-1)$ stays in a compact subinterval of $(0,1)$.

On $|\Im s|,|\Im t|\le5\log(2q)$, apply Cauchy's formula to the
weighted error in~\eqref{ug:parity-recovery}. Integrating
\eqref{ug:exact-inverse} costs at most a fixed power of $\log(2q)$,
so the resulting error is $O_\varepsilon(q^{1-\omega+\varepsilon})$.
Outside this box, bound the signed sum of logarithmic derivatives
by~\eqref{eq:logbound} and the main terms by~\eqref{ug:LDmodel}.
Both are $O_{\eta,\varepsilon}
 \bigl(q^{1+\varepsilon}(1+|\Im s|+|\Im t|)^C\bigr)$.
The decay of $\mathcal W$ makes the tail
$O_\varepsilon(q^{-4+\varepsilon})$.
Thus the integrated error required in~\eqref{ug:weighted-multiples}
holds uniformly in $\lambda,\mu$.

On the line $u+v=x$, direct integration gives
\[
 F(x)=x^{1-\lambda-\mu}w(x)\,
       \mathrm B(1-\lambda,1-\mu),
\]
where $\mathrm B$ is the beta function.
The Mellin transform of $F$ is
$\mathrm B(1-\lambda,1-\mu)\widetilde w(1+z-\lambda-\mu)$,
so~\eqref{ug:w-decay} supplies the decay needed for the
$D_q$ contour shifts. Equations~\eqref{ug:weighted-multiples}
and~\eqref{ug:exact-selection} now give, uniformly for
$0\le\lambda,\mu\le\omega$,
\begin{equation}\label{ug:real-power-test}
 \sum_{m+n=N}\frac{\Lambda(m)\Lambda(n)}{m^\lambda n^\mu}
 =\mathfrak S(N)N^{1-\lambda-\mu}\mathrm B(1-\lambda,1-\mu)
  +O_\varepsilon(N^{1-\omega+2\eta+\varepsilon}).
\end{equation}
Taking $\lambda=\mu=0$ proves the required estimate for the
von Mangoldt sum with the fixed-line saving.
\endgroup

\begingroup
To count ordered pairs of primes, integrate~\eqref{ug:real-power-test}
over $[0,\omega]^2$. For $t>0$ write
\[
 K_\omega(t)=\int_0^\omega t^{-\lambda}\,d\lambda
             =\frac{1-t^{-\omega}}{\log t},
 \qquad K_\omega(1)=\omega .
\]
The prime contribution is
$\sum_{p_1+p_2=N}(1-p_1^{-\omega})(1-p_2^{-\omega})$.
Its difference from the unweighted count is $O(N^{1-\omega})$,
using $\sum_{m\le N}m^{-\omega}\ll N^{1-\omega}$.
Pairs involving proper prime powers contribute
$O(N^{1/2}\log^2(2N))$, which is absorbed by the error in
\eqref{ug:real-power-test}. By the beta integral the integrated
main term is
$\mathfrak S(N)\int_0^N K_\omega(t)K_\omega(N-t)\,dt$.
On $[2,N-2]$, replacing the integrand by
$1/(\log t\log(N-t))$ costs $O(N^{1-\omega})$ before multiplication
by $\mathfrak S(N)$. The integrals over $[0,2]$ and $[N-2,N]$
are $O(1/\log N)$, since $K_\omega$ is integrable at zero.
Since $\mathfrak S(N)\ll_\varepsilon N^\varepsilon$, for even $N\ge6$
we obtain
\[
 \#\{(p_1,p_2):p_1+p_2=N\}
 =\mathfrak S(N)\int_2^{N-2}\frac{dt}{\log t\log(N-t)}
  +O_\varepsilon(N^{1-\omega+2\eta+\varepsilon}).
\]
\endgroup

\subsection{{Interpolation from the critical line}}
\label{ci:section}

{
The sum of the two main terms in~\eqref{ug:ordinary} is}
\begin{equation*}
 \mathcal M_q(a,b;s,t)=\varphi^*(q)Z_q(a,b;s,t)
       +H_q(a,b)Z_q(1-b,1-a;s,t).
\end{equation*}

{
For prime $q$, this uses the two terms of
Wang--Xu~\cite[Conjecture~2.1]{WX}, with the exact primitive parity
counts in the archimedean coefficient. On the critical line,
replacing the two counts by $\varphi^*(q)/2$ changes that coefficient
by $O(1)$, by the bound for $\Delta(q)$ above and
$|X_{q,e}(a)|=|X_{q,e}(b)|=1$. After taking out the common factor
$(q/\pi)^{1-a-b}$, which has modulus one on the critical line,
the difference is independent of $q$ and vanishes on $a+b=1$.
Cauchy's estimate for the gamma factors bounds its quotient
by $a+b-1$ polynomially in the heights. This removes the pole
of the reflected kernel, and the change in the main term is
$O_{\eta,\varepsilon}(q^\varepsilon)$ at the logarithmic heights in
Hypothesis~\ref{hyp:int-D}. The hypothesis uses the exact coefficient
and covers composite conductors as well.}

\begin{proof}[Proof of Lemma~\ref{ug:ratios}]
\begingroup
Fix $s,t$ with $\Re s=\Re t=\sigma$ and
$|\Im s|,|\Im t|\le10\log(2q)$. The normalized error
\[
 F_q(a,b)=q^{-1}\bigl(R_q(a,b;s,t)-\mathcal M_q(a,b;s,t)\bigr)
\]
is holomorphic in a neighborhood of
$|\Re a-1/2|,|\Re b-1/2|\le1/8$. The zeta factors in the
denominators in~\eqref{eq:Zkernel} have real parts at least
$7/8+\eta$ there, and the products in~\eqref{eq:Akernel}
converge normally. Along $a+b=1$ the two kernels have
opposite residues and $H_q(a,1-a)=\varphi^*(q)$, so the
apparent pole is removable.

Throughout these strips, there is an absolute constant $D$ with
\begin{equation}\label{ci:crude}
 |F_q(a,b)|\ll_\eta q^2
 \bigl(1+|\Im a|+|\Im b|+|\Im s|+|\Im t|\bigr)^D .
\end{equation}
For the ratio, use convexity for the numerator $L$-functions and
the bound under GRH
$L(\sigma+iy,\chi)^{-1}\ll_{\eta,\varepsilon}
 (q(2+|y|))^\varepsilon$ for the denominator factors.
For the main terms, use Stirling's formula, the zeta bounds in
fixed strips, and normal convergence of the arithmetic products.
The finite Euler factors cost $q^\varepsilon$. Near $a+b=1$,
the same polynomial bound follows by Cauchy's formula on a
small circle in $a$ contained in a slightly wider strip.
The bounds on the $L$-functions are standard under GRH
\cite[Chapter~5]{IK}.

Let $(a,b)$ range over either of the two pairs of numerator disks
in Lemma~\ref{ug:ratios}, with $\eta+r<1/8$.
Apply the three-lines theorem on $0\le\Re \xi\le1/8$ to
\[
 e^{\xi^2}F_q\!\left(
 \tfrac12+i\Im a+\frac{\Re a-1/2}{\eta+r}\,\xi,\,
 \tfrac12+i\Im b+\frac{\Re b-1/2}{\eta+r}\,\xi
 \right).
\]
The real coefficients of $\xi$ have modulus at most one.
On $\Re \xi=0$ both numerator arguments are on the critical line.
For $|\Im \xi|\le80\log(2q)$ their heights are less than
$100\log(2q)$, since the disk centers have heights at most
$10\log(2q)$ and the radii are fixed and small.
Hypothesis~\ref{hyp:int-D} then bounds the expression by
$O_{\eta,\varepsilon}(q^{-\omega_D+\varepsilon})$.
For larger $|\Im \xi|$ the Gaussian factor and~\eqref{ci:crude}
give the same bound. On $\Re \xi=1/8$, \eqref{ci:crude} gives
$O_{\eta,\varepsilon}(q^{2+\varepsilon})$ after the Gaussian
absorbs the height dependence. Evaluating at $\xi=\eta+r$ gives
\begin{equation}\label{ci:extensionformula}
 |F_q(a,b)|\ll_{\eta,\varepsilon}
 q^{-\omega_D+8(\eta+r)(2+\omega_D)+\varepsilon}.
\end{equation}
Multiplying by $q$ proves Lemma~\ref{ug:ratios} with
$\omega=\omega_D-8(\eta+r)(2+\omega_D)>2\eta$ for sufficiently
small $\eta$. The error is holomorphic on a neighborhood of the
two closed disks.
\endgroup
\end{proof}
%\textcolor{red}{The order is a bit confusing. Why the proof of Lemma 3.1 is at the end now?  }

{
Given $\varepsilon>0$, choose $r=\eta/8$ and then choose $\eta>0$
so small that $\omega_D-\omega+2\eta<\varepsilon/2$.
With the remaining $\varepsilon/2$ in the fixed-line errors,
Section~\ref{ug:exact-section} gives
$O_\varepsilon(N^{1-\omega_D+\varepsilon})$ for the von Mangoldt
sum and the ordered prime count. This proves
Theorem~\ref{thm:goldbach}.}

\section{Prime Pairs}\label{sec:twins}
%\textcolor{red}{A time average} of a short Dirichlet polynomial can distinguish a fixed additive difference. We choose a signed Fourier filter that keeps $n$ and $n+h$ while removing the diagonal. After the ratios prediction evaluates the average, the same filter extracts the coefficient $\mathfrak S(h)$ from $D_1$.

{
We prove Theorem~\ref{thm:prime-pairs} under RH and
Hypothesis~\ref{hyp:int-Z}. Fix $0<\eta<\min\{\eta_*,\omega_Z/2\}$.
We first evaluate a smoothed sum using the Mellin line $\Re a=\eta$.
Proposition~\ref{qt:sharp} below gives the sharp cutoff, and
Proposition~\ref{wk:zeta} gives the qualitative version.\par}

\subsection{A signed Fourier filter}
The filter has Fourier support near $h$ and $-h$. A narrow spatial support then excludes all other integral differences. Its two zero moments have separate purposes: one removes the Dirichlet-polynomial diagonal, and the other removes the pole of the singular-series Dirichlet series.

Use the Fourier convention $\widehat g(v)=\int_\R g(t)e^{-2\pi ivt}\,dt$. Choose a nonzero real even Schwartz function $k$ with smooth Fourier transform supported in $(-\epsilon_0,\epsilon_0)$, where $2\epsilon_0<1/8$. Let $M$ be a sufficiently large fixed integer and set
\begin{equation*}
 \begin{gathered}
 r_0(t)=t^M k(t),\qquad H=\widehat{r_0^2},\qquad
 g_h(t)=r_0(t)^2\cos(2\pi ht),\\
 \phi_h(v)=\widehat g_h(v)=\frac{H(v-h)+H(v+h)}2.
 \end{gathered}
\end{equation*}
Here $g_h$ is even and Schwartz, vanishes to order $2M$ at zero, and
\begin{equation}\label{eq:two-cancellations}
 \begin{gathered}
 \phi_h(0)=0,\qquad \int_0^\infty g_h(t)\,dt=0,\\
 g_h(0)=0,\qquad \int_0^\infty\phi_h(v)\,dv=0.
 \end{gathered}
\end{equation}
Since $H(0)=\int_\R r_0(t)^2\,dt>0$, there is a sufficiently small interval $J_h=(1,1+\epsilon_h)$ on which
\begin{equation*}
 a_h(u)=\tfrac12H(h/u-h)>0.
\end{equation*}
Shrink it so that $\phi_h(j/u)=0$ for every integer $j\geq1$, $j\ne h$, and every $u\in J_h$.

For real coefficients $c(n)$ and $W\in C_c^\infty(J_h)$ real, put
\begin{equation*}
 P_{c,W,X}(t)=\sum_n\frac{c(n)}{n^{1/2+it}}W(n/X),\qquad
 T=2\pi X,\qquad I_c(X)=\int_0^\infty g_h(t/T)|P_{c,W,X}(t)|^2\,dt.
\end{equation*}
Fourier inversion, using evenness of $|P|^2$, gives an exact sum over differences. The diagonal vanishes by $\phi_h(0)=0$. For all sufficiently large $X$, only the difference $h$ remains, and
\begin{align}\label{eq:exact-difference}
 I_c(X)=\frac T2\sum_n\frac{c(n)c(n+h)}{\sqrt{n(n+h)}}
 W(n/X)W((n+h)/X)\\
 {}\times H\bigl(X\log(1+h/n)-h\bigr).\notag
\end{align}
Indeed the positive Fourier lobe is near $h$, and $X\log(1+j/n)$ is below this lobe for $j<h$ and above it for $j>h$, by monotonicity and the choice of $J_h$. If $c(n)\ll_\varepsilon n^\varepsilon$, Taylor expansion in~\eqref{eq:exact-difference} yields
\begin{equation}\label{eq:approx-difference}
 \frac{I_c(X)}T=\sum_n\frac{c(n)c(n+h)}n W(n/X)^2a_h(n/X)
              +O_{W,h,\varepsilon}(X^{-1+\varepsilon}).
\end{equation}
For $c=\Lambda$ the error is $O(X^{-1}\log^2X)$.

\subsection{The swapped moment}
We now evaluate the filtered average for $c=\Lambda$. Mellin inversion expresses its Dirichlet polynomial in terms of $-\zeta'/\zeta$. The unswapped moment vanishes against the time weight; the swapped moment becomes the coefficient integral of Lemma~\ref{co:extraction}.

{
Take $\Re a=\Re b=\eta$ and $|\Im a|,|\Im b|\le U^{1/2}$.
For sufficiently large $U$, Lemma~\ref{co:Z-reduction} gives the
power form of Hypothesis~\ref{hyp:int-Z} on these lines.
The numerator shifts range over disks of fixed radius
$\min\{\rho,\eta/2\}$. Cauchy's formula and the kernel identities in
Section~\ref{sec:prelim} give, with $z=a+b$,
\begin{align}\label{eq:diff-zeta}
 \int_0^U\frac{\zeta'}\zeta(1/2+a+it)
              \frac{\zeta'}\zeta(1/2+b-it)\,dt
 &=UP_1(z)+(2\pi)^z\frac{U^{1-z}}{1-z}C_1(z)
   +O_{\eta,\varepsilon}(U^{1-\omega_Z+\varepsilon}).
\end{align}
\par}

For $c=\Lambda$, shifting a Mellin contour gives
\begin{align}\label{eq:zeta-Mellin}
 P_{\Lambda,W,X}(t)&=Q_X(t)+R_X(t),\\
 Q_X(t)&=-\frac1{2\pi i}\int_{(\eta)}
       \frac{\zeta'}\zeta(1/2+a+it)X^a\widetilde W(a)\,da,\notag\\
 R_X(t)&=X^{1/2-it}\widetilde W(1/2-it).\notag
\end{align}
The residue $R_X$ comes from the pole of $\zeta$ at one. {RH ensures that the shift crosses no zeros.} The bounds of Section~\ref{sec:prelim} give
\begin{equation*}
 Q_X(t)\ll_W X^\eta\log(2+|t|),\qquad
 R_X(t)\ll_{W,A}X^{1/2}(1+|t|)^{-A}.
\end{equation*}
Since $g_h(t/T)\ll(t/T)^{2M}$, the terms containing $R_X$ contribute $O_W(X^{-1})$ after taking $M$ and $A$ large enough.

Write $\widetilde g_h$ for the Mellin transform of the time weight. Apply~\eqref{eq:weighted-moment} to~\eqref{eq:diff-zeta} inside the two Mellin integrals for $|Q_X|^2$. The first cancellation in~\eqref{eq:two-cancellations} gives $\widetilde g_h(1)=0$. Combining the Mellin variables, and using $T=2\pi X$, leaves
\begin{equation}\label{eq:twin-mainintegral}
 I_\Lambda(X)=\frac T{2\pi i}\int_{(2\eta)}
       \widetilde{W^2}(z)C_1(z)\widetilde g_h(1-z)\,dz
       +O_{W,h,\varepsilon}(X^{1-{\omega_Z}+2\eta+\varepsilon}).
\end{equation}
Here the Mellin convolution of $\widetilde W(a)$ and $\widetilde W(z-a)$ is $\widetilde{W^2}(z)$. The factor $X^{2\eta}$ in the error is the cost of the two Mellin lines.

Define $\Phi_h(z)=\int_0^\infty\phi_h(v)v^{z-1}\,dv$. The cosine Mellin transform gives
\begin{equation}\label{eq:cosine-Mellin}
 \Phi_h(z)=\Xi(z)\widetilde g_h(1-z),\qquad 0<\Re z<1,
\end{equation}
initially with an exponential convergence factor. Continuation extends the identity to $0<\Re z<2M+1$. In addition, $\Phi_h(1)=0$ by~\eqref{eq:two-cancellations}. Thus $C_1(z)\widetilde g_h(1-z)=D_1(z)\Phi_h(z)$ and the pole of $D_1$ at one is canceled.

The Dirichlet-series expansion~\eqref{co:series} identifies the coefficients of $D_1$ with the Hardy--Littlewood singular series, including the local factor at two. Lemma~\ref{co:extraction}, with $D=D_1$, now gives the exact model contribution
\begin{align*}
 T\sum_{j\geq1}\mathfrak S(j)\int_0^\infty W(u)^2\phi_h(j/u)\,\frac{du}{u}
 =T\mathfrak S(h)\int_0^\infty W(u)^2a_h(u)\,\frac{du}{u}.
\end{align*}
The support construction selects $j=h$. Equations~\eqref{eq:approx-difference} and~\eqref{eq:twin-mainintegral} now imply
\begin{align}\label{eq:twin-local}
 \sum_n\frac{\Lambda(n)\Lambda(n+h)}n W(n/X)^2a_h(n/X)
 =\mathfrak S(h)\int_0^\infty W(u)^2a_h(u)\,\frac{du}{u}
    +O_{W,h,\varepsilon}(X^{-{\omega_Z}+2\eta+\varepsilon}).
\end{align}

{
For $f\in C_c^\infty(J_h)$ real, take
$W(u)=f(u)\sqrt{u/a_h(u)}$ in~\eqref{eq:twin-local}.
This gives the smoothed asymptotic for the test $f^2$.
If $\Phi$ is real and supported in $J_h$, choose
$\chi\in C_c^\infty(J_h)$ with $\chi=1$ on $\supp\Phi$.
Then $\Phi=((\chi+\Phi)^2-(\chi-\Phi)^2)/4$.
A finite partition of unity and rescaling $X$ give, for every
$\Phi\in C_c^\infty((0,\infty))$,
\begin{equation}\label{eq:twin-smooth}
 \sum_n\Lambda(n)\Lambda(n+h)\Phi(n/X)
 =\mathfrak S(h)X\int_0^\infty\Phi(u)\,du
  +O_{\Phi,h,\eta,\varepsilon}(X^{1-\omega_Z+2\eta+\varepsilon}).
\end{equation}
For complex $\Phi$ we apply the same argument to the real and
imaginary parts. On each fixed compact support the constants are
bounded polynomially in finitely many seminorms of $\Phi$.\par}

{
To pass to sharp cutoffs, we need uniform estimates for the
Mellin transforms of the varying weights. The same estimates apply
to the angular detector in Section~\ref{sec:quadratic}.\par}

\begin{proposition}[Sharp cutoffs without loss of the ratios saving]\label{qt:sharp}
Assume the relevant Riemann Hypothesis and either $\mathrm{RC}_Z(\omega)$ or $\mathrm{RC}_H(\omega)$. For sufficiently small fixed $\eta>0$, set $d=\omega-2\eta>0$. The corresponding sharp weighted error is
\[
 O_\varepsilon(X^{1-d+\varepsilon}).
\]
Consequently it is $O_\varepsilon(X^{1-\omega+\varepsilon})$ after choosing $\eta$ in terms of $\varepsilon$. In the continuous case the implied constants may also depend on the fixed difference $h$.
\end{proposition}
\begin{proof}
We give the dependence on the cutoff, rather than using an unspecified smooth seminorm. On a fixed compact interval in $(0,\infty)$, choose upper and lower smooth cutoffs with transition width $v=X^{-d}$. The radial weights $W_v$ in either detector, after multiplication by fixed smooth factors, may be chosen to satisfy
\[
 \|W_v\|_\infty+\|W_v'\|_1\ll1,\qquad
 \|W_v^{(j)}\|_1\ll_j v^{1-j}\quad(j\ge1).
\]
For prime pairs one can use squares of smooth cutoffs and a fixed partition into squares, so that taking a square root causes no loss in these bounds. Integration by parts in logarithmic coordinates gives, on each of the fixed Mellin lines used in the detectors,
\begin{equation}\label{qt:Mellin-cutoff}
 |\widetilde W_v(b+it)|\ll_A
 \frac1{1+|t|}\min\{1,(v|t|)^{-A}\},\qquad
 \int_\mathbb R|\widetilde W_v(b+it)|\,dt\ll\log(2/v).
\end{equation}
The minimum is interpreted as one at $t=0$.

{ 
Choose $d<b_0<1/2$ and $2b_0<a_0<1$, and truncate the
Mellin heights at $X^{b_0}$. For $Q\ge X^{a_0}$ these heights are
at most $Q^{1/2}$. In the continuous case apply
Hypothesis~\ref{hyp:int-Z} with the power form supplied by
Lemma~\ref{co:Z-reduction}. In the angular case apply
Hypothesis~\ref{hyp:int-H}. Cauchy differentiation costs at most a
fixed power of $\log(2Q)$. Integration by parts against the time or
family weight gives a weighted error $O_\varepsilon(X^{1-\omega+\varepsilon})$
uniformly over the retained Mellin heights. For $Q<X^{a_0}$ use the
direct bounds for the logarithmic derivatives and the main terms.
A sufficiently high fixed order of vanishing of the time or family
weight at zero makes the contribution of this range $O(X^{-A})$
for any prescribed $A$. For the pole residue in~\eqref{eq:zeta-Mellin},
use $R_X(t)\ll X^{1/2}$ and the zero of the time weight at the
origin for $t\le X^{a_0}$. For $t>X^{a_0}$ use
\eqref{qt:Mellin-cutoff} on $\Re s=1/2$, where
$vt\ge X^{a_0-d}$. Together with the polynomial bound for
$Q_X(t)$, these estimates make the residue contribution
$O(X^{-A})$ for any prescribed $A$.\par}

The Mellin tails are also negligible: by~\eqref{qt:Mellin-cutoff} each additional integration by parts supplies a factor $X^{-(b_0-d)}$, which beats the fixed polynomial growth of the actual moments and their models. The two Mellin integrals in the continuous case cost $O(\log^2X)$, and the one integral in the angular case costs $O(\log X)$. Thus~\eqref{eq:twin-mainintegral} and~\eqref{eqc:main-integral} retain the error
$O_\varepsilon(X^{1-\omega+2\eta+\varepsilon})$ for these varying cutoffs.

The elementary Taylor errors in the detectors are at most $O(v^{-1}\log^2X)$ after conversion to the coefficient sums. This is absorbed, since $d\le1/2$ and $v^{-1}=X^d\le X^{1-d}$. The endpoint intervals contain $O(vX+1)$ integers, each with weight $O(\log^2X)$ in the prime-pair problem and $O(\log X)$ in the quadratic problem. Their contribution and the change in the model integral are $O_\varepsilon(X^{1-d+\varepsilon})$. A fixed partition of unity and a dyadic sum finish the argument. Choose $2\eta<\varepsilon/2$ and use the remaining $\varepsilon/2$ to compensate logarithmic factors.
\end{proof}

{
Applying Proposition~\ref{qt:sharp} to~\eqref{eq:twin-smooth}
proves Theorem~\ref{thm:prime-pairs}.\par}

For prime pairs, proper prime powers contribute $O_h(X^{1/2}\log^3X)$, which is absorbed. Partial summation with the exact weight $1/(\log t\log(t+h))$ therefore gives
\begin{equation}\label{eq:twin-primecount}
 \#\{p\le X:p+h\text{ prime}\}
 =\mathfrak S(h)\int_2^X\frac{dt}{\log t\log(t+h)}
 +O_{h,\varepsilon}(X^{1-\omega_Z+\varepsilon})
\end{equation}
for fixed even $h$.
The two cancellations in the detector retain their separate roles: $\phi_h(0)=0$ removes the diagonal, whereas $\Phi_h(1)=0$ removes the singular-series pole.

\section{Primes between consecutive squares}\label{si:sec}

We prove Theorem~\ref{thm:squares} by averaging
$|\psi(x+h)-\psi(x)-h|^2$ over $X\le x\le X+h$. As in the work of
Goldston and Montgomery and of Heath-Brown~\cite{GM,HB}, the explicit
formula relates this mean square to a sum over pairs of zeros. A test
function with nonnegative Fourier transform lets us restrict the
resulting mean square of the zero sum to a short interval about $X$.
Monotonicity of $\psi$ then gives the estimate at $X=n^2$.

\subsection{Differentiating the ratios formula}

For large $H$, write $Z=\zeta'/\zeta$ and take
\begin{equation}\label{si:slice}
 \sigma_H=\frac1{\log H},\qquad d=\tfrac12+\sigma_H,\qquad
 {   r\in\mathbb R.}
\end{equation}
{   For $|u|,|v|\le\rho$, use the power form $\mathcal M_Z$
in~\eqref{co:Zmodel} to define $\Delta_H(u,v;r)$ on $[H,2H]$ by}
\begin{align}
 &\int_H^{2H}\frac{\zeta(d+it+u)\zeta(d-it+ir+v)}
                   {\zeta(d+it)\zeta(d-it+ir)}\,dt\notag \\
&\qquad =\mathcal M_Z(2H;\sigma_H,\sigma_H+ir;u,v)-\mathcal M_Z(H;\sigma_H,\sigma_H+ir;u,v)
   +\Delta_H(u,v;r).\label{si:RC}
\end{align}

\begin{theorem}\label{si:main}
Assume RH. Suppose that for some $0<\omega\le1/2$ and every
$\varepsilon>0$,
\begin{equation}\label{si:powerhyp}
 \sup_{\substack{|r|\le H^\omega\\ {   |u|,|v|\le\sigma_H/2}}}
 |\Delta_H(u,v;r)|\ll_{\omega,\varepsilon}H^{1-\omega+\varepsilon}.
\end{equation}
Then, for every $\varepsilon>0$,
\begin{align}
 \psi((n+1)^2)-\psi(n^2)
 &=2n+1+O_{\omega,\varepsilon}(n^{1-\omega/3+\varepsilon}),
 \label{si:Legendre-quant}\\
 \pi((n+1)^2)-\pi(n^2)
 &=\int_{n^2}^{(n+1)^2}\frac{dt}{\log t}
   +O_{\omega,\varepsilon}(n^{1-\omega/3+\varepsilon}),\notag
\end{align}
and $p_{j+1}-p_j\ll_{\omega,\varepsilon}
p_j^{(1-\omega)/(2-\omega)+\varepsilon}$. 
\end{theorem}

{  
Under $\mathrm{RC}_Z(\omega_Z)$ and RH, apply
Lemma~\ref{co:Z-reduction} at $Q=H$ and $Q=2H$ with
$a=\sigma_H$ and $b=\sigma_H+ir$ at both endpoints.
For $|r|\le H^{\omega_Z}$ these shifts lie in the stated height
range at both endpoints. For large $H$, the numerator disks have
radius $\sigma_H/2$ at both endpoints. Subtracting the two
formulas gives~\eqref{si:powerhyp} with $\omega=\omega_Z$.\par}

Let
\begin{equation}\label{si:height-primitive}
 W_H(z)=\int_H^{2H}(t/2\pi)^{-z}\,dt
       =(2\pi)^z\frac{(2H)^{1-z}-H^{1-z}}{1-z},\qquad
 K_H(z)=H P_1(z)+W_H(z)C_1(z),
\end{equation}
with $P_1$ and $C_1$ as in
\eqref{co:Pcontinuation} and~\eqref{co:Ccontinuation}.

\begin{lemma}\label{si:diff}
Assume RH. For the shifts in~\eqref{si:slice}, we have
\[
 \mathcal Z_H(r):=\int_H^{2H} Z(d+it)Z(d-it+ir)\,dt
 =K_H(2\sigma_H+ir)+E_H(r).
\]
{  
For every $0<\rho'\le\rho$,
\begin{equation}\label{si:differentiated}
 |E_H(r)|\le (\rho')^{-2}
       \sup_{|u|=|v|=\rho'}|\Delta_H(u,v;r)|.
\end{equation}\par}

\end{lemma}
\begin{proof}
For the moment, write $z=2\sigma_H+ir$. By~\eqref{co:Zmodel}, the difference of
the main terms in~\eqref{si:RC} is
\[
 \begin{aligned}
 &\mathcal M_Z(2H;\sigma_H,\sigma_H+ir;u,v)
  -\mathcal M_Z(H;\sigma_H,\sigma_H+ir;u,v)\\
 &\quad=H Z_1(d+u,d+ir+v;d,d+ir)+W_H(z+u+v)
   Z_1(1-d-ir-v,1-d-u;d,d+ir).
 \end{aligned}
\]
The possible numerator poles come from
$\zeta(1+z+u+v)$ and $\zeta(1-z-u-v)$. On $z+u+v=0$,
\[
 d+u=1-d-ir-v,\qquad d+ir+v=1-d-u,\qquad W_H(0)=H.
\]
The remaining factors of the two terms therefore agree by
\eqref{eq:Zkernel}, while the two zeta factors have opposite
residues. The poles cancel. RH excludes poles of the reciprocal
zeta factors on the numerator disks, and the arithmetic products
converge normally there. The integral on the left of
\eqref{si:RC} is also holomorphic on these disks, since its
numerator arguments have real part less than $1$ and its
denominators are nonzero. Hence $\Delta_H(u,v;r)$ is holomorphic
on a neighborhood of the closed disks.

The kernel identities in Section~\ref{sec:kernels} give
\[
 \left.\partial_u\partial_v
 Z_1(d+u,d+ir+v;d,d+ir)\right|_{u=v=0}=P_1(z)
\]
and
\[
 \left.\partial_u\partial_v
 \bigl\{W_H(z+u+v)Z_1(1-d-ir-v,1-d-u;d,d+ir)\bigr\}
 \right|_{u=v=0}=W_H(z)C_1(z).
\]
For the second identity, the reflected kernel contains
\[
 \frac1{\zeta(1-u)\zeta(1-v)}
 =uv+O\!\left(|uv|(|u|+|v|)\right),
\]
so both derivatives must fall on these reciprocal zeta factors.
Differentiating~\eqref{si:RC} at $u=v=0$ thus yields
\[
 \mathcal Z_H(r)
 =HP_1(z)+W_H(z)C_1(z)
  +\partial_u\partial_v\Delta_H(0,0;r)
 =K_H(z)+\partial_u\partial_v\Delta_H(0,0;r).
\]
It follows that $E_H(r)=\partial_u\partial_v\Delta_H(0,0;r)$.
{  
Apply Cauchy's formula on $|u|=|v|=\rho'$ to obtain
\[
 \partial_u\partial_v\Delta_H(0,0;r)
 =\frac1{(2\pi i)^2}
  \int_{|u|=\rho'}\int_{|v|=\rho'}
  \frac{\Delta_H(u,v;r)}{u^2v^2}\,dv\,du.
\]
This gives~\eqref{si:differentiated}.\par}

\end{proof}

\subsection{A pair correlation estimate}

For $M\ge2$, convolve the triangular weight on $[-M,M]$ with
a Gaussian and normalize the resulting test function at zero. Define
\[
 \varphi_M(r)=\frac{c_M}{\sqrt\pi}
     \int_{-M}^M(1-|y|/M)e^{-(r-y)^2}\,dy,\qquad
 \varphi_M(0)=1.
\]
Thus $c_M\asymp1$. The Gaussian permits vertical shifts of the
argument. With $\widehat f(v)=\int_{\mathbb R}f(r)e^{-irv}\,dr$,
the convolution formula gives
\begin{equation}\label{si:kernel}
 \widehat\varphi_M(v)=c_M M
  \left(\frac{\sin(Mv/2)}{Mv/2}\right)^2 e^{-v^2/4}\ge0.
\end{equation}
In particular, $\widehat\varphi_M(v)\gg M$ for $|v|\le1/M$. The
integral defining $\varphi_M$ extends to an even entire function
and gives
\begin{equation}\label{si:kernelbounds}
 |\varphi_M(r+ib)|\le e^{b^2}\varphi_M(r),\qquad
 \int_{\mathbb R}(1+|r|)^j\varphi_M(r)\,dr\ll_j M^{j+1}
 \quad (j\ge0).
\end{equation}
For $|r|\ge 2M$ and bounded $b$, it also gives
$|\varphi_M(r+ib)|\ll_b M e^{-r^2/4}$.

\begin{lemma}\label{si:tested}
Assume RH. Fix $K_0>1$, $0<c<1/2$, and $\varepsilon>0$. For
sufficiently large $H$, let $H/4\le V\le H^{K_0}$ and
$\log H\le M\le H^{1/2}/4$. Then we have
\begin{align}
 P_H(V,M)&:=\sum_{H<\gamma,\gamma'\le2H}
         V^{i(\gamma-\gamma')}\varphi_M(\gamma-\gamma')
 =\frac1{2\pi}\int_H^{2H}\log(t/2\pi)\,dt
 \label{si:testedformula}\\
 &\qquad+O_{K_0,c,\varepsilon}\left(
 M\sup_{|r|\le2M}|E_H(r)|
 +HV^{-c}M^{1+\varepsilon}
 +H(H/V)^cM^\varepsilon+M^2(\log H)^B\right),\notag
\end{align}
where $\gamma,\gamma'$ are ordinates of nontrivial zeros counted with
multiplicity, and $B$ is an absolute constant.
\end{lemma}

\begin{proof}
Put $f(r)=\varphi_M(r)\cos(r\log V)$. By symmetry,
\[
 P_H(V,M)=\sum_{H<\gamma,\gamma'\le2H}f(\gamma-\gamma').
\]
Choose $T_1\in[H-1,H+1]$ and $T_2\in[2H-1,2H+1]$ at distance
$\gg1/\log H$ from the zero ordinates, and set $I=[T_1,T_2]$.
Let $\mathcal R$ be the rectangle with vertical sides at $1-d,d$
and horizontal sides at $T_1,T_2$, oriented positively.

By the classical estimates for the function $Z=\zeta'/\zeta$ and RH, we have
\[
 Z(s)\ll(\log H)^2
 \quad\bigl(s\in\partial\mathcal R
 \text{ or } d\le\Re s\le2,\ \Im s\in I\bigr).
\]
Moreover, by ~\eqref{si:kernelbounds} and the condition $V\le H^{K_0}$, we have
\[
 |f(r+ib)|\ll_{K_0}\varphi_M(r)\quad(|b|\le2\sigma_H),
 \qquad \int_{\mathbb R}\varphi_M(r)\,dr\ll M.
\]
Hence the double contour integral with at least one variable on
a horizontal side is $O_{K_0}(M(\log H)^4)$.

The zero-counting bound and the Gaussian decay of $\varphi_M$
also give
\[
 \sum_{H-1<\gamma\le2H+1}\varphi_M(t-\gamma)
 \ll M\log H\qquad(H-1\le t\le2H+1).
\]
Applying the residue theorem to both variables and replacing
$T_1,T_2$ by $H,2H$ in the resulting sum, we obtain
\[
 P_H(V,M)=\frac1{(2\pi i)^2}
 \int_{\partial\mathcal R}\int_{\partial\mathcal R}
 Z(s)Z(w)f(-i(s-w))\,ds\,dw
 +O\!\left(M(\log H)^2\right).
\]

{   Write $G(t)=\chi_\zeta'/\chi_\zeta(1-d+it)$.}
By the functional equation and Stirling's formula,
\[
 Z(1-d+it)=G(t)-Z(d-it),\quad
 G(t)=-\log(t/2\pi)+O(H^{-1}),\quad G'(t)\ll H^{-1}
 \quad(t\in I).
\]
Parametrize the vertical sides by their
ordinates $t,u$ and put $r=t-u$. The two mixed terms containing
$Z(d+it)Z(d-iu)$ agree after interchanging $t$ and $u$.
In each term the minus sign in the functional equation cancels the
reversed orientation on the left side. Together they contribute
\[
 \frac1{2\pi^2}\int_{\mathbb R} f(r-2i\sigma_H)
   \int_{I\cap(I+r)} Z(d+it)Z(d-i(t-r))\,dt\,dr .
\]
The intersection expresses the conditions $t,t-r\in I$.
For $|r|\le H/4$, replacing $I\cap(I+r)$ by $[H,2H]$
changes the inner integral by $O((1+|r|)(\log H)^B)$.
For $|r|>H/4$ the Gaussian decay gives a negligible error,
using the local bound for $Z$ on the line $\Re s=d$.
Integrating the endpoint loss against $\varphi_M(r)$ gives
$O_{K_0}(M^2(\log H)^B)$ by~\eqref{si:kernelbounds}.

For the term with $Z(d+it)Z(d+i(t-r))$, fix $r$ and translate the path
of integration in $t$ down by $i(2-d)$. Both arguments of $Z$
move to real part $2$ while the test $f(r)$ is fixed. The two
connecting segments contribute $O((\log H)^B)$ for each $r$,
since the heights stay in $I$ and both real parts are at least
$d$. On the translated path the integrand is
\[
 Z(2+it)Z(2+i(t-r))
 =\sum_{m,n\ge2}\frac{\Lambda(m)\Lambda(n)n^{ir}}{m^2n^2}
   (mn)^{-it}.
\]
Integrating $(mn)^{-it}$ over any interval costs at most
$2/\log(mn)$. Thus the integral in $t$ is $O(1)$ on the new
path, and integrating against $f(r)$ gives
$O_{K_0}(M(\log H)^B)$ including the connecting segments.
The conjugate term has the same bound. A shift of one factor also gives
$\int_b^{b'}Z(d+it)\,dt\ll(\log H)^B$ uniformly for $[b,b']\subseteq I$.
Integration by parts against $G(t-r)$ or $G(t)$ therefore bounds
the terms with one factor $G$ by $O_{K_0}(M(\log H)^B)$.

If $t,t-r\in I$, then
\[
 G(t)G(t-r)=\log^2(t/2\pi)+O((1+|r|)\log H/H).
\]
Using this approximation and extending the inner integral to
$[H,2H]$ costs $O_{K_0}(M^2(\log H)^B)$ by~\eqref{si:kernelbounds}.
The resulting term is negligible, since
$\int_{\mathbb R}f(r)\,dr=\widehat\varphi_M(\log V)
 \ll M e^{-(\log V)^2/4}$. We have proved
\begin{equation}\label{si:contour}
 P_H(V,M)=\frac1{2\pi^2}
   \int_{\mathbb R}f(r-2i\sigma_H)\mathcal Z_H(r)\,dr
   +O_{K_0}(M^2(\log H)^B).
\end{equation}

Insert~\eqref{si:differentiated} for $|r|\le2M$.
The Gaussian bound and the polynomial growth of the integrands
allow us to extend the integral of the main term to the full
line. The ratios error contributes
$O_{K_0}(M\sup_{|r|\le2M}|E_H(r)|)$ by~\eqref{si:kernelbounds}.
With $z=2\sigma_H+ir$, the cosine in $f(r-2i\sigma_H)$ is
$\tfrac12(V^z+V^{-z})$, and evenness gives
$\varphi_M(r-2i\sigma_H)=\varphi_M(iz)$. The resulting main term is
\begin{equation}\label{si:Jpm}
 \frac1{4\pi^2 i}\int_{(2\sigma_H)}
       \varphi_M(iz)(V^z+V^{-z})K_H(z)\,dz .
\end{equation}

By~\eqref{co:Pcontinuation} and~\eqref{co:Ccontinuation},
$P_1(z)=z^{-2}+O(1)$ and $C_1(z)=-z^{-2}+O(1)$ at zero.
For the latter expansion, $\mathcal A(0)=1$ and $\mathcal A'(0)=0$.
Under RH the only possible pole in $|\Re z|\le c$ is at zero,
since $\Re(1\pm z)>1/2$ and the remaining arithmetic factors are
holomorphic there. As $W_H(0)=H$, the coefficient of $z^{-2}$
cancels. The coefficient of $z^{-1}$ is
$-W_H'(0)=\int_H^{2H}\log(t/2\pi)\,dt$. On $\Re z=\pm c$,
the bound for $Z'$ to the right of the critical line and the
Lindel\"of bound under RH give
$P_1(z),C_1(z)\ll_{c,\varepsilon}(1+|\Im z|)^\varepsilon$,
using the normal convergence in~\eqref{co:Pcontinuation}
and~\eqref{co:Ccontinuation}.
In~\eqref{si:height-primitive}, the numerator has size
$O_c(H^{1\mp c})$ on $\Re z=\pm c$, and
$|1-(\pm c+ir)|\asymp_c1+|r|$. Hence
\begin{equation}\label{si:height-decay}
 |W_H(\pm c+ir)|\ll_c H^{1\mp c}(1+|r|)^{-1}.
\end{equation}
The integral defining $\varphi_M$ gives $\varphi_M(r)\ll1$ on
the real line. Use this on $|r|\le2M$ and the Gaussian bound
outside to obtain, uniformly for bounded $b$,
\begin{align}
 \int_{\mathbb R}|\varphi_M(r+ib)|(1+|r|)^\varepsilon\,dr
 &\ll_{b,\varepsilon} M^{1+\varepsilon},\notag\\
 \int_{\mathbb R}|\varphi_M(r+ib)|(1+|r|)^{-1+\varepsilon}\,dr
 &\ll_{b,\varepsilon} M^\varepsilon.\label{si:two-moments}
\end{align}
Move the line for $V^z$ to $-c$ and the line for $V^{-z}$ to $c$.
The first shift crosses zero and gives the main term of
\eqref{si:testedformula}, since $\varphi_M(0)=1$. Both powers of
$V$ have absolute value $V^{-c}$ on their new lines. The first
bound in~\eqref{si:two-moments} gives
$O_{c,\varepsilon}(HV^{-c}M^{1+\varepsilon})$ for the terms with
$P_1$. For the terms with $C_1$, the factor $(1+|r|)^{-1}$
in~\eqref{si:height-decay} gives the second bound in
\eqref{si:two-moments}. Their total contribution is
$O_{c,\varepsilon}(V^{-c}(H^{1+c}+H^{1-c})M^\varepsilon)$,
which proves~\eqref{si:testedformula}.
\end{proof}

\subsection{A mean square in short intervals}

Let
\[
 B_H(u)=\sum_{H<\gamma\le2H}u^{i\gamma}
 \qquad\text{and}\qquad
 S_t(u)=\sum_{0<\gamma\le t}u^{i\gamma}.
\]
Fourier inversion in the sum defining $P_H(V,M)$ gives
\[
 P_H(V,M)=\frac1{2\pi}\int_{\mathbb R}
       \widehat\varphi_M(v)|B_H(Ve^v)|^2\,dv.
\]
Suppose that $0<Y\le V/(2M)$. Then
$\log(1+Y/V)\le 1/M$, so by~\eqref{si:kernel},
\[
 \begin{aligned}
 P_H(V,M)
 &\gg M\int_0^{\log(1+Y/V)}|B_H(Ve^v)|^2\,dv\\
 &= M\int_V^{V+Y}|B_H(u)|^2\,\frac{du}{u}
 \gg \frac{M}{V}\int_V^{V+Y}|B_H(u)|^2\,du .
 \end{aligned}
\]
In particular,
\begin{equation}\label{si:localize}
 \int_V^{V+Y}|B_H(u)|^2\,du\ll \frac{V}{M}P_H(V,M)
 \qquad(0<Y\le V/(2M)).
\end{equation}

For $X>1$ and $h>0$, define
\[
 J(X,h)=\int_X^{X+h}|\psi(x+h)-\psi(x)-h|^2\,dx.
\]
The change in the phase of $x^{i\gamma}$ from $x$ to $x+h$ is
$\gamma\log(1+h/x)\asymp\gamma h/X$ for $x\asymp X$ and
$h\le X/10$. We split the sum over zeros in the explicit formula
at ${X/h}$.

\begin{lemma}\label{si:transfer}
Assume RH and let $2\le h\le X/10$. Suppose that for some
$0<b<2$ and $D\ge0$,
\begin{equation}\label{si:zeronorm}
 \int_V^{V+h}|S_t(u)|^2\,du\le Dt^b
 \qquad(X\le V\le X+h,\quad {X/h}\le t\le2X).
\end{equation}
Then
\begin{equation}\label{si:Jtransfer}
 J(X,h)\ll_b D\frac{h^2}{X}{(X/h)^b}+h(\log X)^4.
\end{equation}
\end{lemma}

\begin{proof}
Choose a good height $T\in[X/2,X]$. Under RH the truncated
explicit formula (see~\cite{IK}) gives, for $X\le u\le X+2h$,
\[
 \psi(u)=u-\sum_{|\gamma|\le T}
        \frac{u^{1/2+i\gamma}}{1/2+i\gamma}
        +O((\log X)^2).
\]
Indeed, the truncation error is
$O((u/T)\log^2(uT)+\log u)$, including the change from the
midpoint value to $\psi(u)$ at a prime power. Since $h\ge2$,
we have ${X/h}\le T$.

For $0<\gamma\le {X/h}$, we sum the identity
\[
 \frac{(x+h)^{1/2+i\gamma}-x^{1/2+i\gamma}}{1/2+i\gamma}
 =\int_0^h(x+y)^{-1/2+i\gamma}\,dy.
\]
Minkowski's inequality and~\eqref{si:zeronorm}, with $V=X+y$,
give
\[
 \begin{aligned}
 \left\|\int_0^h (x+y)^{-1/2}S_{{X/h}}(x+y)\,dy
       \right\|_{L^2([X,X+h])}
 &\le X^{-1/2}\int_0^h
       \|S_{{X/h}}\|_{L^2([X+y,X+h+y])}\,dy\\
 &\le hX^{-1/2}D^{1/2}{(X/h)^{b/2}}.
 \end{aligned}
\]

For ${X/h}<\gamma\le T$, partial summation gives
\[
 \sum_{{X/h}<\gamma\le T}\frac{u^{i\gamma}}{1/2+i\gamma}
 =\frac{S_T(u)}{1/2+iT}-\frac{S_{{X/h}}(u)}{1/2+i{X/h}}
   +i\int_{{X/h}}^T\frac{S_t(u)}{(1/2+it)^2}\,dt .
\]
Thus, for $X\le V\le X+h$, the $L^2([V,V+h])$ norm of this
sum is, by~\eqref{si:zeronorm}, at most 
\[
 \ll D^{1/2}\left(T^{b/2-1}+{(X/h)^{b/2-1}}
                +\int_{{X/h}}^T t^{b/2-2}\,dt\right)
 \ll_b D^{1/2}{(X/h)^{b/2-1}},
\]
where the last estimate uses $b<2$. In the difference of the
two explicit formulas, the high zeros give this sum multiplied by
$u^{1/2}$ at $u=x+h$ and at $u=x$. Using the norm bound with
$V=X+h$ and $V=X$ gives a total norm
\[
 \ll_b X^{1/2}D^{1/2}{(X/h)^{b/2-1}}
 =hX^{-1/2}D^{1/2}{(X/h)^{b/2}}.
\]
The negative ordinates give the conjugate sums. The difference
of the explicit-formula errors has norm $O(h^{1/2}(\log X)^2)$,
and hence
\[
 J(X,h)^{1/2}
 \ll_b hX^{-1/2}D^{1/2}{(X/h)^{b/2}}
       +h^{1/2}(\log X)^2.
\]
Squaring proves~\eqref{si:Jtransfer}.
\end{proof}

Let $0<\delta<1/2$ and $X-\delta h>1$. For
$X\le y\le X+\delta h$ we have
$(y,y+(1-\delta)h]\subseteq(X,X+h]$, and for
$X-\delta h\le y\le X$ we have
$(X,X+h]\subseteq(y,y+(1+\delta)h]$. Consequently,
\[
 \begin{aligned}
 \psi(X+h)-\psi(X)-h
 &\ge \psi(y+(1-\delta)h)-\psi(y)-(1-\delta)h-\delta h,
       &&(X\le y\le X+\delta h),\\
 \psi(X+h)-\psi(X)-h
 &\le \psi(y+(1+\delta)h)-\psi(y)-(1+\delta)h+\delta h,
       &&(X-\delta h\le y\le X).
 \end{aligned}
\]
Average each inequality over its range of $y$. By Cauchy--Schwarz,
the average of the absolute value of each centered increment
is at most $(\delta h)^{-1/2}$ times the square root of its
mean-square integral over that range. Since
$\delta h\le(1-\delta)h$, these integrals are bounded by the
two mean squares below. We obtain
\begin{equation}\label{si:monotonicity}
 |\psi(X+h)-\psi(X)-h|
 \le \delta h+\left(
 \frac{J(X,(1-\delta)h)+J(X-\delta h,(1+\delta)h)}{\delta h}
 \right)^{1/2}.
\end{equation}

\subsection{Proof of Theorem~\ref{si:main}}
{
Let $2\le h\le4\sqrt X$. For $X/h\le t\le2X$,
write $(0,t]$ as the union of dyadic blocks $(H,2H]$, with
$H=t/2,t/4,\ldots$. On the blocks with $H\ge X^{1/100}$ take
$M_H=H^\omega/100$. Since $H\le X$ and $\omega\le1/2$, for
sufficiently large $X$ we have
\[
 \log H\le M_H\le H^{1/2}/4,\qquad
 2M_H\le H^\omega,\qquad
 2M_Hh\le \frac{2}{25}H^\omega\sqrt X\le \frac{2X}{25}.
\]
Also $H/4\le V\le H^{200}$ for $X\le V\le X+h$. Hence we can
apply~\eqref{si:localize} with $Y=h$ and Lemma~\ref{si:tested}
with $K_0=200$.

{  
Apply Lemma~\ref{si:diff} with $\rho'=\sigma_H/2$
and use~\eqref{si:powerhyp} with $\varepsilon/8$.
The factor $(\rho')^{-2}=4(\log H)^2$ is absorbed in
$H^{\varepsilon/8}$, so
\[
 \sup_{|r|\le2M_H}|E_H(r)|
 \ll_{\omega,\varepsilon}H^{1-\omega+\varepsilon/4}.
\]\par}

Using~\eqref{si:testedformula} with $\varepsilon/4$ in place of
$\varepsilon$, and noting that $V\asymp X$, we obtain
\[
 \begin{aligned}
 \int_V^{V+h}|B_H(u)|^2\,du
 &\ll_{\omega,c,\varepsilon}
 XH^{1-\omega}\log H+XH^{1-\omega+\varepsilon/4}
     +X^{1-c}H M_H^{\varepsilon/4}\\
 &\qquad +X^{1-c}H^{1+c}M_H^{-1+\varepsilon/4}
     +XM_H(\log H)^B\\
 &\ll_{\omega,c,\varepsilon}
 X^{1+\varepsilon/2}H^{1-\omega}
 \left(1+\frac{H^\omega}{X^c}
           +\left(\frac{H}{X}\right)^c+H^{2\omega-1}\right).
 \end{aligned}
\]
The last two terms in parentheses are at most one. Choose
$0<c<1/2$ with $c>\omega-\varepsilon/2$, so that
$H^\omega/X^c\le X^{\omega-c}\le X^{\varepsilon/2}$. It follows
that each retained block contributes at most
$O_{\omega,\varepsilon}(X^{(1+\varepsilon)/2}H^{(1-\omega)/2})$
to the $L^2([V,V+h])$ norm of $S_t$. The remaining zeros have
ordinates at most $2X^{1/100}$. By the zero count, their contribution
to this norm is $O(h^{1/2}X^{1/100}\log X)=o(X^{1/2})$.
Minkowski's inequality therefore gives
\[
 \begin{aligned}
 \|S_t\|_{L^2([V,V+h])}
 &\ll_{\omega,\varepsilon} X^{(1+\varepsilon)/2}
    \sum_{H=t/2,t/4,\ldots} H^{(1-\omega)/2}
     +o(X^{1/2})\\
 &\ll_{\omega,\varepsilon}X^{(1+\varepsilon)/2}t^{(1-\omega)/2}
 \qquad(X\le V\le X+h,\quad X/h\le t\le2X).
 \end{aligned}
\]
Squaring and applying Lemma~\ref{si:transfer} with $b=1-\omega$, we get}
{
\begin{equation}\label{si:Jbound}
 J(X,h)\ll_{\omega,\varepsilon}
       X^{1-\omega+\varepsilon}h^{1+\omega}+h(\log X)^4.
\end{equation}
}

{
To bound the prime gaps, take $h=X^{(1-\omega)/(2-\omega)+\varepsilon}$
with $0<\varepsilon<\omega/(2(2-\omega))$. By~\eqref{si:Jbound},
\[
 \frac{J(X,h)}{h^3}
 \ll_{\omega,\varepsilon}
 X^{1-\omega+\varepsilon}h^{\omega-2}+h^{-2}(\log X)^4
 =X^{-(1-\omega)\varepsilon}+h^{-2}(\log X)^4=o(1).
\]
In an interval of length $O(\sqrt X)$ starting at $X$, the proper
prime powers have total von Mangoldt weight $O((\log X)^2)$.
For each $k\ge2$, the corresponding interval of $k$th roots
has length $O(X^{1/k-1/2})=O(1)$, and there are $O(\log X)$
possible exponents. If $(X,X+2h]$ contained no prime, then
$0\le\psi(x+h)-\psi(x)\ll(\log X)^2$ for all $X\le x\le X+h$.
Consequently,
\[
 \begin{aligned}
 J(X,h)
 &\ge h^3-2h\int_X^{X+h}\bigl(\psi(x+h)-\psi(x)\bigr)\,dx\\
 &\ge h^3-O\!\left(h^2(\log X)^2\right),
 \end{aligned}
\]
contradicting $J(X,h)=o(h^3)$. Taking $X=p_j$ gives
$p_{j+1}-p_j\le 2p_j^{(1-\omega)/(2-\omega)+\varepsilon}$ for
all sufficiently large $j$.

Suppose now that $\sqrt X\le h\le \tfrac52\sqrt X$. For
$0<\delta<1/2$, the lengths $(1-\delta)h$ and $(1+\delta)h$ in
\eqref{si:monotonicity} are both comparable to $\sqrt X$. The first
is at most $4\sqrt X$, and the second satisfies
\[
 (1+\delta)h\le\tfrac{15}{4}\sqrt X\le4\sqrt{X-\delta h}
\]
for sufficiently large $X$. Since $X-\delta h\asymp X$,
\eqref{si:Jbound} gives
\[
 J(X,(1-\delta)h)+J(X-\delta h,(1+\delta)h)
 \ll_{\omega,\varepsilon} X^{3/2-\omega/2+\varepsilon}.
\]
Thus
\[
 \begin{aligned}
 |\psi(X+h)-\psi(X)-h|
 &\ll_{\omega,\varepsilon}\delta X^{1/2}
     +\left(\frac{X^{3/2-\omega/2+\varepsilon}}
                  {\delta X^{1/2}}\right)^{1/2}\\
 &=O_{\omega,\varepsilon}\!\left(
   \delta X^{1/2}+\delta^{-1/2}X^{1/2-\omega/4+\varepsilon/2}\right).
 \end{aligned}
\]
Taking $\delta=X^{-\omega/6}$,
\begin{equation}\label{si:pointwise-power}
 \psi(X+h)-\psi(X)
 =h+O_{\omega,\varepsilon}(X^{1/2-\omega/6+\varepsilon/2}).
\end{equation}
For $X<p\le X+h$, we have $\log p=\log X+O(h/X)$. Summing the
error over at most $h+1$ primes costs $O(1)$ since $h\le\tfrac52\sqrt X$.
Due to the bound for proper prime powers, we therefore have 
\[
 (\log X)\bigl(\pi(X+h)-\pi(X)\bigr)
 =\psi(X+h)-\psi(X)+O((\log X)^2).
\]
Also,
\[
 \int_X^{X+h}\frac{du}{\log u}
 =\frac{h}{\log X}
   +O\!\left(\frac{h^2}{X(\log X)^2}\right).
\]
It follows from~\eqref{si:pointwise-power} that
\[
 \pi(X+h)-\pi(X)
 =\int_X^{X+h}\frac{du}{\log u}
  +O_{\omega,\varepsilon}\!\left(
      \frac{X^{1/2-\omega/6+\varepsilon/2}}{\log X}\right).
\]
Taking $X=n^2$ and $h=2n+1$ proves both estimates
in~\eqref{si:Legendre-quant}.\par}

\section{Quadratic prime values}\label{sec:quadratic}

{
We prove Theorem~\ref{thm:quadratic-main} under GRH for the
angular Hecke family and Hypothesis~\ref{hyp:int-H}. We use the
characters $\xi_k$, the functions $L_k$, and the model $G$ from
Section~\ref{co:angular-input}. The average over $k$ selects
Gaussian integers $n\pm i$. We evaluate the resulting smooth
sum of $\Lambda(n^2+1)$ on a fixed Mellin line and then apply
Proposition~\ref{qt:sharp} to pass to the sharp cutoff.\par}

\subsection{An angular detector}
We use the Fourier transform $\widehat g(v)=\int_{\mathbb R}g(t)e^{-2\pi itv}\,dt$.

\begin{lemma}\label{leqc:window}
For every integer $J\ge1$ and every sufficiently small fixed $\delta>0$, there is a real even Schwartz function $g$ such that $H=\widehat g$ is smooth,
\begin{equation}\label{eqc:window}
 \operatorname{supp}H\subset
 (1-\delta,1+\delta)\cup(-1-\delta,-1+\delta),
 \qquad H(1)\ne0,
\end{equation}
and $g^{(j)}(0)=0$ for $0\le j<J$. Moreover,
\begin{equation}\label{eqc:vanishing}
 \int_0^\infty g(u)\,du=0,
 \qquad \int_0^\infty H(u)\,du=0.
\end{equation}
\end{lemma}
\begin{proof}
{
Choose a real $H_+\in C_c^\infty((1-\delta,1+\delta))$ with
$H_+(1)\ne0$. Add smooth bumps of integral one supported near
distinct points of $(1-\delta,1+\delta)\setminus\{1\}$, with
coefficients chosen so that
\[
 \int_0^\infty v^{2j}H_+(v)\,dv=0\qquad(0\le j\le J).
\]
With $J+1$ bumps, the moment matrix tends to the Vandermonde
matrix in the squares of their centers as the supports shrink.
It is therefore invertible for sufficiently small supports.
Extend $H_+$ evenly to $H$ and take its inverse Fourier transform
$g$. The function $g$ is real and even, and
\[
 g^{(2j)}(0)=2(2\pi i)^{2j}\int_0^\infty v^{2j}H_+(v)\,dv.
\]
Together with evenness this gives the required vanishing of
$g^{(j)}(0)$. Finally, $2\int_0^\infty g(u)\,du=H(0)=0$ and
$2\int_0^\infty H(u)\,du=g(0)=0$.\par}

\end{proof}

{
Let $W\in C_c^\infty((0,\infty))$ be supported sufficiently close
to $1$, and write $\Lambda_F$ for the von Mangoldt function on
ideals. Set}
\begin{equation}\label{eqc:detector-definition}
 K=\frac{\pi X}{2},\qquad
 P_X(k)=\sum_{\mathfrak a\ne0}
 \frac{\Lambda_F(\mathfrak a)\xi_k(\mathfrak a)}{\sqrt{N\mathfrak a}}
 W\left(\frac{\sqrt{N\mathfrak a}}X\right),\qquad
 I_X=\sum_{k\ge1}g(k/K)P_X(k).
\end{equation}

\begin{lemma}\label{leqc:detector}
For all sufficiently large $X$,
\begin{equation}\label{eqc:exact-detector}
 I_X=K\sum_{n\ge1}\frac{\Lambda(n^2+1)}{\sqrt{n^2+1}}
 W\left(\frac{\sqrt{n^2+1}}X\right)
 H\left(X\arctan\frac1n\right).
\end{equation}
\end{lemma}
\begin{proof}
Conjugation gives $P_X(-k)=P_X(k)$. Since $g$ is even and $g(0)=0$, the sum defining $I_X$ is half the corresponding sum over $k\in\mathbb Z$. If $\vartheta$ is a generator angle of $\mathfrak a$, Poisson summation gives
\begin{equation}\label{eqc:Poisson}
 \sum_{k\in\mathbb Z}g(k/K)e^{4ik\vartheta}
 =K\sum_{j\in\mathbb Z}H\left(K\left(j-\frac{2\vartheta}{\pi}\right)\right).
\end{equation}
{
Choose a generator $a+bi$ with angle $\vartheta\in[-\pi/4,\pi/4)$.
For large $X$, only $j=0$ can contribute to~\eqref{eqc:Poisson}.
Since $K=\pi X/2$, the support of $H$ forces
\[
 |X\vartheta|\in(1-\delta,1+\delta).
\]
Take the support of $W$ inside $(1-\delta,1+\delta)$. Then
\[
 |b|=\sqrt{a^2+b^2}\,|\sin\vartheta|
     =1+O(\delta+X^{-2}).
\]
For small fixed $\delta$ this gives $|b|=1$, with $a=n>0$.
The ideals $(n+i)$ and $(n-i)$ are distinct for $n>1$ and give
the same contribution by evenness of $H$. Their multiplicity
cancels the factor $1/2$ from extending the sum to all $k$.

If $n$ is even, any common prime ideal divisor of $(n+i)$ and
$(n-i)$ would divide $(2i)$, whereas $n^2+1$ is odd. Every rational
prime dividing $n^2+1$ then splits in $\mathbb Z[i]$, since $-1$
is a square modulo that prime. The coprimality shows that
$(n+i)$ contains at most one of the two prime ideals above each
such prime. Thus $(n+i)$ is a prime ideal power exactly when
$n^2+1$ is a rational prime power, and its ideal von Mangoldt
weight is $\Lambda(n^2+1)$. For odd $n>1$, we have $n^2+1\equiv2
\pmod 8$ and $n^2+1>2$, so both weights vanish. The radial
support excludes $n=1$ for large $X$. This proves
\eqref{eqc:exact-detector}.\par}

\end{proof}

{
For $\Phi$ supported sufficiently close to $1$, the function
\begin{equation}\label{eqc:weight-choice}
 W(u)=\frac{u\Phi(u)}{H(1/u)}
\end{equation}
is well defined on a neighborhood of its support, since $H(1)\ne0$.
We extend it by zero. For $n\asymp X$ we have
\[
 \frac{\sqrt{n^2+1}}X=\frac nX+O(X^{-2}),\qquad
 X\arctan(1/n)=\frac Xn+O(X^{-2}).
\]
Taylor's formula and~\eqref{eqc:weight-choice} give
\[
 \frac{1}{\sqrt{n^2+1}}W\left(\frac{\sqrt{n^2+1}}X\right)
 H\left(X\arctan\frac1n\right)
 =\frac{\Phi(n/X)}{X}+O_\Phi(X^{-3}).
\]
Only $O(X)$ integers occur, all with $n\asymp X$. Multiplying by
$K\Lambda(n^2+1)$ and summing gives
\begin{equation}\label{eqc:detector-approx}
 I_X=\frac KX\sum_{n\ge1}\Lambda(n^2+1)\Phi(n/X)
       +O_\Phi(X^{-1}\log X).
\end{equation}\par}

\subsection{The differentiated angular average}
{
Fix $0<\eta<\min\{\eta_*,\omega_H/2\}$ and use the line
$r=\eta+it$. For $Y\ge e^{2/\eta}$ and $|t|\le Y^{3/4}$ this
line lies in the shift range of Hypothesis~\ref{hyp:int-H}.
Under GRH, the local zero expansion gives
\begin{equation}\label{eqc:logder-bound}
 \frac{L_k'}{L_k}(u+iv)\ll_\eta \log^2(3+k+|v|)
 \qquad\bigl(u\ge\tfrac12+\eta,\ k\ge1\bigr).
\end{equation}
For $u\le2$, the zeros within one unit of height $v$ number
$O(\log(3+k+|v|))$, and the distance from each is at least
$\eta$. The logarithmic conductor term is also
$O(\log(3+k+|v|))$. For $u>2$ one can use the absolutely
convergent series for $L_k'/L_k$.

Write
\begin{equation}\label{eqc:BC}
 \mathcal B(r)=-\left.\partial_aG(a,r)\right|_{a=r},
 \qquad
 \mathcal C(r)=\zeta(1-2r)\frac{L(1+2r,\chi)}{L(1,\chi)}
                       A_1(-r,r)A_3(-r,r).
\end{equation}
The factor $1/\zeta(1-a+r)$ in $G(-a,r)$ has expansion
$-(a-r)+O((a-r)^2)$ at $a=r$. Hence
\[
 G(-a,r)=-(a-r)\mathcal C(r)+O_r((a-r)^2).
\]
Apply $-\partial_a$ at $a=r$ to the ratios formula in
Hypothesis~\ref{hyp:int-H}. Cauchy's formula on a fixed circle
about $r$, of radius at most $\min\{\rho,\eta/2\}$, preserves
the power of $Y$ in the error. The model is holomorphic on the
circle and its interior for a sufficiently small fixed radius,
by the local analyticity of the Euler products and of
$1/\zeta(s)$ and $1/L(s,\chi)$ near $s=1$. We obtain
\begin{equation}\label{eqc:differentiated}
 \sum_{1\le k\le Y}-\frac{L_k'}{L_k}(1/2+r)
 =Y\mathcal B(r)+\left(\frac\pi2\right)^{2r}
   \frac{Y^{1-2r}}{1-2r}\mathcal C(r)
   +O_{\eta,\varepsilon}(Y^{1-\omega_H+\varepsilon}).
\end{equation}

\begin{lemma}\label{leqc:weighted-derivative}
Fix $0<\tau<3/4$. For sufficiently large $J$ in
Lemma~\ref{leqc:window}, uniformly for $r=\eta+it$ and
$|t|\le K^\tau$,
\begin{equation}\label{eqc:weighted-derivative}
 \sum_{k\ge1}g(k/K)\left(-\frac{L_k'}{L_k}(1/2+r)\right)
 =K\left(\frac{\pi}{2K}\right)^{2r}\mathcal C(r)\widetilde g(1-2r)
   +O_{\eta,g,\varepsilon}(K^{1-\omega_H+\varepsilon}).
\end{equation}
\end{lemma}
\begin{proof}
By Abel summation the left side is
\[
 -\frac1K\int_0^\infty g'(Y/K)
          \sum_{1\le k\le Y}-\frac{L_k'}{L_k}(1/2+r)\,dY.
\]
For $Y\ge\max\{e^{2/\eta},|t|^{4/3}\}$ insert
\eqref{eqc:differentiated}. Its error contributes
\[
 \ll_{\eta,\varepsilon} K^{1-\omega_H+\varepsilon}
    \int_0^\infty |g'(u)|u^{1-\omega_H+\varepsilon}\,du.
\]
Below this range we use~\eqref{eqc:logder-bound}. Both
$\mathcal B(r)$ and $\mathcal C(r)$ are of at most polynomial
growth on the fixed line, by~\eqref{eqc:G} and~\eqref{eqc:BC}.
Since $g'(u)\ll_g u^{J-1}$ near zero, the integral of the actual
sum and of both model terms over the omitted range is bounded
in absolute value by
\[
 \ll_{\eta,g} K^{-J+C_\eta}
     \int_0^{\max\{e^{2/\eta},|t|^{4/3}\}}
              Y^{J-1}\bigl(Y+Y^{1-2\eta}\bigr)\,dY
 \ll_{\eta,g} K^{C_\eta+4\tau/3-J(1-4\tau/3)}
\]
for some fixed $C_\eta$. This is $O(K^{-1})$ for large enough
$J$. We may therefore integrate the two model terms over all
$Y>0$. Integration by parts gives $K\mathcal B(r)
\int_0^\infty g(u)\,du=0$ for the first and
\[
 \left(\frac\pi2\right)^{2r}\mathcal C(r)
    \int_0^\infty g(Y/K)Y^{-2r}\,dY
 =K\left(\frac{\pi}{2K}\right)^{2r}\mathcal C(r)\widetilde g(1-2r)
\]
for the second.
\end{proof}\par}

\subsection{The quadratic singular series}
{
The singular series is the first coefficient of the Dirichlet
series associated with the swapped term. Define}
\begin{equation}\label{eqc:D}
 \mathcal D(z)=\zeta(z)\frac{L(1+z,\chi)}{L(1,\chi)}
 A_1(-z/2,z/2)A_3(-z/2,z/2).
\end{equation}
For $\Re z>0$, local simplification yields
\begin{equation}\label{eqc:D-product}
 \mathcal D(z)=\mathfrak S_E\zeta(z)
 \prod_{p\equiv1\ (4)}\left(1+\frac{p^{-z}}{p-2}\right)
 \prod_{p\equiv3\ (4)}(1-p^{-1-z}).
\end{equation}
Indeed the respective local factors outside $\zeta(z)$ are
\begin{equation*}
 \frac{1-2/p+p^{-1-z}}{1-1/p}
 =\frac{p-2}{p-1}\left(1+\frac{p^{-z}}{p-2}\right),
 \qquad
 \frac{1-p^{-1-z}}{1-1/p}
 =\frac p{p-1}(1-p^{-1-z}).
\end{equation*}
The constant factors, ordered by increasing primes, give~\eqref{eq:quadratic-series}; equivalently one may factor out $L(1,\chi)^{-1}=4/\pi$ to obtain the absolutely convergent expression there. There is no extra local factor at $2$.

Both products in~\eqref{eqc:D-product} converge absolutely for $\Re z>0$. Thus $\mathcal D$ is meromorphic there with its only pole at $z=1$. For $\Re z>1$ it has the absolutely convergent expansion
\begin{equation}\label{eqc:D-series}
 \mathcal D(z)=\sum_{b\ge1}\frac{S_b}{b^z},\qquad
 S_b=\mathfrak S_E
 \prod_{\substack{p\mid b\\p\equiv1\ (4)}}\frac{p-1}{p-2}
 \prod_{\substack{p\mid b\\p\equiv3\ (4)}}\frac{p-1}{p}.
\end{equation}
{
In particular $S_1=\mathfrak S_E$. To identify $S_b$, use the
local density $(1-\nu_p/p)/(1-1/p)$, where $\nu_p$ counts the
roots of $n^2+b^2$ modulo $p$. For odd $p\nmid b$ we have
$\nu_p=1+\chi(p)$. For odd $p\mid b$ we have $\nu_p=1$, so the
density equals $1$. Dividing this by the density for $b=1$
gives $(p-1)/(p-2)$ at split primes and $(p-1)/p$ at inert
primes. At $2$ the density is $1$ for every $b$. Thus $S_b$
is the Bateman--Horn singular series for $n^2+b^2$.\par}

\subsection{Extracting the first coefficient}
Evenness and the Fourier convention give the Fourier--Mellin identity
\begin{equation}\label{eqc:Fourier-Mellin}
 \widetilde H(z)=2(2\pi)^{-z}\Gamma(z)
 \cos(\pi z/2)\widetilde g(1-z),\qquad 0<\Re z<1.
\end{equation}
For $0<\Re z<1$, insert
$H(u)=2\int_0^\infty g(t)\cos(2\pi tu)\,dt$ into the Mellin
integral. With a factor $e^{-\epsilon u}$, the inner integral is
\[
 \frac{\Gamma(z)}2\left((\epsilon-2\pi it)^{-z}
                        +(\epsilon+2\pi it)^{-z}\right).
\]
As $\epsilon\to0^+$ it tends to
$\Gamma(z)(2\pi t)^{-z}\cos(\pi z/2)$, and the bound
$O_z(t^{-\Re z})$ allows passage to the limit against $g(t)$.
This gives~\eqref{eqc:Fourier-Mellin}. We will also use
$\widetilde H(1)=0$, by~\eqref{eqc:vanishing}.

{
Take $\Phi$ supported sufficiently close to $1$, and choose $W$
by~\eqref{eqc:weight-choice}. Mellin inversion first gives the
following integral on a line $\Re z>1$. GRH and
\eqref{eqc:logder-bound} allow us to shift it to $\Re z=2\eta$.
\begin{equation}\label{eqc:Mellin-P}
 P_X(k)=\frac1{2\pi i}\int_{(2\eta)}
    -\frac{L_k'}{L_k}\left(\frac{1+z}{2}\right)
       X^z\widetilde W(z)\,dz.
\end{equation}
Sum against $g(k/K)$. For fixed $0<\tau<3/4$ the part with
$|\Im z|>K^\tau$ is $O_{W,g,A}(X^{-A})$ for each fixed $A$,
by~\eqref{eqc:logder-bound} and repeated integration by parts
in $\widetilde W$. On the remaining range apply
Lemma~\ref{leqc:weighted-derivative} with $r=z/2$. Its error
contributes
\begin{equation}\label{eqc:Mellin-error}
 \ll_{\eta,W,g,\varepsilon} X^{2\eta}K^{1-\omega_H+\varepsilon}
           \int_{\mathbb R}|\widetilde W(2\eta+it)|\,dt
 \ll_{\eta,W,g,\varepsilon} X^{1-\omega_H+2\eta+\varepsilon}.
\end{equation}
In the main term, the factor $X^z(\pi/(2K))^z$ equals $1$.
The zeta functional equation
$\zeta(1-z)=2(2\pi)^{-z}\Gamma(z)\cos(\pi z/2)\zeta(z)$
and~\eqref{eqc:Fourier-Mellin} now give
\begin{equation}\label{eqc:FE}
 \mathcal C(z/2)\widetilde g(1-z)
 =\mathcal D(z)\widetilde H(z).
\end{equation}
Rapid Mellin decay and polynomial growth of $\mathcal D$
allow us to extend the main term to the full line. Thus
\begin{equation}\label{eqc:main-integral}
 I_X=\frac K{2\pi i}\int_{(2\eta)}
 \widetilde W(z)\widetilde H(z)\mathcal D(z)\,dz
 +O_{\eta,\Phi,\varepsilon}(X^{1-\omega_H+2\eta+\varepsilon}).
\end{equation}
The factor $\widetilde H(1)=0$ cancels the possible pole at
$z=1$. Apply Lemma~\ref{co:extraction} with
\eqref{eqc:D-series} to obtain
\begin{equation}\label{eqc:coefficient-extraction}
 I_X=K\sum_{b\ge1}S_b\int_0^\infty W(u)H(b/u)\,\frac{du}{u}
       +O_{\eta,\Phi,\varepsilon}(X^{1-\omega_H+2\eta+\varepsilon}).
\end{equation}
The product of the support of $W$ and the positive support of
$H$ contains only one positive integer, namely $1$.
By~\eqref{eqc:weight-choice}, the integral for $b=1$ is
$\int_0^\infty\Phi(u)\,du$.

Combining this with~\eqref{eqc:detector-approx} and $K\asymp X$
gives
\begin{equation}\label{eqc:quadratic-smooth}
 \sum_{n\ge1}\Lambda(n^2+1)\Phi(n/X)
 =\mathfrak S_E X\int_0^\infty\Phi(u)\,du
   +O_{\eta,\Phi,\varepsilon}(X^{1-\omega_H+2\eta+\varepsilon}).
\end{equation}
For general $\Phi\in C_c^\infty((0,\infty))$, take a fixed
smooth partition of unity on its support into intervals of
sufficiently small relative width. On an interval centered at
$u_0>0$, use the scale $u_0X$ and rescale that piece of the
weight by $u\mapsto u_0u$. Summing over the pieces gives
\eqref{eqc:quadratic-smooth}. On a fixed compact support, the
estimates use only finitely many derivatives of $\Phi$.

Proposition~\ref{qt:sharp}, with $\eta$ sufficiently small in
terms of $\varepsilon$, gives~\eqref{eq:main-E}. This proves
Theorem~\ref{thm:quadratic-main}. Lebesgue's
theorem~\cite{Lebesgue1850} rules out proper prime powers of
the form $n^2+1$ with $n\ge1$. Partial summation with
$1/\log(t^2+1)$ therefore gives
\[
 \#\{n\le X:n^2+1\text{ is prime}\}
 =\mathfrak S_E\int_1^X\frac{dt}{\log(t^2+1)}
   +O_\varepsilon(X^{1-\omega_H+\varepsilon}).
\]\par}

\section{Binary prime problems over function fields}\label{ff:section}

In this section we prove Theorems~\ref{ff:reflection} and~\ref{ff:twins},
following the proof sketch given in Section~\ref{SKETCH-Fpt}.

\subsection{Uniform notation for the proofs}

For a monic Goldbach target~$N$ of degree~$d\ge3$, let
\begin{align*}
 S_{-1}&=\sum_{F\in\mathcal M_{d-1}}\Lambda(F)\Lambda(N-F)\,,&
 S_{-1}^\times&=\sum_{\substack{F\in\mathcal M_{d-1}\\(F,N)=1}}
                 \Lambda(F)\Lambda(N-F)\,.
\end{align*}
For each nonzero $h\in \F_q[T]$ and integer $n>\deg h$, let
\begin{align*}
 S_1&=\sum_{F\in\mathcal M_n}\Lambda(F)\Lambda(F+h)\,,&
 S_1^\times&=\sum_{\substack{F\in\mathcal M_n\\(F,h)=1}}
                 \Lambda(F)\Lambda(F+h)\,.
\end{align*}
Given $\sigma\in \{-1,1\}$ and a polynomial $F$ appearing in the sum $S_\sigma$, write
$$G=D+\sigma F.$$
For $\sigma=-1$, let $D=N$, $n=d-1$ and~$m=d$.
For $\sigma=1$, let $D=h$ and~$m=n$.
Then in either case, we have $F\in \mathcal M_n$, $G\in \mathcal M_m$, and
\begin{align*}
 S_\sigma&=\sum_{F\in\mathcal M_n}\Lambda(F)\Lambda(G)\,,&
 S_\sigma^\times&=\sum_{\substack{F\in\mathcal M_n\\(F,D)=1}}
                 \Lambda(F)\Lambda(G)\,.
\end{align*}

\subsection{Short divisors and the main term}
For each $\sigma\in \{-1,1\}$, it suffices to estimate~$S_\sigma^\times$,
since $S_\sigma - S_\sigma^\times \ll_{q,\varepsilon} q^{\varepsilon n}$
by the divisor bound.
Expanding $\Lambda(G)=-\sum_{C\mid G}\mu(C)\deg C$ gives
\begin{equation}\label{ff:direct-divisors}
 S_\sigma^\times=-\sum_{j=1}^m j
 \sum_{\substack{C\in\mathcal M_j\\(C,D)=1}}\mu(C)
 \sum_{\substack{F\in\mathcal M_n\\ F\equiv-\sigma D\bmod C\\(F,D)=1}}\Lambda(F),
\end{equation}
since $G = D+\sigma F$.
Here the condition $(C,D)=1$ follows from $(F,D)=1$.
Put $K=\lfloor5n/6\rfloor$ and let $S_{\sigma,K}^\times$ be the contribution to $S_\sigma^\times$ from $j\le K$. For each $j\le K$ and $C\in\mathcal M_j$ with $(C,D)=1$ and $\mu(C)\ne 0$, Sawin's estimate~\cite[Proposition~8.16 at level $\theta=5/6$]{Sawin} gives
\begin{equation}\label{ff:progressionbound}
 \sum_{\substack{F\in\mathcal M_n\\F\equiv -\sigma D\bmod C}}\Lambda(F)
 =\frac{q^n}{\Phi(C)}+O_q(q^{(n-j)(1-\delta_q)}),
 \qquad \delta_q=\frac12-\log_q(352\mathrm e),
\end{equation}
where $\Phi(C)=\#(\F_q[T]/C)^\times$ and $\mathrm e=\exp(1)$.
The restriction of the left-hand side of \eqref{ff:progressionbound} 
to $F$ with $(F,D)\ne 1$
is $O_{q,\varepsilon}(q^{\varepsilon n})$ by the divisor bound,
which is negligible.
Subtracting off this contribution
and summing over $C\in\mathcal M_j$ for each $j\le K$ gives
\begin{equation}\label{ff:direct-AP-error}
S_{\sigma,K}^\times
+ \sum_{j=1}^K j
 \sum_{\substack{C\in\mathcal M_j\\(C,D)=1}}\mu(C) \frac{q^n}{\Phi(C)}
\ll_q m^2 q^{n-\delta_q(n-K)}
   \ll_q m^2 (352\mathrm e)^{n/6}q^{11n/12}.
\end{equation}

The main term of \eqref{ff:direct-AP-error}
was already computed in \cite[Proposition 6.2]{SSAnnals}:
\begin{equation*}
\sum_{j=1}^K j
 \sum_{\substack{C\in\mathcal M_j\\(C,D)=1}} \frac{\mu(C)}{\Phi(C)}
 = -\mathfrak S_q(D)
 + O_{q,\varepsilon}((q-1)^{\varepsilon(K+\deg D) - K}).
\end{equation*}
Multiplying this by $q^n$ and plugging it into \eqref{ff:direct-AP-error} gives
\begin{equation}\label{ff:direct-short}
S_{\sigma,K}^\times
- q^n \mathfrak S_q(D)
   \ll_q m^2 (352\mathrm e)^{n/6}q^{11n/12}
   + q^{n - K/2},
\end{equation}
since $\deg D \le n+1 \ll K$.
It remains to bound the quantity
\begin{equation*}
S_\sigma^\times - S_{\sigma,K}^\times=-\sum_{K<j\le m}j
 \sum_{\substack{C\in\mathcal M_j\\(C,D)=1}}\mu(C)
 \sum_{\substack{F\in\mathcal M_n\\ F\equiv-\sigma D\bmod C\\(F,D)=1}}\Lambda(F).
\end{equation*}

\subsection{Long divisors and their degrees of freedom}
For $j>K$, write $G=D+\sigma F=BC$, with $B\in\mathcal M_{m-j}$.
Moreover, by divisor switching, we may write
\begin{equation*}
\begin{split}
 \Lambda(F)
 &=-
 \sum_{Q\mid F}\mu(Q)\deg Q \\
 &=-
 \sum_{\substack{Q\mid F\\ \deg Q\le \lfloor n/2\rfloor}}\mu(Q)\deg Q
 -
 \sum_{\substack{Q\mid F\\ \deg Q<n-\lfloor n/2\rfloor}}\mu(F/Q)(n-\deg Q).
\end{split}
\end{equation*}
Writing $F=QZ$, where $Q\in\mathcal M_k$ with $k\le n/2$, it follows that
\begin{equation}
\label{ff:direct-hyperbola}
S_\sigma^\times - S_{\sigma,K}^\times
= S_{\sigma,K,1}^\times
+ S_{\sigma,K,2}^\times,
\end{equation}
where
\begin{equation*}
S_{\sigma,K,1}^\times
= \sum_{K<j\le m}j
 \sum_{C\in\mathcal M_j}\mu(C)
 \sum_{B\in \mathcal M_{m-j}}
 \sum_{k\le n/2} \sum_{Q\in \mathcal M_k}
 \sum_{\substack{Z\in\mathcal M_{n-k}\\ BC=D+\sigma QZ\\(QZ,D)=1}}
 \mu(Q) k
\end{equation*}
and
\begin{equation*}
S_{\sigma,K,2}^\times
= \sum_{K<j\le m}j
 \sum_{C\in\mathcal M_j}\mu(C)
 \sum_{B\in \mathcal M_{m-j}}
 \sum_{k\le (n-1)/2} \sum_{Q\in \mathcal M_k}
 \sum_{\substack{Z\in\mathcal M_{n-k}\\ BC=D+\sigma QZ\\(QZ,D)=1}}
 \mu(Z) (n-k),
\end{equation*}
since $n-\lfloor n/2\rfloor - 1 = \lfloor (n-1)/2\rfloor$.
Here we have dropped the condition $(C,D)=1$,
which we may do because the conditions $BC=D+\sigma QZ$ and $(QZ,D)=1$
already imply $(BC,D)=1$.

The strategy is now to fix $B$ and $Q$,
and to seek cancellation over $C$ and $Z$
in each sum $S_{\sigma,K,i}^\times$.
For this it will be helpful to parameterize solutions to $BC=D+\sigma QZ$.
Clearly $(B,Q)=1$, since $(B,Q)\mid (D,Q)=1$.
Given $B$ and $Q$ with $(B,Q)=1$,
there is a unique $u\in \F_q[T]$ with $\deg u<k$ such that $Bu\equiv D\bmod Q$.
Then $C\equiv u\bmod{Q}$, so there exists $V\in \mathcal M_{j-k}$ for which
\begin{equation}\label{ff:direct-param}
 C=u+QV,\qquad Z=\sigma(Bu-D)/Q+\sigma BV .
\end{equation}

Conversely, for any $V\in \mathcal M_{j-k}$,
the polynomials defined by \eqref{ff:direct-param}
satisfy $C\in \mathcal M_j$ and $Z\in \mathcal M_{n-k}$ automatically if $\sigma=1$,
since $\deg Bu/Q < m-j\le m-k$ and $\deg D/Q < n-k = m-k = \deg BV$.
If $\sigma=-1$, then $C\in \mathcal M_j$ is automatic
and $\deg D = m=d=n+1$, so the leading terms of $(Bu-D)/Q$ and $BV$
are $-T^{n+1-k}$ and $T^{n+1-k}$, respectively.
Thus the condition $Z\in \mathcal M_{n-k}$ cuts out an additional linear equation
in the coefficients of $V$ if $\sigma=-1$.
In any case, there exists $V_0\in \mathcal M_{j-k}$ depending only on $B,Q,\sigma,D$,
such that $V=V_0+w$, where $w\in \F_q[T]$ runs through all polynomials of degree less than $L=j-k-(m-n)$.
If $\sigma=1$, then we can simply take $V_0 = T^{j-k}$;
if $\sigma=-1$, then $V_0$ must satisfy the aforementioned linear constraint on its coefficients.
Observe that
\begin{equation}\label{ff:direct-length}
 j-k\ge L\ge n/3-O(1).
\end{equation}

\subsection{Uniformity in the Das--Howe sieve under degeneration}

To estimate the M\"obius sums over $C$ and $Z$ above, we will want a suitable homological vanishing result in a linear range, in view of the Lefschetz trace formula.
To this end, we first prove a useful degeneration result in \'etale cohomology, using a homological inclusion--exclusion theorem of Das and Howe \cite[Theorem~6.4.1(ii)]{DasHowe}.

Let~$k$ be algebraically closed, let~$V=\mathbb A_k^b$ and~$X=\mathbb A_k^N\times_k V$, and let~$Z\to Y\hookrightarrow X$ be a finite surjective map followed by a closed immersion. Put~$W=X\setminus Y$, and write~$D_0=D\times_V\{0\}$ for every~$D\to V$. Let~$m(x)$ be the number of distinct geometric points of~$Z_{\bar x}$, for any geometric point $\bar x\to x\in |X|$. The ordered and unordered incidence schemes are
\begin{equation*}
 \widetilde I_j=Z_X^j\setminus\bigcup_{1\le a<b\le j}\Delta_{ab},
 \qquad I_j=\widetilde I_j/\mathfrak S_j\,.
\end{equation*}
Here~$Z_X^j$ is the~$j$-fold fiber product of $Z$ over $X$,
and~$\Delta_{ab}$ is the diagonal where coordinates~$a,b$ agree.
The ordering cover~$\widetilde I_j\to I_j$ is a finite \'etale~$\mathfrak S_j$-torsor;
denote its sign local system by~$o_j$, so that $o_j$ is trivialized by the sign double cover of $I_j$.

\begin{proposition}\label{ff:configcompare}
Fix~$K\geq1$. Suppose that, for each $1\leq j\leq K$,
there exists a $\mathfrak S_j$-equivariant isomorphism of reduced schemes
\begin{equation}\label{ff:incidence-product}
 (\widetilde I_j)_{red} \cong ((\widetilde I_j)_0\times_k V)_{red}
\end{equation}
over~$V$.
Suppose also that the set $B=\{x\in|X|:m(x)>K\}$ satisfies
\begin{equation}\label{ff:incidence-codim}
 \dim X - \dim\overline B\geq K+1,
 \qquad\dim{X_0} - \dim\overline{B\cap|X_0|} \geq K+1\,,
\end{equation}
where the closures $\overline B$ and $\overline{B\cap|X_0|}$
are taken in $X$ and $X_0$, respectively.
For~$\Lambda=\F_\ell$ or~$\Q_\ell$, with~$\ell$ invertible in~$k$, restriction induces
an isomorphism
\begin{equation*}
 H^a(W,\Lambda)\xrightarrow{\ \sim\ } H^a(W_0,\Lambda)\qquad(0\leq a\leq K)\,.
\end{equation*}
If~$\operatorname{char}k\neq2$ and~$\mathcal E$ is a rank-one~$\Z_2$-local system on~$W$ with monodromy in~$\{1,-1\}$, restriction with coefficients~$\mathcal E$ is an isomorphism for~$a<K$ and injective for~$a=K$. The same holds with coefficients~$\mathcal E\otimes\Q_2$.
\end{proposition}

\begin{proof}
First take constant field coefficients~$\Lambda$.
Let $P=\Gamma_Y^{\le K}(Z)$ be the relative poscheme of nonzero
effective zero-cycles of degree at most $K$ in the fibers of $Z\to Y$,
ordered by inclusion. Its nerve $NP$ parametrizes chains of such cycles.
Let $\epsilon_Y:NP\to Y$ send each chain to its base point, and put
$\mathcal R_{K,Y}=R(\epsilon_Y)_*\Lambda$
\cite[Example~6.2.5 and Theorem~6.4.1(ii)]{DasHowe}.

Filter $NP$ by the degree of the largest cycle in a chain.
For $p_j:I_j\to Y$ the structure map, the induced filtration satisfies
\[
  \operatorname{Gr}^{j-1}\mathcal R_{K,Y}
  \cong R(p_j)_!o_j[1-j] \qquad (1\le j\le K)
\]
by \cite[Lemma~4.3.3 and Proposition~6.3.2]{DasHowe}.
% For $\Lambda=\F_\ell$, we apply these results to the constructible
% $\Z_\ell$-module $\Z_\ell/\ell$.
% 
% First take constant field coefficients~$\Lambda$.
% Let~$\mathcal R_{K,Y}=R(\epsilon_Y)_*\Lambda$, where~$\epsilon_Y:N\Gamma_Y^{\leq K}(Z)\to Y$ is the augmentation of Das--Howe's truncated effective-cycle construction.
% Its functorial cycle-degree filtration has graded pieces~$(I_j\to Y)_!o_j[1-j]$, $1\leq j\leq K$~\cite[\S4.3 and Proposition~6.3.2]{DasHowe}.
Moreover, properness of the map $Z\to Y$ lets us identify the restriction $\mathcal R_{K,Y}\vert_{Y_0}$ with the corresponding construction $\mathcal R_{K,Y_0}=R(\epsilon_{Y_0})_*\Lambda$ on $Y_0$.
% Relative divided powers commute with base change, so proper base change
% on each level of the nerve identifies the restriction of the filtered
% complex $\mathcal R_{K,Y}$ to $Y_0$ with $\mathcal R_{K,Y_0}$.
% For $i:Y_0\hookrightarrow Y$, pulling the supported class of
% $0\hookrightarrow V$ back to $Y$ and tensoring with each filtered
% term $\mathcal F$ gives $i_*i^*\mathcal F(-b)[-2b]\to\mathcal F$.
The Gysin map associated to the regular immersion $0\to V$ is natural
with respect to the aforementioned filtration.
On graded pieces $\operatorname{Gr}^{j-1}\mathcal R_{K,Y}$ it induces compact-support cohomology isomorphisms by~\eqref{ff:incidence-product} and the K\"unneth formula; induction on the degree $j$ in the filtration then gives compact-support cohomology isomorphisms on~$\mathcal R_{K,Y}$ itself.
Since~$Z\to Y$ is finite, the auxiliary dimension hypothesis of~\cite[Theorem~6.4.1(ii)]{DasHowe} is automatic, and the augmentation map $\epsilon_Y$ identifies~$H_c^q(Y,\Lambda)$ with~$H_c^q(Y,\mathcal R_{K,Y})$ for~$q\geq K+2\dim\overline B+1$.
Similarly, $\epsilon_{Y_0}$ identifies $H_c^{q-2b}(Y_0,\Lambda)$ with~$H_c^{q-2b}(Y_0,\mathcal R_{K,Y_0})$ for~$q-2b\geq K+2\dim\overline{B\cap|X_0|}+1$. The two bounds in~\eqref{ff:incidence-codim} therefore give
\begin{equation*}
 H_c^{q-2b}(Y_0,\Lambda)(-b)\xrightarrow{\ \sim\ } H_c^q(Y,\Lambda)
 \qquad(q\geq2N+2b-K-1)\,.
\end{equation*}
The localization exact sequence in compactly supported cohomology for~$Y\subset X$, together with the fact that $X=X_0\times_k V$, now gives the Gysin isomorphism~$H_c^{h-2b}(W_0,\Lambda)(-b)\xrightarrow{\ \sim\ } H_c^h(W,\Lambda)$ for~$h\geq2N+2b-K$. Poincar\'e duality gives the asserted restriction isomorphism.
% ~\cite[\S\S15--19,22--24]{Milne}.

For quadratic coefficients, put~$\mathcal E_m=\mathcal E/2^m\mathcal E$.
Using the short exact sequence of coefficients
$0\to\mathcal E_{m-1}\xrightarrow{2}\mathcal E_m\to\F_2\to0$,
and the isomorphism $\mathcal E_1\cong\F_2$, we may induct on $m\ge 1$
using the five lemma and four lemma
to obtain restriction isomorphisms $H^a(W,\mathcal E_m)\to H^a(W_0,\mathcal E_m\vert_{W_0})$ for~$a<K$ and injectivity for~$a=K$. Finiteness of these \'etale cohomology groups permits passage to the inverse limit,
% ~\cite[\S19]{Milne},
proving the~$\Z_2$-assertion; inverting~$2$ gives the~$\Q_2$-assertion.
\end{proof}

\subsection{Cancellation on affine families}\label{ff:geometry-appendix}

Recall from \eqref{ff:direct-param} and the subsequent discussion that
$C$ and $Z$ are affine functions of a polynomial $w\in \F_q[T]_{<L}$.
Moreover, we are interested in $\mu(C) \mu(Z) = \mu(CZ)$, since $(C,Z)=1$.
Thus we are interested in the locus of polynomials $w\in \F_q[T]_{<L}$
for which $CZ$ is square-free, and moreover we seek cancellation over $w$.
Via the Lefschetz trace formula, this setting will all fit into
the following geometric framework.

Let~$k$ be algebraically closed of characteristic different from two, and fix~$0\ne D\in k[T]$. For~$1\le i\le s$, with~$s=1$ or~$2$, put
\begin{equation}\label{ff:families}
 f_i=a_iw+b_i,\qquad w\in k[T]_{<L},\qquad a_i\ne0,\qquad
 n_i=\deg b_i\ge\deg a_i+L\,.
\end{equation}
Here~$k[T]_{<L}$ is the affine space of polynomials of degree less than~$L$.
Assume $(a_1,b_1)=1$ if $s=1$, and if $s=2$ then assume
\begin{equation}\label{ff:determinant}
 (a_1,a_2)=1,\qquad (a_1a_2,D)=1,\qquad a_1b_2-a_2b_1=D\,.
\end{equation}
These conditions automatically imply~$(a_i,b_i)=1$ if $s=2$. Let~$U\subset\mathbb A^L$ be the open set where the~$f_i$ are squarefree and prime to~$D$, and let~$\mathcal L$ be the (nontrivial) rank-one local system over $\Q_2$ given by the square root of~$\prod_i\Disc(f_i)$ on~$U$.

\begin{proposition}\label{ff:vanishing}
For the families~\eqref{ff:families}, put~$\kappa(L)=\max(0,\lfloor L/2\rfloor-5)$. For every~$L\ge1$,
\begin{equation}\label{ff:lowvanishing}
 H^j(U,\mathcal L)=0\quad(j<\kappa(L)),\qquad
 H_c^j(U,\mathcal L)=0\quad(j>2L-\kappa(L))\,.
\end{equation}
\end{proposition}

\begin{proof}
The two assertions are dual since~$U$ is smooth and~$\mathcal L$ is self-dual. If~$L<12$ then~$\kappa(L)=0$ and the assertions hold by Artin's vanishing theorem. Assume~$L\ge12$ and set~$K=\kappa(L)$. Write~$r=\sum_i n_i$, $V=k[T]_{<r}$, and~$X=\mathbb A^L_w\times V$. For~$v\in V$ set~$g_i=f_i'+a_iv$, where the prime $'$ in $f_i'$ denotes differentiation in~$T$. On~$X$ define the polynomial
\begin{equation}\label{ff:norm-deformation}
 R= \prod_{1\le i\le s} \Res(f_i,g_i)\,,
\end{equation}
a product of resultant polynomials.
Let~$W\subset X$ be the open set where~$R\ne0$, and all~$f_i$ are prime to~$D$. On~$W$ let~$\mathcal E$ be the (nontrivial) quadratic $\Z_2$-valued local system associated to the square root of $R$, and write~$\mathcal E_{\Q_2}=\mathcal E\otimes\Q_2$. At~$v=0$ the resultant--discriminant identity gives
\begin{equation}\label{ff:zero-fiber}
 W\vert_{v=0}=U,\qquad \mathcal E_{\Q_2}|_U\cong\mathcal L\,.
\end{equation}

We first bound the compactly supported cohomology of~$W$ by projecting to~$w$. Over a geometric point of $\mathbb A^L_w$ with residue field~$k^\prime$ and with each~$f_i$ prime to~$D$, the~$a_i$ are units modulo~$f_i$, and in the two-polynomial case $s=2$, the~$f_i$ are coprime by~\eqref{ff:determinant}. The Chinese remainder theorem gives an affine linear isomorphism $V\to \prod_i k^\prime[T]/(f_i)$, $v\mapsto (g_i\bmod f_i)_i$. The fiber of~$W$ over $w$ is isomorphic to the unit group of $\prod_i k^\prime[T]/(f_i)$.
By the Chinese remainder theorem
and the multiplicative properties of resultants, a root $x$ of $f_i$ with multiplicity~$m_x$ contributes a factor~$\mathbb G_m\times\mathbb A^{m_x-1}$ on which the resultant $\Res((T-x)^{m_x},g_i)$ equals~$g_i(x)^{m_x}$ up to a fixed sign independent of $x$ and $g_i$.
Here $g_i(x)\ne 0$ is the coordinate of the $\mathbb G_m$ factor.
However, it is well known that $H^\ast_c(\mathbb G_m,\chi_2) = 0$, where $\chi_2$ is the (nontrivial) quadratic $\Q_2$-valued local system on $\mathbb G_m$.
If any of the integers $m_x$ is odd, then it follows from the K\"unneth formula
that $H^\ast_c(W_w,\mathcal E_{\Q_2}\vert_{W_w}) = 0$.
Therefore, $H^\ast_c(W_w,\mathcal E_{\Q_2}\vert_{W_w}) = 0$ unless
every~$f_i$ is a square.

The exceptional locus $\mathcal B \subset \mathbb A^L_w$, consisting of points $w$ for which every~$f_i$ is a square,
has dimension at most~$\lfloor L/2\rfloor$. Indeed it is empty for odd~$n_1$. For even~$n_1$ write~$f_1=c h^2$ with~$h$ monic and~$c$ its fixed leading coefficient. Squaring on monic polynomials is a closed immersion since the top coefficients of the square recover those of~$h$ recursively by division by two. Every tangent vector to the inverse image of the square locus thus satisfies~$a_1\dot w=2ch\dot h$. Since~$(a_1,h)=1$, this forces~$h\mid\dot w$, leaving at most~$\max(0,L-n_1/2)\le\lfloor L/2\rfloor$ tangent dimensions. By proper base change, the cohomology sheaves of the $w$-coordinate projection from $W$ are supported on the locus $\mathcal B$ and vanish in degrees above~$2r$. The compact-support Leray spectral sequence gives
\begin{equation}\label{ff:generic-norm}
 H_c^i(W,\mathcal E_{\Q_2})=0\qquad(i>2r+2\lfloor L/2\rfloor)\,.
\end{equation}

To transfer vanishing to $U$,
we wish to apply Proposition~\ref{ff:configcompare}.
Thus consider the closed incidence correspondence $Z \subset X\times\mathbb A^1$ marking a location~$x\in \mathbb A^1$ where either~$f_i(x)=g_i(x)=0$ for some~$i$, or~$D(x)=f_1(x)=0$.
Note that if $D(x)=f_1(x)=0$ and $s=2$, then $f_2(x)=0$ by \eqref{ff:determinant}.
The locus $Z$ is clearly finite over~$X$,
and its image in $X$ is $X\setminus W$ by the definition of $W$.

For an ordered configuration of~$j$ distinct marked locations~$(x_1,\ldots,x_j)$ in $Z$, put~$H(T)=\prod_{\nu=1}^j(T-x_\nu)$, and set
\begin{equation}\label{ff:hermite}
 h_H=\left(\frac{v}{H'}\bmod H\right) \, H,\qquad \deg h_H<2j\,.
\end{equation}
Here the remainder $\frac{v}{H'}\bmod H$ is viewed as a polynomial of degree less than~$j$.
We have $h_H(x)=0$ and~$h_H'(x)=v(x)$ at each mark $x=x_\nu$.
For~$2j\le L$ the translation~$w_0=w+h_H$ preserves~$f_i(x)$ and gives~$f_i'(w_0)(x)=g_i(w)(x)$ there. Thus $(\widetilde I_j)_{red} \cong ((\widetilde I_j)\vert_{v=0}\times_k V)_{red}$ over~$V$.

On the zero section $(\widetilde I_j)\vert_{v=0}$, fix~$s'$ marks outside~$D$ and~$t'$ distinct marked roots of~$D$, with~$s'+t'=j$. At a mark~$x$ outside $D$, the index~$i$ for which $f_i(x)=g_i(x)=0$ is unique by~\eqref{ff:determinant} if $s=2$.
Since~$(a_i,b_i)=1$, one has~$a_i(x)\ne0$, so the two equations prescribe~$w(x)$ and~$w'(x)$.
At each marked root~$y$ of $D$, the condition $f_1(y)=0$ prescribes~$w(y)$ since~$a_1(y)\ne0$. These~$2s'+t'$ conditions are independent whenever~$2s'+t'\le L$, by polynomial interpolation modulo~$\prod_x(T-x)^2\prod_y(T-y)$.
Allowing the marks outside $D$ to vary adds at most~$s'$ dimensions,
whereas the choices of indices, and marked roots of $D$, are finite.
Thus
$$\dim(\widetilde I_j)\vert_{v=0}\le L-(2s'+t')+s' = L-j$$
for~$2j\le L$. Since~$2(K+1)\le L$, taking $j=K+1$ shows that if $B=\{x\in X:\#Z_{\bar x}>K\}$, then
\begin{equation}\label{ff:manymarks}
 \operatorname{codim}_X\overline B\ge K+1,\qquad
 \operatorname{codim}_{X_0}\overline{B\cap X_0}\ge K+1\,,
\end{equation}
where each closure is taken in its indicated ambient space.
Proposition~\ref{ff:configcompare} now gives
\begin{equation}\label{ff:signcompare}
 H^a(W,\mathcal E_{\Q_2})\xrightarrow{\ \sim\ }H^a(U,\mathcal L)
 \qquad(a<K)\,.
\end{equation}
By duality on the smooth open~$W$ of dimension~$L+r$, \eqref{ff:generic-norm} implies that the left side vanishes for~$a<2L-2\lfloor L/2\rfloor$, which includes~$a<K$. This proves~\eqref{ff:lowvanishing}.
\end{proof}

\subsection{From cohomology to the affine sum}
Suppose now that the data of~\eqref{ff:families} descend from $k=\overline{\F}_q$ to~$\F_q$, with~$q$ odd, and put~$e=\deg D$.
Proposition~\ref{ff:vanishing} gives a linear vanishing range for the homology of $U$,
but in the Lefschetz trace formula we still need to control the dimensions of the remaining homology groups.
Applying an affine Betti bound directly to the discriminants in \eqref{ff:norm-deformation} would not work, because the degrees of the discriminants grow with~$n_i$.
Instead we express the open conditions and the sign summand by equations of degrees at most four.
To do this, we use auxiliary polynomials to impose coprimality conditions, following Sawin and Shusterman~\cite[Appendix B, especially Lemma B.10]{SSBetti}.

\begin{proposition}\label{ff:affine-bound}
Under the hypotheses of Proposition~\ref{ff:vanishing}, with~$r=\sum_i n_i$ and~$e=\deg D$,
\begin{equation}\label{ff:affine-sum}
 \left|\sum_{\substack{w\in\F_q[T]_{<L}\\(f_i(w),D)=1\ \forall i}}
        \prod_i\mu(f_i(w))\right|
 \le 10^{L+2r+e+1}\,q^{3L/4+3}\,.
\end{equation}
\end{proposition}

\begin{proof}
% Take geometric cohomology over $k=\overline{\F}_q$.
% If $U$ is empty, the sum is zero.
For each $i$, squarefreeness of $f_i$ is equivalent to the existence of
polynomials $A_i,B_i$ satisfying
\[
 A_i f_i+B_i f_i'=1,\qquad \deg A_i<n_i-1,\quad \deg B_i<n_i.
\]
Indeed $B_i$ is uniquely determined as the inverse of $f_i'$ modulo
$f_i$, and the identity then determines $A_i$. Their coefficients are
regular on the squarefree locus by Cramer's rule. If $e>0$, impose also
$f_1E\equiv1\bmod D$ with $\deg E<e$. This inverse is again unique
and regular, and in the two-polynomial case~\eqref{ff:determinant}
makes the two coprimality conditions equivalent. The identity for $f_i$
gives $2n_i-1$ coefficient equations and the congruence gives $e$.
They have degree at most two in $w$ and the auxiliary coefficients.
Thus $U$ is a smooth complete intersection of codimension $c=2r-s+e$
in $\mathbb A^M$, where $s$ is the number of polynomials and $M=L+c$.

Choose successive general affine hyperplane sections
$U=Y_0\supset Y_1\supset\cdots\supset Y_L$ in this presentation.
Write $\mathcal L_j=\mathcal L|_{Y_j}$. On the discriminant double
cover of $Y_j$ the deck involution acts freely and has order prime to
$\operatorname{char}k$. Its alternating trace on compactly supported
$\Q_2$-cohomology is therefore zero by
\cite[Theorem~3.2]{DeligneLusztig}. The two eigenspaces compute
cohomology with coefficients $\Q_2$ and $\mathcal L_j$, respectively.
Hence $\chi_c(Y_j,\mathcal L_j)=\chi_c(Y_j,\Q_2)$.

Put $b_j=\sum_i\dim H^i(Y_j,\mathcal L_j)$. Each nonempty $Y_j$ is
smooth affine of dimension $d=L-j$, so duality and self-duality of
$\mathcal L_j$ give the same total Betti number and Euler characteristic
with compact supports. For $d>0$, affine vanishing gives $H^i=0$ for
$i>d$, and weak Lefschetz applied to $\mathcal L_j[d]$ gives injections
into $H^i(Y_{j+1},\mathcal L_{j+1})$ for $i<d$
\cite[Theorem~2.3.2, arXiv version]{WanZhang}. Solving for the middle
Betti number in the Euler characteristic gives
\[
 b_j\le |\chi_c(Y_j,\Q_2)|+2b_{j+1},
 \qquad b_L=|\chi_c(Y_L,\Q_2)|.
\]
In its ambient linear space $\mathbb A^{M-j}$, the slice $Y_j$ has codimension
$c$ whenever nonempty and is cut out by $c$ equations of degree at most two.
Therefore \cite[Theorem~1.2.2, arXiv version]{WanZhang} gives
\[
 |\chi_c(Y_j,\Q_2)|
 \le \binom{M-j-1}{c-1}3^{M-j}\le \tfrac12\,6^{M-j}.
\]
Empty slices contribute zero. Iterating gives
$b_0\le \tfrac12\,6^M\sum_{j=0}^{L}3^{-j}<6^M$.

By Pellet's identity the sum in~\eqref{ff:affine-sum} is $(-1)^r$
times the trace sum of $\mathcal L$ on $U(\F_q)$. The sheaf has
weight zero, and Proposition~\ref{ff:vanishing} gives
$H_c^i(U,\mathcal L)=0$ for $i>2L-\kappa(L)$. The trace formula and
Deligne's weight bound~\cite{Deligne} thus bound the sum in absolute
value by $b_0 q^{L-\kappa(L)/2}$. Now
$L-\kappa(L)/2\le 3L/4+3$ and $6^M<10^{L+2r+e+1}$,
which proves~\eqref{ff:affine-sum}.
\end{proof}

\subsection{Estimating the long-divisor sum}

In the first sum of~\eqref{ff:direct-hyperbola} only $\mu(C)$ varies, and in the second both $\mu(C)$ and $\mu(Z)$ vary. We return to the families~\eqref{ff:direct-param}, with $V=V_0+w$, and check the hypotheses of Proposition~\ref{ff:affine-bound},
with $f_1=C=u+QV$ and $f_2=Z=\sigma(Bu-D)/Q+\sigma BV$.
The coefficient of $w$ in $f_1$ is $a_1=Q$, which is prime to $b_1=u+QV_0$ because $Bu\equiv D\bmod Q$.
Moreover $a_2=\sigma B$, $b_2=\sigma(Bu-D)/Q+\sigma BV_0$, and
\begin{equation*}
a_1b_2 - a_2b_1 =
 Q\bigl(\sigma(Bu-D)/Q+\sigma BV_0\bigr)-\sigma B(u+QV_0)=-\sigma D .
\end{equation*}
Moreover, $(a_1,a_2)=1$ and $(a_1a_2,D)=1$, so \eqref{ff:determinant} holds.
The degree inequalities
$\deg b_i \ge \deg a_i + L$ required in \eqref{ff:families} also hold,
because $j\ge k+L$ and $n-k=(m-j)+L$.
We can therefore apply~\eqref{ff:affine-sum} with $-\sigma D$ in place of $D$.

Observe that $\deg b_1+\deg b_2\le j+n-k\le2m$, and $\deg D\le m$, $L\le m$. Consequently~\eqref{ff:affine-sum} bounds the inner sum over $C,Z$ in $S_{\sigma,K,i}^\times$ by $O_q(10^{6m}q^{3L/4})$, after switching the order of summation so that $j,k,B,Q$ are fixed on the outside. At fixed $j,k$ there are at most $q^{m-j+k}=q^{n-L}$ choices of $B,Q$. Summing over $j,k,B,Q$ and then using~\eqref{ff:direct-length} gives
\begin{equation*}
S_{\sigma,K,i}^\times \ll_q m^4 10^{6m} \sup_{j,k} q^{n-L/4}
\ll_q m^4 10^{6m}q^{11n/12}.
\end{equation*}
Together with~\eqref{ff:direct-short} and \eqref{ff:direct-hyperbola} this yields
\begin{equation}\label{ff:weightedunit}
 S_\sigma^\times=q^n\mathfrak S_q(D)+O_q(m^4 10^{6m}q^{11m/12})
             =q^n\mathfrak S_q(D)+O_q(q^{23m/24}),
\end{equation}
since $11/12+6\log_q10<23/24$.
This completes the proof of Theorems~\ref{ff:reflection} and~\ref{ff:twins}.

\appendix

\section{Ratios with qualitative errors}\label{wk:section}
For prime differences and quadratic values, the power-saving ratios errors in the main theorems can be weakened if only qualitative asymptotics are required. The condition $\omega>2\eta$ compensates for the fixed Mellin-line factor $X^{2\eta}$. Taking $\eta$ of order $1/\log X$ keeps that factor bounded. Uniform $o(Q)$ primitive errors then give $o(X)$ after smoothing and Mellin inversion. The character argument continues to use a positive power saving.

\begin{lemma}[Moving Mellin lines]\label{wk:transfer}
Assume RH for the relevant $L$-functions and the weak logarithmic ratios Hypothesis~\ref{hyp:WZ} or~\ref{hyp:WH}. Put $c=4$ and
\begin{equation}\label{wk:line}
 \eta_X=\frac c{\log X}.
\end{equation}
Let $W\in C_c^\infty((0,\infty))$ be a fixed radial weight allowed in the corresponding detector. In the continuous case use the time weight $g_h$ of Section~\ref{sec:twins} and $T=2\pi X$; in the angular case use the weight $g$ of Lemma~\ref{leqc:window} and $K=\pi X/2$. Choose their fixed vanishing orders at zero sufficiently large. Then~\eqref{eq:twin-mainintegral} and~\eqref{eqc:main-integral}, respectively, hold on $\Re z=2\eta_X$ with error $o(X)$ in place of their displayed error terms.
\end{lemma}

\begin{proof}
\emph{Step 1: Bounds near the critical line.}
The local zero expansion, together with RH, gives the bounds
\begin{equation}\label{wk:ldbound}
 \left|\frac{L'}L(\tfrac12+\eta_X+it)\right|
 \ll_c (\log X)\log^2(3+Q_L+|t|),
\end{equation}
where $Q_L=1$ for $\zeta$ and $Q_L=k$ for the angular character of index $k$. The pole of $\zeta$ contributes a bounded term on this line. Indeed, the zeros within unit distance contribute at most their number, $O(\log(3+Q_L+|t|))$, times $\eta_X^{-1}$; the remaining local expansion is bounded by a logarithmic conductor term. We have allowed an extra logarithm in~\eqref{wk:ldbound}. Thus logarithmic derivatives have only a fixed power of $\log X$ as growth when the other parameters are bounded by fixed powers of $X$.

The differentiated arithmetic models enjoy the same type of bound, with polynomial growth in the Mellin heights. For the continuous kernel this follows by factoring
\begin{equation*}
 \prod_p F_p(z)=\zeta(1+z)\mathcal E(z),
 \qquad \mathcal E(z)=\prod_p(1-p^{-1-z})F_p(z).
\end{equation*}
The latter product is uniformly absolutely convergent on $0\le\Re z\le1/4$, since its local factors are $1+O(p^{-2})$ there. Consequently $P_1(z)$ and $C_1(z)$ on $\Re z=2\eta_X$ are bounded by a fixed power of $\log X$ times a fixed polynomial in $|\Im z|$. For the angular model, the products $A_1(-r,r)$ and $A_3(-r,r)$ are uniformly absolutely convergent down to $\Re r=0$; the zeta and Dirichlet factors in~\eqref{eqc:BC} give the same bound for $\mathcal B(r)$ and $\mathcal C(r)$ on $\Re r=\eta_X$.

\emph{Step 2: Truncation of the Mellin variables.}
For any fixed $W\in C_c^\infty((0,\infty))$ and every integer $A\ge0$, integration by parts gives
\begin{equation}\label{wk:Mellinbound}
 |\widetilde W(u+iv)|\ll_{W,A}(1+|v|)^{-A},
 \qquad 0\le u\le1/4.
\end{equation}
The constants are uniform in $u$. We can therefore truncate every Mellin integral at $|v|\le D\log X$, with any fixed $D>0$. Equations
\eqref{wk:ldbound} and~\eqref{wk:Mellinbound} show that the omitted
parts, after integration against the fixed time or family weight, are $o(X)$: choosing $A$ sufficiently large beats every fixed logarithmic loss. The same argument permits the eventual extension of the explicit model integrals back to complete vertical lines.

\emph{Step 3: Transfer of the qualitative error.}
We transfer the ratios error to the smooth time or family weight. Take $\kappa=1/2$ and $R=2$. On
\begin{equation}\label{wk:middle}
 X^\kappa\le Q\le X^R,\qquad |v_j|\le D\log X,
\end{equation}
the denominator shift $\eta_X$ belongs to $[2/\log Q,8/\log Q]$, and $|v_j|\le2D\log Q$. Choose the fixed $D>0$ small enough for the height window in Hypothesis~\ref{hyp:WH}. Hypothesis~\ref{hyp:WZ} permits every fixed logarithmic window. Their uniformity, followed by Cauchy's formula on the fixed numerator disks, therefore gives a differentiated error
\begin{equation}\label{wk:epsilon}
 E_X(Q;\boldsymbol v)=o(Q)
\end{equation}
uniformly on~\eqref{wk:middle}. No logarithmic differentiation loss occurs. The two main terms are interpreted together at removable singularities before differentiation.

Let $T\asymp X$ and let $g$ be a fixed Schwartz weight with a zero of order $J$ at the origin. In the continuous case, subtract the explicit primitive in~\eqref{eq:diff-zeta} from the actual primitive. In the angular case use the corresponding difference in
\eqref{eqc:differentiated}. Integration by parts, or Abel summation,
expresses the weighted error as
\begin{equation}\label{wk:Abelerror}
 -\frac1T\int_0^\infty g'(Q/T)E_X(Q;\boldsymbol v)\,dQ.
\end{equation}
On the middle range~\eqref{wk:middle}, its absolute value is at most
\begin{equation*}
 o(1)\,T\int_0^\infty v|g'(v)|\,dv=o(X),
\end{equation*}
uniformly in the retained Mellin variables. Their Mellin transforms have uniformly bounded $L^1$ norms by~\eqref{wk:Mellinbound}, so this estimate is preserved by the remaining integrations.

\emph{Step 4: The remaining family ranges.}
Outside the middle range, the ratios hypothesis is replaced by the direct bounds from Step~1. There the actual primitive and its explicit model are bounded by $Q$ times fixed logarithmic factors for $Q\ge1$, by
\eqref{wk:ldbound} and the model bounds above. For $0<Q<1$ the model
can instead have the harmless power $Q^{1-2\eta_X}$. Since $g'(v)\ll v^{J-1}$ at zero, the contribution below $X^\kappa$ to
\eqref{wk:Abelerror} is
\begin{equation*}
 \ll X^{-J+\kappa(J+1)}(\log X)^B+X^{-J}(\log X)^B.
\end{equation*}
Increasing the fixed vanishing order $J$ makes it negligible. Above $X^R$, Schwartz decay of $g'$ gives a negligible contribution by the same argument. Boundary terms vanish. These estimates also justify the integration by parts when the angular primitive is a step function.

Finally, the radial Mellin factors cost at most $X^{2\eta_X}=e^{2c}$. Thus the weighted error $o(X)$ obtained above remains $o(X)$ after Mellin inversion. This proves the lemma.
\end{proof}

\begin{proposition}[The two zeta applications]\label{wk:zeta}
Assume RH and Hypothesis~\ref{hyp:WZ}. For every fixed $h\ge1$ and $\Phi\in C_c^\infty((0,\infty))$,
\begin{equation}\label{wk:pairs}
 \sum_n\Lambda(n)\Lambda(n+h)\Phi(n/X)
 =\mathfrak S(h)X\int_0^\infty\Phi(u)\,du+o_{\Phi,h}(X).
\end{equation}
Moreover,
\begin{equation}\label{wk:squares}
 \pi((n+1)^2)-\pi(n^2)\sim\frac n{\log n},
 \qquad p_{n+1}-p_n=o(\sqrt{p_n}).
\end{equation}
\end{proposition}

\begin{proof}
For prime pairs, use the exact Fourier detector
\eqref{eq:exact-difference} and move the two Mellin lines to
$\Re a=\Re b=\eta_X$. RH justifies these shifts. The pole residue in
\eqref{eq:zeta-Mellin} is negligible against the high-order zero of the
time weight, exactly as in the fixed-line argument. Lemma~\ref{wk:transfer} now evaluates the moment with error $o(X)$. The diagonal term vanishes because $\widetilde g_h(1)=0$. The swapped term is still the same integral of $D_1(z)\Phi_h(z)$; its only pole in $\Re z>0$ is canceled by $\Phi_h(1)=0$. Moving the $z$-line to the right and inserting~\eqref{co:series} select the coefficient $\mathfrak S(h)$, independently of the location of the original line. The polarization and partition of unity following
\eqref{eq:twin-local} give~\eqref{wk:pairs} for every fixed $\Phi$.

For the interval assertion, apply Hypothesis~\ref{hyp:WZ} at $Q=H$
and $Q=2H$ with $a=\sigma_H$ and $b=\sigma_H+ir$ at both endpoints.
For $|r|\le A\log H$, subtracting the two formulas gives the
$o_A(H)$ error on the fixed numerator disks in
Proposition~\ref{si:weak}. This proves~\eqref{wk:squares}.
\end{proof}

\begin{proposition}[The angular application]\label{wk:Hecke}
Assume GRH for the angular Hecke family and Hypothesis~\ref{hyp:WH}. For every fixed $\Phi\in C_c^\infty((0,\infty))$,
\begin{equation}\label{wk:quadratic}
 \sum_{n\ge1}\Lambda(n^2+1)\Phi(n/X)
 =\mathfrak S_E X\int_0^\infty\Phi(u)\,du+o_\Phi(X).
\end{equation}
\end{proposition}

\begin{proof}
In the angular detector set $K=\pi X/2$ and move the Mellin line to $\Re z=2\eta_X$. Lemma~\ref{wk:transfer} replaces the weighted logarithmic-derivative sum by its explicit swapped term, with total error $o(X)$. The unswapped term vanishes since $\int_0^\infty g=0$. The identity~\eqref{eqc:FE}, the canceled pole at $z=1$, and the coefficient extraction in~\eqref{eqc:coefficient-extraction} are exact arithmetic calculations and are unchanged by the moving line. They select $S_1=\mathfrak S_E$. The exact detector and a finite partition of unity prove~\eqref{wk:quadratic}.
\end{proof}

Both~\eqref{wk:pairs} and~\eqref{wk:quadratic} yield their sharp weighted asymptotics without a rate of convergence in the ratios hypotheses. A fixed smooth majorant of $[1/2,1]$ first bounds each dyadic block by $O(X)$. Summing blocks bounds the part below $\epsilon X$ by $O(\epsilon X)+o(X)$. Positivity and fixed smooth upper and lower approximations on $[\epsilon,1]$ then allow $X\to\infty$ first and $\epsilon\to0$ afterward. Removing proper prime powers and partial summation give the unweighted prime counts. For Goldbach, the positive rational reference weight in Section~\ref{ug:section} supplies the analogous endpoint control directly from the ordinary character ratios formula.

The preceding propositions prove the qualitative counterparts of Theorems~\ref{thm:prime-pairs}--\ref{thm:quadratic-main}. The quantitative assertions require the stronger hypotheses stated in the introduction and are proved in Sections~\ref{sec:twins}--\ref{sec:quadratic}; the Goldbach saving is obtained separately in Section~\ref{ug:exact-section}.

{
For intervals of length $h\asymp\sqrt X$, a ratios asymptotic
for logarithmic shifts suffices. With $M_H=L\log H$ for fixed
$L\ge1$, the residue contributes $O(XH/L)$ to the local
mean-square bound for $B_H$, while an $o_L(H)$ ratios error
contributes $o_L(XH)$. We may therefore let $X\to\infty$ first
and then let $L\to\infty$.\par}

\begin{proposition}\label{si:weak}
Assume RH. Suppose that for every fixed $A>0$, the error
$\Delta_H(u,v;r)$ in~\eqref{si:RC} is $o_A(H)$ uniformly for
$|r|\le A\log H$ and $|u|,|v|\le\rho$. Then
{
\begin{equation}\label{si:Legendre-weak}
 \psi((n+1)^2)-\psi(n^2)\sim 2n,\qquad
 \pi((n+1)^2)-\pi(n^2)\sim\frac{n}{\log n}.
\end{equation}
}
and $p_{j+1}-p_j=o(\sqrt{p_j})$.
\end{proposition}

\begin{proof}
Fix $L\ge1$ and suppose $h\asymp\sqrt X$.
% {\hypersetup{linkcolor=red}
For $X/h\le t\le2X$, use the dyadic decomposition of $(0,t]$ in
the proof of Theorem~\ref{si:main} and retain the blocks with
$X^{1/100}\le H\le X$.
% }
On these blocks take $M_H=L\log H$. For sufficiently large $X$,
we have $\log H\le M_H\le H^{1/2}/4$, $2M_Hh\le V$, and
$H/4\le V\le H^{200}$ for $X\le V\le X+h$.
Apply Lemma~\ref{si:diff} with $\rho'=\rho$ and the hypothesis
with $A=2L$. This gives
\[
 \sup_{|r|\le2M_H}|E_H(r)|=o_L(H)
\]
uniformly on the retained blocks. After multiplication by
$V/M_H$ in~\eqref{si:localize}, the residue in
\eqref{si:testedformula} is $O(XH/L)$ with an absolute implied
constant, and the ratios error is $o_L(XH)$. For fixed
$0<c<1/2$ and $0<\varepsilon<1$, the three remaining terms,
divided by $XH$, are bounded by a constant times
\[
 X^{-c}M_H^\varepsilon
 +\left(\frac{H}{X}\right)^c M_H^{-1+\varepsilon}
 +\frac{M_H(\log H)^B}{H}.
\]
All three tend to zero uniformly for $X^{1/100}\le H\le X$.
We conclude that
\[
 \int_V^{V+h}|B_H(u)|^2\,du\ll \frac{XH}{L}+o_L(XH),
\]
where the implied constant on the first term is absolute.
The sum of $H^{1/2}$ over the dyadic blocks of $(0,t]$ is
$O(t^{1/2})$. { Including the unretained zeros
as in the proof of Theorem~\ref{si:main}, we obtain}
\[
 \begin{aligned}
 \|S_t\|_{L^2([V,V+h])}
 &\ll (L^{-1/2}+o_L(1))X^{1/2}
       \sum_{\substack{(H,2H]\ \text{retained}}}H^{1/2}
       +o(X^{1/2})\\
 &\ll (L^{-1/2}+o_L(1))(Xt)^{1/2}
 \qquad(X\le V\le X+h,\quad {X/h}\le t\le2X).
 \end{aligned}
\]
Squaring and applying Lemma~\ref{si:transfer} with $b=1$ gives
\[
 J(X,h)\ll\frac{Xh}{L}+o_L(Xh).
\]
The constant on the first term is independent of $L$.
Since $Xh\asymp h^3$, taking the upper limit as $X\to\infty$
and then letting $L\to\infty$ proves $J(X,h)=o(h^3)$,
uniformly when $h/\sqrt X$ lies in a fixed compact subinterval
of $(0,\infty)$. For each fixed $0<\delta<1/2$, both mean
squares in~\eqref{si:monotonicity} are then $o_\delta(h^3)$, so
\[
 \frac{|\psi(X+h)-\psi(X)-h|}{h}\le\delta+o_\delta(1).
\]
Letting $\delta\to0$ gives $\psi(X+h)-\psi(X)\sim h$.
{
The bound for proper prime powers and the comparison of $\log p$
with $\log X$ in the proof of Theorem~\ref{si:main} give}
$\pi(X+h)-\pi(X)\sim h/\log X$. Taking $X=n^2$ and $h=2n+1$
proves~{{\eqref{si:Legendre-weak}}}. For every fixed $\lambda>0$,
taking $X=p_j$ and $h=\lambda\sqrt{p_j}$ shows that
$p_{j+1}-p_j\le\lambda\sqrt{p_j}$ for all sufficiently large $j$.
This proves the gap assertion.
\end{proof}

\section{Reciprocal ratios and Chowla}
\label{rc:section}
The signed filter also applies to M\"obius and Liouville coefficients, provided a reciprocal-ratios prediction is available. This is an additional hypothesis, separate from the three inputs to Theorems~\ref{thm:goldbach}--\ref{thm:quadratic-main}. Its model has only a diagonal, which the filter removes. For Liouville, truncating the square-divisor expansion reduces the required uniformity to small square twists.

Put
\begin{equation*}
 D_\mu(z)=\frac{\zeta(1+z)}{\zeta(2+2z)}
          =\sum_{n\ge1}\frac{\mu(n)^2}{n^{1+z}},\qquad \Re z>0.
\end{equation*}
Retain $\sigma=1/2+\eta$ and the zeta zero-free margin from Section~\ref{sec:prelim}.
\begin{hypothesis}[Square-twisted reciprocal ratios]\label{rc:hyp}
There are fixed $0<\omega_C\le1/2$, $\theta_C,\nu>0$, and $A,B\ge0$ such that, for every $\varepsilon>0$, uniformly for $U\ge2$, $\Re a=\Re b=\eta$, $|\Im a|,|\Im b|\le U^{\theta_C}$, and integers $1\le k,l\le U^\nu$,
\begin{align}\label{rc:moment}
 \int_0^U\frac{(l/k)^{2it}}
 {\zeta(1/2+a+it)\zeta(1/2+b-it)}\,dt
 ={}&\ind_{k=l}\,U D_\mu(a+b)\\
 &+O_\varepsilon\bigl(U^{1-\omega_C+\varepsilon}(kl)^A
              (1+|\Im a|+|\Im b|)^B\bigr).\notag
\end{align}
\end{hypothesis}
The diagonal is the denominator-only prediction of the ratios recipe~\cite{CFZ,CS}.  Indeed its formal coefficient constraint is $k^2m=l^2n$.  If $k=gk_0$, $l=gl_0$, and $(k_0,l_0)=1$, then $m=l_0^2j$ and $n=k_0^2j$.  Unless $k=l$, one M\"obius coefficient vanishes; when $k=l$ the diagonal is $D_\mu(a+b)$. The stated square-twist range and error cost are additional uniformity assumptions, not consequences here of a published untwisted formula.

\begin{theorem}\label{rc:theorem}
Assume the zeta zero-free margin and Hypothesis~\ref{rc:hyp}, with
\begin{equation*}
 \Delta_C=\omega_C-2\eta>0,\qquad
 \rho_C=\min\left\{\frac\nu4,\frac{\Delta_C}{2(2A+1)}\right\}.
\end{equation*}
For each fixed $h\ge1$ and $\Phi\in C_c^\infty((0,\infty))$,
\begin{align}\label{rc:smooth}
 \sum_n\mu(n)\mu(n+h)\Phi(n/X)&\ll_{\Phi,h,\varepsilon}
                                  X^{1-\Delta_C+\varepsilon},\\
 \sum_n\lambda(n)\lambda(n+h)\Phi(n/X)&\ll_{\Phi,h,\varepsilon}
                                  X^{1-\rho_C+\varepsilon}.\notag
\end{align}
On each fixed compact support the constants depend polynomially on finitely many seminorms of $\Phi$.  For each $c\in\{\mu,\lambda\}$ there is $\kappa>0$ such that
\begin{equation}\label{rc:sharp}
 \sum_{n\le X}c(n)c(n+h)\ll_h X^{1-\kappa}.
\end{equation}
For the M\"obius conclusions only the untwisted case $k=l=1$ of Hypothesis~\ref{rc:hyp} is needed.
\end{theorem}

\begin{proof}
Use $g_h$, $W$, $P_{c,W,X}$, $T=2\pi X$, and $I_c(X)$ as in Section~\ref{sec:twins}; choose the vanishing order $2M$ of $g_h$ sufficiently large in terms of the fixed hypothesis parameters. Mellin inversion and the zero-free margin give
\begin{equation}\label{rc:mu-Mellin}
 P_{\mu,W,X}(t)=\frac1{2\pi i}\int_{(\eta)}
       \frac{X^a\widetilde W(a)}{\zeta(1/2+a+it)}\,da.
\end{equation}
There is no pole residue, since $1/\zeta$ vanishes at one.  The polynomial bound for $1/\zeta$ justifies the contour shift. Insert~\eqref{rc:moment} with $k=l=1$ into the two Mellin integrals and apply the weighted transfer~\eqref{eq:weighted-moment}. Its model vanishes because $\widetilde g_h(1)=\int_0^\infty g_h(v)\,dv=0$, by
\eqref{eq:two-cancellations}.  Therefore
\begin{equation*}
 I_\mu(X)\ll_{W,h,\varepsilon}X^{1-\Delta_C+\varepsilon}.
\end{equation*}
Equation~\eqref{eq:approx-difference}, followed by the same polarization and finite partition of unity used for prime pairs, proves the first line of~\eqref{rc:smooth}.

For Liouville use
\begin{equation*}
 \lambda(n)=\sum_{d^2\mid n}\mu(n/d^2),\qquad
 \lambda_D(n)=\sum_{\substack{d^2\mid n\\d\le D}}\mu(n/d^2).
\end{equation*}
There is exactly one $d$ for which $n/d^2$ is squarefree.  Hence $\lambda_D(n)$ is either $\lambda(n)$ or zero; in particular $|\lambda_D(n)|\le1$.  Moreover
\begin{equation*}
 \sum_{n\le Y}|\lambda(n)-\lambda_D(n)|
 \le\sum_{d>D}\left\lfloor\frac Y{d^2}\right\rfloor\ll\frac YD.
\end{equation*}
For each fixed smooth $\Phi$ this implies
\begin{equation}\label{rc:truncate}
 \sum_n\bigl(\lambda(n)\lambda(n+h)
       -\lambda_D(n)\lambda_D(n+h)\bigr)\Phi(n/X)
 \ll_{\Phi,h}X/D.
\end{equation}
Writing $Z_D(w)=\sum_{d\le D}d^{-w}$, the Mellin formula is
\begin{equation*}
 P_{\lambda_D,W,X}(t)=\frac1{2\pi i}\int_{(\eta)}
 \frac{Z_D(1+2a+2it)}{\zeta(1/2+a+it)}
          X^a\widetilde W(a)\,da.
\end{equation*}
On this line $|Z_D(1+2a+2it)|\le\zeta(1+2\eta)$ uniformly in $D$. Expanding the two finite numerators produces~\eqref{rc:moment}, with weights $k^{-1-2a}l^{-1-2b}$.  Every diagonal model again vanishes against $g_h$, while the ratios error contributes at most
\begin{equation}\label{rc:largeheight}
 X^{1-\Delta_C+\varepsilon}
 \left(\sum_{k\le D}k^{-1-2\eta+A}\right)^2
 \ll_{\eta}X^{1-\Delta_C+\varepsilon}D^{2A}.
\end{equation}

For completeness, the twist restriction does not require a uniform formula down to time zero.  Put $L=2+|\Im a|+|\Im b|$ and $U_0=C\max(2,L^{1/\theta_C},D^{1/\nu})$.  Above $U_0$ the hypothesis applies.  Below $U_0$, the integrand and the model have polynomial growth, uniformly in the real twists.  The integration-by-parts argument of Section~\ref{sec:weight} bounds this part, after the two Mellin integrals, by
\begin{equation*}
 O_W\left(X^{2\eta-2M}
              \{1+D^{(2M+C_0+1)/\nu}\}\right)
\end{equation*}
for a fixed $C_0$; powers of $L$ are absorbed by rapid Mellin decay. Take $D=\lfloor X^{\rho_C}\rfloor$.  Since $\rho_C\le\nu/4$, choosing $M$ large makes the last display $O_W(X^{-2})$. Thus~\eqref{rc:largeheight} bounds $I_{\lambda_D}$ up to this negligible error.  The Taylor error in~\eqref{eq:approx-difference} is uniform in $D$, because $|\lambda_D|\le1$.  Polarization, partition of unity, and~\eqref{rc:truncate} now give the Liouville estimate, since
\begin{equation*}
 \Delta_C-2A\rho_C\ge\rho_C+\Delta_C/2.
\end{equation*}
Only finitely many integrations by parts and a fixed vanishing order were used, proving the seminorm assertion.

Finally approximate the indicator of a fixed dyadic interval by a smooth cutoff with transition width $\epsilon$.  Since $|c|\le1$, the endpoint error is $O(\epsilon X+1)$.  If the smooth seminorm cost is $O(\epsilon^{-K})$, choose $\epsilon=X^{-a}$ with $a>0$ sufficiently small, and choose the $\varepsilon$ in
\eqref{rc:smooth} smaller than the saving.  This yields a power saving
on each dyadic interval.  Summing these intervals proves
\eqref{rc:sharp}.
\end{proof}

Partial summation of~\eqref{rc:sharp} shows that, for $c=\mu$ or $c=\lambda$, the series $\sum_{n\ge1}c(n)c(n+h)/n$ converges, with a power-saving tail. No mean value for $\zeta(2s)/\zeta(s)$ at mixed heights is assumed: the numerator needed for Liouville was truncated before the spectral average was taken.


\begingroup
\small
\begin{thebibliography}{MPPRW24}

\bibitem[ACD+26]{ACD26}
T. C. Anderson, S. Chow, R. Dietmann, D. Hokken, P. Koymans, R. J. Lemke Oliver, W. Sawin and M. Shusterman,
\emph{The strong form of Van der Waerden's conjecture via twisted Chowla},
\href{https://arxiv.org/abs/2608.22297}{arXiv:2608.22297v2}.

\bibitem[Aum25]{Aumonier}
A. Aumonier,
\emph{An $h$-principle for complements of discriminants},
Geom. Topol. \textbf{29} (2025), no. 3, 1441--1488.
\href{https://doi.org/10.2140/gt.2025.29.1441}{doi:10.2140/gt.2025.29.1441};
\href{https://arxiv.org/abs/2112.00326}{arXiv:2112.00326}.

\bibitem[BBR15]{BBR}
E. Bank, L. Bary-Soroker and L. Rosenzweig,
\emph{Prime polynomials in short intervals and in arithmetic progressions},
Duke Math. J. \textbf{164} (2015), 277--295.
\href{https://doi.org/10.1215/00127094-2856728}{doi:10.1215/00127094-2856728}.

\bibitem[Baz00]{Bazzanella}
D. Bazzanella,
\emph{Primes between consecutive squares},
Arch. Math. (Basel) \textbf{75} (2000), no. 1, 29--34.
\href{https://doi.org/10.1007/s000130050469}{doi:10.1007/s000130050469}.

\bibitem[BDPW23]{BDPW}
J. Bergstr\"om, A. Diaconu, D. Petersen and C. Westerland,
\emph{Hyperelliptic curves, the scanning map, and moments of families of quadratic $L$-functions},
\href{https://arxiv.org/abs/2302.07664}{arXiv:2302.07664}.

\bibitem[BFK23]{BFK}
H. M. Bui, A. Florea and J. P. Keating,
\emph{The Ratios Conjecture and upper bounds for negative moments of $L$-functions over function fields},
Trans. Amer. Math. Soc. \textbf{376} (2023), no. 6, 4453--4510.
\href{https://doi.org/10.1090/tran/8907}{doi:10.1090/tran/8907};
\href{https://arxiv.org/abs/2109.10396}{arXiv:2109.10396}.

\bibitem[BF24]{BFZeta}
H. M. Bui and A. Florea,
\emph{Negative moments of the Riemann zeta-function},
J. Reine Angew. Math. \textbf{806} (2024), 247--288.
\href{https://doi.org/10.1515/crelle-2023-0091}{doi:10.1515/crelle-2023-0091};
\href{https://arxiv.org/abs/2302.07226}{arXiv:2302.07226}.

\bibitem[BFH26]{BuiFloreaHeycock}
H. M. Bui, A. Florea and A. Heycock,
\emph{Small values of the Riemann zeta-function near the critical line},
\href{https://arxiv.org/abs/2610.27591}{arXiv:2610.27591}.

\bibitem[BFM24]{BFM}
H. M. Bui, A. Florea and M. B. Milinovich,
\emph{Negative discrete moments of the derivative of the Riemann zeta-function},
Bull. London Math. Soc. \textbf{56} (2024), no. 8, 2680--2703.
\href{https://doi.org/10.1112/blms.13092}{doi:10.1112/blms.13092};
\href{https://arxiv.org/abs/2310.03949}{arXiv:2310.03949}.

\bibitem[BGW26]{BGWSums}
T. D. Browning, J. Glas and V. Y. Wang,
\emph{Sums of three cubes over a function field},
Forum Math. Sigma \textbf{14} (2026), Paper No. e112, 58 pp.
\href{https://doi.org/10.1017/fms.2026.10259}{doi:10.1017/fms.2026.10259};
\href{https://arxiv.org/abs/2402.07146}{arXiv:2402.07146}.

% \bibitem[Caz12]{Cazanave}
% B. Cazanave, \emph{Algebraic homotopy classes of rational functions},
% Ann. Sci. \'Ec. Norm. Sup\'er. (4) \textbf{45} (2012), 511--534.
% \href{https://arxiv.org/abs/0912.2227}{arXiv:0912.2227}.

\bibitem[CFKRS05]{CFKRS}
J. B. Conrey, D. W. Farmer, J. P. Keating, M. O. Rubinstein and N. C. Snaith,
\emph{Integral moments of $L$-functions},
Proc. London Math. Soc. (3) \textbf{91} (2005), 33--104.
\href{https://doi.org/10.1112/S0024611504015175}{doi:10.1112/S0024611504015175};
\href{https://arxiv.org/abs/math/0206018}{arXiv:math/0206018}.

\bibitem[CFZ08]{CFZ}
J. B. Conrey, D. W. Farmer and M. R. Zirnbauer,
\emph{Autocorrelation of ratios of $L$-functions},
Commun. Number Theory Phys. \textbf{2} (2008), 593--636.
\href{https://arxiv.org/abs/0711.0718}{arXiv:0711.0718}.

\bibitem[CS07]{CS}
J. B. Conrey and N. C. Snaith,
\emph{Applications of the $L$-functions ratios conjectures},
Proc. Lond. Math. Soc. (3) \textbf{94} (2007), 594--646.
\href{https://doi.org/10.1112/plms/pdl021}{doi:10.1112/plms/pdl021}.

\bibitem[DH24]{DasHowe}
R. Das and S. Howe, \emph{Cohomological and motivic inclusion--exclusion}, Compos. Math. \textbf{160} (2024), 2228--2283.
\href{https://arxiv.org/abs/2204.04165}{arXiv:2204.04165}.

\bibitem[Del80]{Deligne}
P. Deligne, \emph{La conjecture de Weil. II}, Publ. Math. Inst. Hautes \'Etudes Sci. \textbf{52} (1980), 137--252.

\bibitem[DL76]{DeligneLusztig}
P. Deligne and G. Lusztig, \emph{Representations of reductive groups
over finite fields}, Ann. of Math. (2) \textbf{103} (1976), 103--161.

\bibitem[Flo24]{FloreaSmallShifts}
A. Florea,
\emph{Negative moments of $L$-functions with small shifts over function fields},
Int. Math. Res. Not. IMRN \textbf{2024}, no. 3, 2298--2337.
\href{https://arxiv.org/abs/2111.10477}{arXiv:2111.10477}.

\bibitem[GM87]{GM}
D. A. Goldston and H. L. Montgomery,
\emph{Pair correlation of zeros and primes in short intervals},
in \emph{Analytic Number Theory and Diophantine Problems}, Progr. Math. \textbf{70}, Birkh\"auser, Boston, 1987, 183--203.

\bibitem[GM25]{GrimmeltMerikoski}
L. Grimmelt and J. Merikoski,
\emph{On the greatest prime factor and uniform equidistribution of quadratic polynomials},
\href{https://arxiv.org/abs/2505.00493}{arXiv:2505.00493}.

\bibitem[Har23]{HarperBTSC}
A. J. Harper,
\emph{The typical size of character and zeta sums is $o(\sqrt{x})$},
\href{https://arxiv.org/abs/2301.04390}{arXiv:2301.04390}.

\bibitem[Har26]{HarperICM}
A. J. Harper,
\emph{Better than squareroot cancellation in number theory},
\href{https://arxiv.org/abs/2512.23681}{arXiv:2512.23681};
to appear in the \emph{Proceedings of the ICM 2026}.

\bibitem[HB82]{HB}
D. R. Heath-Brown, \emph{Gaps between primes, and the pair correlation of zeros of the zeta-function}, Acta Arith. \textbf{41} (1982), 85--99.

\bibitem[Hel15]{HelfgottTernary}
H. A. Helfgott,
\emph{The ternary Goldbach problem},
\href{https://arxiv.org/abs/1501.05438}{arXiv:1501.05438}, 2015, 317 pp.

\bibitem[IK04]{IK}
H. Iwaniec and E. Kowalski, \emph{Analytic Number Theory}, Amer. Math. Soc. Colloq. Publ. \textbf{53}, American Mathematical Society, Providence, RI, 2004.

\bibitem[JT23]{JT}
O. J\"arviniemi and J. Ter\"av\"ainen,
\emph{Composite values of shifted exponentials},
Adv. Math. \textbf{429} (2023), Paper No. 109187.
\href{https://doi.org/10.1016/j.aim.2023.109187}{doi:10.1016/j.aim.2023.109187};
\href{https://arxiv.org/abs/2010.01789}{arXiv:2010.01789}.


\bibitem[Kat01]{Katz}
N. M. Katz,
\emph{Sums of Betti numbers in arbitrary characteristic},
Finite Fields Appl. \textbf{7} (2001), no. 1, 29--44.
\href{https://doi.org/10.1006/ffta.2000.0303}{doi:10.1006/ffta.2000.0303}.


\bibitem[Lee26]{Lee}
S. Lee, \emph{The Hardy--Littlewood asymptotic for prime pairs of fixed even difference from a uniform zeta-ratios conjecture}, preprint, revised 3 November 2026.
\href{https://arxiv.org/abs/2610.32672v2}{arXiv:2610.32672v2}.

\bibitem[LLOP26]{LLOP}
M. Lal\'in, K.-H. Lee, T. Oliver and A. Pozdnyakov,
\emph{Murmurations of quadratic and cubic characters over function fields},
\href{https://arxiv.org/abs/2608.01337}{arXiv:2608.01337}.

\bibitem[LPZ16]{LPZ}
A. Languasco, A. Perelli and A. Zaccagnini,
\emph{An extension of the pair-correlation conjecture and applications},
Math. Res. Lett. \textbf{23} (2016), 201--220.
\href{https://arxiv.org/abs/1308.3934}{arXiv:1308.3934}.

\bibitem[May19]{MaynardTakagi}
J. Maynard, \emph{The twin prime conjecture},
Jpn. J. Math. \textbf{14} (2019), 175--206.
\href{https://doi.org/10.1007/s11537-019-1837-z}{doi:10.1007/s11537-019-1837-z}.

\bibitem[MPPRW24]{MPPRW}
J. Miller, P. Patzt, D. Petersen and O. Randal-Williams,
\emph{Uniform twisted homological stability},
\href{https://arxiv.org/abs/2402.00354}{arXiv:2402.00354}.

% \bibitem[Mil13]{Milne}
% J. S. Milne, \emph{Lectures on \'Etale Cohomology}, version 2.21, 2013.
% \url{https://www.jmilne.org/math/CourseNotes/LEC.pdf}.

\bibitem[Mon73]{Montgomery}
H. L. Montgomery,
\emph{The pair correlation of zeros of the zeta function},
in \emph{Analytic Number Theory} (Proc. Sympos. Pure Math., Vol. XXIV,
St. Louis Univ., St. Louis, Mo., 1972),
Amer. Math. Soc., Providence, RI, 1973, 181--193.

\bibitem[Leb50]{Lebesgue1850}
V.-A. Lebesgue, \emph{Sur l'impossibilit\'e, en nombres entiers, de l'\'equation $x^m=y^2+1$},
Nouv. Ann. Math. (1) \textbf{9} (1850), 178--181.
\url{https://www.numdam.org/item/NAM_1850_1_9__178_0/}.

\bibitem[MS02]{MontgomerySoundararajan}
H. L. Montgomery and K. Soundararajan,
\emph{Beyond pair correlation},
in \emph{Paul Erd\H{o}s and his mathematics, I} (Budapest, 1999),
Bolyai Soc. Math. Stud. \textbf{11}, J\'anos Bolyai Math. Soc., Budapest, 2002, 507--514.

\bibitem[Pin09]{Pintz2009}
J. Pintz, \emph{Landau's problems on primes},
J. Th\'eor. Nombres Bordeaux \textbf{21} (2009), 357--404.
\href{https://doi.org/10.5802/jtnb.676}{doi:10.5802/jtnb.676}.


\bibitem[Poo03]{Poonen}
B. Poonen,
\emph{Squarefree values of multivariable polynomials},
Duke Math. J. \textbf{118} (2003), no. 2, 353--373.
\href{https://doi.org/10.1215/S0012-7094-03-11826-8}{doi:10.1215/S0012-7094-03-11826-8}.


\bibitem[Ros02]{Rosen}
M. Rosen, \emph{Number Theory in Function Fields}, Grad. Texts in Math. \textbf{210}, Springer, New York, 2002.

\bibitem[Saw24]{Sawin}
W. Sawin, \emph{Square-root cancellation for sums of factorization functions over squarefree progressions in $\mathbb F_q[t]$}, Acta Math. \textbf{233} (2024), 285--418.
\href{https://arxiv.org/abs/2102.09730}{arXiv:2102.09730}.

\bibitem[Saw21]{SawinIntervals}
W. Sawin, \emph{Square-root cancellation for sums of factorization functions over short intervals in function fields}, Duke Math. J. \textbf{170} (2021), 997--1026.
\href{https://arxiv.org/abs/1809.05137}{arXiv:1809.05137}.

\bibitem[SS22a]{SSAnnals}
W. Sawin and M. Shusterman,
\emph{On the Chowla and twin primes conjectures over $\mathbb F_q[T]$},
Ann. of Math. (2) \textbf{196} (2022), 457--506.
\href{https://doi.org/10.4007/annals.2022.196.2.1}{doi:10.4007/annals.2022.196.2.1}.

\bibitem[SS22b]{SSQuadratic}
W. Sawin and M. Shusterman,
\emph{M\"obius cancellation on polynomial sequences and the quadratic Bateman--Horn conjecture over function fields},
Invent. Math. \textbf{229} (2022), 751--927.
\href{https://doi.org/10.1007/s00222-022-01115-y}{doi:10.1007/s00222-022-01115-y}.

\bibitem[SS25]{SSBetti}
W. Sawin and M. Shusterman, \emph{Betti number bounds for the moduli space of maps from a rational curve},
Appendix B to R. Das, B. Lehmann, S. Tanimoto and P. Tosteson,
\emph{Homological stability and Manin's conjecture for rational curves on quartic del Pezzo surfaces},
\href{https://arxiv.org/abs/2506.17071v1}{arXiv:2506.17071v1}, 79--84.

\bibitem[TB26]{TurnageButterbaugh}
C. L. Turnage-Butterbaugh,
\emph{A decades-long breakthrough in zero-density estimates and primes in short intervals},
\href{https://arxiv.org/abs/2607.04632}{arXiv:2607.04632};
to appear in Bull. Amer. Math. Soc.

\bibitem[WZ26]{WanZhang}
D. Wan and D. Zhang,
\emph{Betti number bounds for varieties and exponential sums},
Adv. Math. \textbf{490} (2026), article 110852.
\href{https://doi.org/10.1016/j.aim.2026.110852}{doi:10.1016/j.aim.2026.110852};
\href{https://arxiv.org/abs/2501.12623v1}{arXiv:2501.12623v1}.

\bibitem[Wan24]{WangRatios}
V. Y. Wang, \emph{Notes on zeta ratio stabilization},
\href{https://arxiv.org/abs/2402.01214}{arXiv:2402.01214}.

\bibitem[WX25]{WX}
V. Y. Wang and M. W. Xu, \emph{Harper's beyond square-root conjecture}, Int. Math. Res. Not. IMRN \textbf{2025} (2025), no. 18, rnaf279.
\href{https://arxiv.org/abs/2405.04094}{arXiv:2405.04094}.

\bibitem[Wax21]{Waxman}
E. Waxman, \emph{Lower order terms for the one-level density of a symplectic family of Hecke $L$-functions}, J. Number Theory \textbf{221} (2021), 447--483.
\href{https://arxiv.org/abs/1905.10362}{arXiv:1905.10362}.

\end{thebibliography}
\endgroup
\end{document}